\documentclass[11pt]{amsart}
\usepackage{mathrsfs,stmaryrd}
\newtheorem{theorem}{Theorem}[section]
\newtheorem{proposition}[theorem]{Proposition} 
\newtheorem{lemma}[theorem]{Lemma} 
\newtheorem{corollary}[theorem]{Corollary} 
\theoremstyle{definition}
\newtheorem{definition}[theorem]{Definition} 
\newtheorem{remark}[theorem]{Remark} 

\begin{document}

\title{The rigidity statement in the Horowitz-Myers conjecture}
\author{Simon Brendle and Pei-Ken Hung}
\address{Columbia University \\ 2990 Broadway \\ New York NY 10027 \\ USA}
\address{University of Illinois \\ 1409 W. Green Street \\ Urbana IL 61801 \\ USA}
\thanks{The authors are grateful to Professors Piotr Chru\'sciel, Michael Eichmair, and Gerhard Huisken for discussions. The first author was supported by the National Science Foundation under grant DMS-2403981 and by the Simons Foundation. He acknowledges the hospitality of T\"ubingen University, where part of this work was carried out.}
\maketitle
\begin{abstract}
In this paper, we give an alternative proof of the Horowitz-Myers conjecture in dimension $3 \leq N \leq 7$. Moreover, we show that a metric that achieves equality in the Horowitz-Myers conjecture is locally isometric to a Horowitz-Myers metric. 
\end{abstract}

\section{Introduction}

In the 1990s, Horowitz and Myers \cite{Horowitz-Myers} proposed a new positive energy theorem for certain asymptotically locally hyperbolic manifolds of dimension $N \geq 3$ with scalar curvature at least $-N(N-1)$. The positive energy theorem conjectured in \cite{Horowitz-Myers} is sharp. Equality holds for the so-called Horowitz-Myers metrics. These are static metrics with scalar curvature $-N(N-1)$ that can be written as a multiply warped product over a halfline; see e.g. \cite{Woolgar}. 

The Horowitz-Myers conjecture has been studied by various authors; see e.g. \cite{Barzegar-Chrusciel-Horzinger-Maliborski-Nguyen}, \cite{Chrusciel-Delay-Wutte}, \cite{Chrusciel-Galloway-Nguyen-Paetz}, \cite{Galloway-Woolgar}, \cite{Liang-Zhang}, and \cite{Woolgar}. In particular, Barzegar, Chru\'sciel, H\"orzinger, Maliborski, and Nguyen \cite{Barzegar-Chrusciel-Horzinger-Maliborski-Nguyen} confirmed the conjecture for multiply warped products over a halfline. 

In a recent paper \cite{Brendle-Hung}, we proved the Horowitz-Myers conjecture in dimension $3 \leq N \leq 7$. The proof in \cite{Brendle-Hung} is based on a new geometric inequality for compact manifolds with boundary. Given a compact manifold with scalar curvature at least $-N(N-1)$, this inequality relates the boundary mean curvature to the systole of the boundary. The Horowitz-Myers conjecture follows from this inequality by a limiting process. 

In this paper, we develop an alternative approach to the Horowitz-Myers conjecture in dimension $3 \leq N \leq 7$. This alternative approach requires stronger asymptotic assumptions than the one in \cite{Brendle-Hung}. It does not give a new proof of the systolic inequality, but it does allow us to prove a rigidity statement. 

While the approach in our earlier paper \cite{Brendle-Hung} used slicings by free boundary minimal hypersurfaces, the arguments in this paper hew more closely to the original dimension descent scheme of Schoen and Yau (see \cite{Schoen}, Section 4, and \cite{Eichmair-Huang-Lee-Schoen}). 

In order to state the main result of this paper, we need the following definition. 

\begin{definition}
Let $N \geq 3$ be an integer. We define a metric $g_{\text{\rm HM},N}$ on $(2^{-\frac{2}{N}},\infty) \times S^1 \times \mathbb{R}^{N-2}$ by 
\begin{align*}
g_{\text{\rm HM},N} &= r^{-2} \, dr \otimes dr + \frac{4}{N^2} \, r^2 \, (1 + \frac{1}{4} \, r^{-N})^{\frac{4}{N}-2} \, (1 - \frac{1}{4} \, r^{-N})^2 \, d\tau_0 \otimes d\tau_0 \\ 
&+ r^2 \, (1 + \frac{1}{4} \, r^{-N})^{\frac{4}{N}} \, \sum_{k=1}^{N-2} d\tau_k \otimes d\tau_k, 
\end{align*}
where $r \in (2^{-\frac{2}{N}},\infty)$, $\tau_0 \in S^1 = \mathbb{R}/(2\pi \mathbb{Z})$, and $(\tau_1,\hdots,\tau_{N-2}) \in \mathbb{R}^{N-2}$. The metric $g_{\text{\rm HM},N}$ extends to a smooth metric on $\mathbb{R}^2 \times \mathbb{R}^{N-2}$. The resulting metric on $\mathbb{R}^2 \times \mathbb{R}^{N-2}$ is complete and has scalar curvature $-N(N-1)$ (see \cite{Woolgar}). 
\end{definition}

\begin{theorem} 
\label{main.theorem}
Let us fix an integer $N$ with $3 \leq N \leq 7$ and a collection of positive real numbers $b_0,\hdots,b_{N-2}$. Let $\theta_0,\hdots,\theta_{N-2}$ denote the coordinate functions on $T^{N-1}$, which take values in $S^1 = \mathbb{R}/(2\pi \mathbb{Z})$. We define a flat metric $\gamma$ on $T^{N-1}$ by $\gamma = \sum_{k=0}^{N-2} b_k^2 \, d\theta_k \otimes d\theta_k$. Given a positive real number $r_0$, we define a hyperbolic metric $\bar{g}$ on $(r_0,\infty) \times T^{N-1}$ by $\bar{g} = r^{-2} \, dr \otimes dr + r^2 \, \gamma$. Let $(M,g)$ be a complete, connected, orientable Riemannian manifold of dimension $N$ with the following properties: 
\begin{itemize} 
\item There exists a bounded open domain $E \subset M$ with smooth boundary such that the complement $M \setminus E$ is diffeomorphic to $[r_0,\infty) \times T^{N-1}$.
\item The map 
\[(r_0,\infty) \times T^{N-1} \to T^{N-2}, \quad (r,\theta_0,\hdots,\theta_{N-2}) \mapsto (\theta_1,\hdots,\theta_{N-2})\] 
extends to a globally defined smooth map from $M$ to $T^{N-2}$.
\item For every nonnegative integer $m$, the metric $g$ satisfies 
\[|\bar{D}^m (g-\bar{g})|_{\bar{g}} \leq O(r^{-N}),\] 
where $\bar{D}^m$ denotes the covariant derivative of order $m$ with respect to the hyperbolic metric $\bar{g}$. 
\item The metric $g$ satisfies 
\[|g-\bar{g} - r^{2-N} \, Q|_{\bar{g}} \leq O(r^{-N-2\delta}).\] 
Here, $\delta$ is a small positive real number and $Q$ is a smooth symmetric $(0,2)$-tensor on $T^{N-1}$.
\item We have 
\[\int_{T^{N-1}} \Big ( N \, \text{\rm tr}_\gamma(Q) + \Big ( \frac{2}{Nb_0} \Big )^N \Big ) \, d\text{\rm vol}_\gamma \leq 0.\] 
\end{itemize}
If the scalar curvature of $(M,g)$ is at least $-N(N-1)$, then there exists a smooth map $\Phi: \mathbb{R}^2 \times \mathbb{R}^{N-2} \to M$ such that $\Phi^* g = g_{\text{\rm HM},N}$.
\end{theorem}

\begin{remark} 
A local isometry between two complete Riemannian manifolds of the same dimension is a covering map and in particular is surjective (see \cite{Cheeger-Ebin}, Lemma 1.38).
\end{remark}

The case of equality in the Horowitz-Myers conjecture was studied by Woolgar \cite{Woolgar}. He showed that a metric that achieves equality in the Horowitz-Myers conjecture is static. The classification of static metrics is a challenging open problem.

We refer to \cite{Andersson-Cai-Galloway}, \cite{Andersson-Dahl}, \cite{Chrusciel-Galloway}, \cite{Lee-Neves}, \cite{Min-Oo}, and \cite{Wang} for other rigidity results for asymptotically hyperbolic manifolds.

The proof of Theorem \ref{main.theorem} is based on an inductive scheme. In the following, we give an overview of the main ideas. Let us fix an integer $N$ with $3 \leq N \leq 7$ and a collection of positive real numbers $b_0,\hdots,b_{N-2}$. 

In Section \ref{definitions}, we introduce the notion of an $(N,n)$-dataset, where $n$ is an integer with $2 \leq n \leq N$. An $(N,n)$-dataset consists of a complete, connected, orientable manifold $M$ of dimension $n$ together with a Riemannian metric $g$ and a positive smooth function $\rho$ satisfying certain conditions (see Definition \ref{definition.dataset}). In particular, we require that there exists a bounded open domain $E \subset M$ with smooth boundary such that the complement $M \setminus E$ is diffeomorphic to $[r_0,\infty) \times T^{n-1}$. We define a flat metric $\gamma$ on $T^{n-1}$ by $\gamma = \sum_{k=0}^{n-2} b_k^2 \, d\theta_k \otimes d\theta_k$. We define a hyperbolic metric $\bar{g}$ on $(r_0,\infty) \times T^{n-1}$ by $\bar{g} = r^{-2} \, dr \otimes dr + r^2 \, \gamma$. With this understood, we require that the metric $g$ satisfies 
\[|g-\bar{g} - r^{2-N} \, Q|_{\bar{g}} \leq O(r^{-N-2\delta}),\] 
where $\delta$ is a small positive real number and $Q$ is a smooth symmetric $(0,2)$-tensor on $T^{n-1}$. Moreover, we require that the function $\rho$ satisfies 
\[|\rho - r^{N-n} - r^{-n} \, P| \leq O(r^{-n-2\delta}),\] 
where $\delta$ is a small positive number and $P$ is a real-valued function on $T^{n-1}$ which is H\"older continuous with exponent $2\delta$. Finally, we require that 
\[\int_{T^{n-1}} \Big ( N \, \text{\rm tr}_\gamma(Q) + 2N \, P + \Big ( \frac{2}{Nb_0} \Big )^N \Big ) \, d\text{\rm vol}_\gamma \leq 0.\] 
For each integer $n$ with $2 \leq n \leq N$, we introduce a rigidity property for $(N,n)$-datasets, which we denote by $(\star_{N,n})$ (see Definition \ref{star.N.n}). This condition forms the basis of our inductive scheme. 

In Section \ref{2D}, we show that the condition $(\star_{N,2})$ holds. 

In Section \ref{properties.of.stable.hypersurfaces}, we consider an $(N,n)$-dataset $(M,g,\rho)$, where $n$ is an integer with $3 \leq n \leq N$. We fix a function $u: T^{n-1} \to \mathbb{R}$ such that 
\[\Delta_\gamma u + \frac{N}{2} \, \text{\rm tr}_\gamma(Q) + N \, P = \text{\rm constant}.\] 
We assume that $\Sigma$ is a properly embedded, connected, orientable hypersurface in $M$ satisfying certain asymptotic conditions near infinity (see Definition \ref{definition.tame}). We further assume that $\Sigma$ is $(g,\rho)$-stationary in the sense of Definition \ref{definition.stationary}; that is, $H_\Sigma + \rho^{-1} \, \langle \nabla \rho,\nu_\Sigma \rangle = 0$ at each point on $\Sigma$. Finally, we assume that $\Sigma$ is $(g,\rho,u)$-stable in the sense of Definition \ref{definition.stability}. Under these assumptions, we construct a positive smooth function $\check{v}: \Sigma \to \mathbb{R}$ such that $\mathbb{L}_\Sigma \check{v} \geq 0$, where $\mathbb{L}_\Sigma$ denotes the weighted Jacobi operator of $\Sigma$ (see Definition \ref{definition.weighted.Jacobi.operator}). Moreover, we show that $(\Sigma,\check{g},\check{\rho})$ is an $(N,n-1)$-dataset, where $\check{g}$ denotes the induced metric on $\Sigma$ and the weight function $\check{\rho}: \Sigma \to \mathbb{R}$ is defined to be a constant multiple of $\check{v} \, \rho|_\Sigma$.

In Section \ref{foliation}, we consider an asymptotically locally hyperbolic manifold $(M,g)$ of dimension $n$, where $3 \leq n \leq N$. We define an exponential map from infinity and use it to construct a foliation near infinity. This foliation will play an important role in Section \ref{proof.of.main.theorem}.

In Section \ref{barrier}, we consider an $(N,n)$-dataset $(M,g,\rho)$, where $n$ is an integer with $3 \leq n \leq N$. We construct barriers for $(g,\rho)$-stationary hypersurfaces. These barriers play a crucial role in the existence theory in Section \ref{existence.theory}.

In Section \ref{existence.theory}, we again consider an $(N,n)$-dataset $(M,g,\rho)$, where $n$ is an integer with $3 \leq n \leq N$. Moreover, we fix an arbitrary point $p_* \in M$. Under the assumption that condition $(\star_{N,n-1})$ holds, we construct a hypersurface $\Sigma$ passing through $p_*$ with the property that $\Sigma$ is $(g,\rho)$-stationary in the sense of Definition \ref{definition.stationary} and $(g,\rho,u)$-stable in the sense of Definition \ref{definition.stability}. Moreover, we show that $\Sigma$ satisfies suitable asymptotic estimates near infinity. To prove the existence of $\Sigma$, we use the barriers from Section \ref{barrier}. We also use ideas from Gang Liu's work \cite{Liu} in an important way (see also \cite{Carlotto-Chodosh-Eichmair} and \cite{Chodosh-Eichmair-Moraru} for related work).

Finally, in Section \ref{proof.of.main.theorem}, we show that the condition $(\star_{N,n})$ holds for all $2 \leq n \leq N$. The proof is by induction on $n$, and uses the results established in the previous sections. Theorem \ref{main.theorem} follows from the condition $(\star_{N,n})$ for $n=N$.

\section{Definitions and preliminary results}

\label{definitions}

Throughout this paper, we fix an integer $N$ with $N \geq 3$ and a collection of positive real numbers $b_0,\hdots,b_{N-2}$. 

\begin{definition}
Let us fix an integer $n$ with $2 \leq n \leq N$. We define a metric $g_{\text{\rm HM},N,n}$ on $(2^{-\frac{2}{N}},\infty) \times S^1 \times \mathbb{R}^{n-2}$ by 
\begin{align*}
g_{\text{\rm HM},N,n} &= r^{-2} \, dr \otimes dr + \frac{4}{N^2} \, r^2 \, (1 + \frac{1}{4} \, r^{-N})^{\frac{4}{N}-2} \, (1 - \frac{1}{4} \, r^{-N})^2 \, d\tau_0 \otimes d\tau_0 \\ 
&+ r^2 \, (1 + \frac{1}{4} \, r^{-N})^{\frac{4}{N}} \, \sum_{k=1}^{n-2} d\tau_k \otimes d\tau_k, 
\end{align*}
where $r \in (2^{-\frac{2}{N}},\infty)$, $\tau_0 \in S^1 = \mathbb{R}/(2\pi \mathbb{Z})$, and $(\tau_1,\hdots,\tau_{n-2}) \in \mathbb{R}^{n-2}$. We define a function $\Upsilon_{\text{\rm HM},N,n}: (2^{-\frac{2}{N}},\infty) \times S^1 \times \mathbb{R}^{n-2} \to [1,\infty)$ by 
\[\Upsilon_{\text{\rm HM},N,n} = r \, (1+\frac{1}{4} \, r^{-N})^{\frac{2}{N}}.\] 
Moreover, we define a function $\rho_{\text{\rm HM},N,n}: (2^{-\frac{2}{N}},\infty) \times S^1 \times \mathbb{R}^{n-2} \to (0,\infty)$ by $\rho_{\text{\rm HM},N,n} = \Upsilon_{\text{\rm HM},N,n}^{N-n}$. The metric $g_{\text{\rm HM},N,n}$ extends to a smooth metric on $\mathbb{R}^2 \times \mathbb{R}^{n-2}$, and the resulting metric on $\mathbb{R}^2 \times \mathbb{R}^{n-2}$ is complete. Finally, the functions $\Upsilon_{\text{\rm HM},N,n}$ and $\rho_{\text{\rm HM},N,n}$ extend to smooth functions on $\mathbb{R}^2 \times \mathbb{R}^{n-2}$. 
\end{definition}

\begin{proposition}
\label{properties.of.HM.metrics} 
Let us fix an integer $n$ with $2 \leq n \leq N$. We define a symmetric $(0,2)$-tensor $T$ on $(2^{-\frac{2}{N}},\infty) \times S^1 \times \mathbb{R}^{n-2}$ by 
\[T = r^{-2} \, dr \otimes dr + \frac{4}{N^2} \, r^2 \, (1 + \frac{1}{4} \, r^{-N})^{\frac{4}{N}-2} \, (1 - \frac{1}{4} \, r^{-N})^2 \, d\tau_0 \otimes d\tau_0.\] 
Then the following statements hold: 
\begin{itemize}
\item The eigenvalues of $T$ with respect to the metric $g_{\text{\rm HM},N,n}$ are $1$ and $0$, and the corresponding multiplicities are $2$ and $n-2$, respectively.
\item The Hessian of $\Upsilon_{\text{\rm HM},N,n}$ with respect to the metric $g_{\text{\rm HM},N,n}$ is given by 
\[\Upsilon_{\text{\rm HM},N,n} \, (1-\Upsilon_{\text{\rm HM},N,n}^{-N}) \, g_{\text{\rm HM},N,n} + \frac{N}{2} \, \Upsilon_{\text{\rm HM},N,n}^{1-N} \, T.\] 
\item The Riemann curvature tensor of $g_{\text{\rm HM},N,n}$ is given by 
\[-\frac{1}{2} \, (1-\Upsilon_{\text{\rm HM},N,n}^{-N}) \, g \owedge g - \frac{N}{2} \, \Upsilon_{\text{\rm HM},N,n}^{-N} \, T \owedge g + \frac{N(N-1)}{4} \, \Upsilon_{\text{\rm HM},N,n}^{-N} \, T \owedge T,\] 
where $\owedge$ denotes the Kulkarni-Nomizu product (see \cite{Besse}, Definition 1.110). 
\end{itemize}
\end{proposition}

\textbf{Proof.}
The metric $g_{\text{\rm HM},N,n}$ can be written in the form 
\begin{align*} 
g_{\text{\rm HM},N,n} 
&= \Upsilon_{\text{\rm HM},N,n}^{-2} \, (1-\Upsilon_{\text{\rm HM},N,n}^{-N})^{-1} \, d\Upsilon_{\text{\rm HM},N,n} \otimes d\Upsilon_{\text{\rm HM},N,n} \\ 
&+ \frac{4}{N^2} \, \Upsilon_{\text{\rm HM},N,n}^2 \, (1-\Upsilon_{\text{\rm HM},N,n}^{-N}) \, d\tau_0 \otimes d\tau_0 + \Upsilon_{\text{\rm HM},N,n}^2 \, \sum_{k=1}^{n-2} d\tau_k \otimes d\tau_k, 
\end{align*} 
and the tensor $T$ can be written in the form 
\begin{align*} 
T 
&= \Upsilon_{\text{\rm HM},N,n}^{-2} \, (1-\Upsilon_{\text{\rm HM},N,n}^{-N})^{-1} \, d\Upsilon_{\text{\rm HM},N,n} \otimes d\Upsilon_{\text{\rm HM},N,n} \\ 
&+ \frac{4}{N^2} \, \Upsilon_{\text{\rm HM},N,n}^2 \, (1-\Upsilon_{\text{\rm HM},N,n}^{-N}) \, d\tau_0 \otimes d\tau_0. 
\end{align*}
The assertion now follows from a straightforward calculation. \\

\begin{proposition}
\label{scalar.curvature.of.HM.metrics}
Let $n$ be an integer with $2 \leq n < N$. Then 
\begin{align*} 
&-2 \, \Delta_{g_{\text{\rm HM},N,n}} \log \rho_{\text{\rm HM},N,n} - \frac{N-n+1}{N-n} \, |d\log \rho_{\text{\rm HM},N,n}|_{g_{\text{\rm HM},N,n}}^2 \\ 
&+ R_{g_{\text{\rm HM},N,n}} + N(N-1) = 0. 
\end{align*}
\end{proposition}

\textbf{Proof.} 
By Proposition \ref{properties.of.HM.metrics}, the scalar curvature of $g_{\text{\rm HM},N,n}$ is given by 
\[R_{g_{\text{\rm HM},N,n}} = -n(n-1) + (N-n+1)(N-n) \, \Upsilon_{\text{\rm HM},N,n}^{-N}.\] 
Moreover, Proposition \ref{properties.of.HM.metrics} implies that 
\[\Upsilon_{\text{\rm HM},N,n}^{-1} \, \Delta_{g_{\text{\rm HM},N,n}} \Upsilon_{\text{\rm HM},N,n} = n + (N-n) \, \Upsilon_{\text{\rm HM},N,n}^{-N}.\] 
Finally, 
\[\Upsilon_{\text{\rm HM},N,n}^{-2} \, |d\Upsilon_{\text{\rm HM},N,n}|_{g_{\text{\rm HM},N,n}}^2 = 1-\Upsilon_{\text{\rm HM},N,n}^{-N}.\] 
Since $\rho_{\text{\rm HM},N,n} = \Upsilon_{\text{\rm HM},N,n}^{N-n}$, we conclude that 
\begin{align*} 
&-2 \, \Delta_{g_{\text{\rm HM},N,n}} \log \rho_{\text{\rm HM},N,n} - \frac{N-n+1}{N-n} \, |d\log \rho_{\text{\rm HM},N,n}|_{g_{\text{\rm HM},N,n}}^2 \\ 
&= -2(N-n) \, \Upsilon_{\text{\rm HM},N,n}^{-1} \, \Delta_{g_{\text{\rm HM},N,n}} \Upsilon_{\text{\rm HM},N,n} \\ 
&- (N-n-1)(N-n) \, \Upsilon_{\text{\rm HM},N,n}^{-2} \, |d\Upsilon_{\text{\rm HM},N,n}|_{g_{\text{\rm HM},N,n}}^2 \\ 
&= -(N+n-1)(N-n) - (N-n+1)(N-n) \, \Upsilon_{\text{\rm HM},N,n}^{-N}. 
\end{align*}
Putting these facts together, the assertion follows. This completes the proof of Proposition \ref{scalar.curvature.of.HM.metrics}. \\

\begin{definition} 
\label{definition.dataset}
Let $n$ be an integer with $2 \leq n \leq N$. Let $\theta_0,\hdots,\theta_{n-2}$ denote the coordinate functions on $T^{n-1}$, which take values in $S^1 = \mathbb{R}/(2\pi \mathbb{Z})$. We define a flat metric $\gamma$ on $T^{n-1}$ by $\gamma = \sum_{k=0}^{n-2} b_k^2 \, d\theta_k \otimes d\theta_k$. Given a positive real number $r_0$, we define a hyperbolic metric $\bar{g}$ on $(r_0,\infty) \times T^{n-1}$ by $\bar{g} = r^{-2} \, dr \otimes dr + r^2 \, \gamma$. 

An $(N,n)$-dataset is a triplet $(M,g,\rho)$ consisting of a complete, connected, orientable manifold $M$ of dimension $n$, a Riemannian metric $g$ on $M$, and a positive smooth function $\rho$ on $M$ satisfying the following conditions: 
\begin{itemize} 
\item There exists a bounded open domain $E \subset M$ with smooth boundary such that the complement $M \setminus E$ is diffeomorphic to $[r_0,\infty) \times T^{n-1}$.
\item The map 
\[(r_0,\infty) \times T^{n-1} \to T^{n-2}, \quad (r,\theta_0,\hdots,\theta_{n-2}) \mapsto (\theta_1,\hdots,\theta_{n-2})\] 
extends to a globally defined smooth map from $M$ to $T^{n-2}$.
\item For every nonnegative integer $m$, the metric $g$ satisfies 
\[|\bar{D}^m (g-\bar{g})|_{\bar{g}} \leq O(r^{-N}),\] 
where $\bar{D}^m$ denotes the covariant derivative of order $m$ with respect to the hyperbolic metric $\bar{g}$. 
\item The metric $g$ satisfies 
\[|g-\bar{g} - r^{2-N} \, Q|_{\bar{g}} \leq O(r^{-N-2\delta}).\] 
Here, $\delta$ is a small positive real number and $Q$ is a smooth symmetric $(0,2)$-tensor on $T^{n-1}$.
\item For every nonnegative integer $m$, the function $\rho$ satisfies 
\[|\bar{D}^m(\rho - r^{N-n})|_{\bar{g}} \leq O(r^{-n}),\] 
where $\bar{D}^m$ denotes the covariant derivative of order $m$ with respect to the hyperbolic metric $\bar{g}$. 
\item The function $\rho$ satisfies 
\[|\rho - r^{N-n} - r^{-n} \, P| \leq O(r^{-n-2\delta}).\] 
Here, $\delta$ is a small positive number and $P$ is a real-valued function on $T^{n-1}$ which is H\"older continuous with exponent $2\delta$.
\item We have 
\[\int_{T^{n-1}} \Big ( N \, \text{\rm tr}_\gamma(Q) + 2N \, P + \Big ( \frac{2}{Nb_0} \Big )^N \Big ) \, d\text{\rm vol}_\gamma \leq 0.\] 
\end{itemize}
\end{definition}

\begin{definition}
\label{model.dataset}
Let $n$ be an integer with $2 \leq n \leq N$ and let $(M,g,\rho)$ be an $(N,n)$-dataset. We say that $(M,g,\rho)$ is a model $(N,n)$-dataset if there exists a smooth map $\Phi: \mathbb{R}^2 \times \mathbb{R}^{n-2} \to M$ such that $\Phi^* g = g_{\text{\rm HM},N,n}$ and the function $\rho \circ \Phi$ is a constant multiple of $\rho_{\text{\rm HM},N,n}$.
\end{definition}

\begin{definition} 
\label{star.N.n}
Let $n$ be an integer with $2 \leq n \leq N$. We say that condition $(\star_{N,n})$ is satisfied if the following holds. Let $(M,g,\rho)$ be an arbitrary $(N,n)$-dataset. If $n=N$, we assume that $\rho=1$ and $R \geq -N(N-1)$ at each point in $M$. If $n<N$, we assume that  
\[-2 \, \Delta \log \rho - \frac{N-n+1}{N-n} \, |\nabla \log \rho|^2 + R + N(N-1) \geq 0\] 
at each point in $M$. Then $(M,g,\rho)$ is a model $(N,n)$-dataset. 
\end{definition}

We next state several lemmata that will be needed later.

\begin{lemma}
\label{derivatives.of.g}
Let $n$ be an integer with $2 \leq n \leq N$. Let $(M,g,\rho)$ be an $(N,n)$-dataset. Then $|\bar{D}^m(g-\bar{g}-r^{2-N} \, Q)|_{\bar{g}} \leq O(r^{-N-\delta})$ for every nonnegative integer $m$. 
\end{lemma}

\textbf{Proof.} 
By assumption, $|\bar{D}^m(g-\bar{g})|_{\bar{g}} \leq O(r^{-N})$ for every nonnegative integer $m$. Since $Q$ is a smooth tensor on $T^{n-1}$, it follows that 
\[|\bar{D}^m(g-\bar{g}-r^{2-N} \, Q)|_{\bar{g}} \leq O(r^{-N})\] 
for every nonnegative integer $m$. Moreover, our assumptions imply that 
\[|g-\bar{g}-r^{2-N} \, Q|_{\bar{g}} \leq O(r^{-N-2\delta}).\] 
The assertion now follows from standard interpolation inequalities. \\

\begin{lemma}
\label{Lie.derivative.of.metric.along.angular.vector.field}
Let $n$ be an integer with $2 \leq n \leq N$. Let $(M,g,\rho)$ be an $(N,n)$-dataset. Let $V$ be a smooth vector field on $M$ with the property that $V = \frac{\partial}{\partial \theta_{n-2}}$ outside a compact set. Then $|\mathscr{L}_V g| \leq O(r^{1-N-\delta})$ and $|\mathscr{L}_V \mathscr{L}_V g| \leq O(r^{2-N-\delta})$.
\end{lemma} 

\textbf{Proof.} 
It follows from Lemma \ref{derivatives.of.g} that 
\[|\mathscr{L}_V(g-\bar{g}-r^{2-N} \, Q)|_{\bar{g}} \leq O(r^{1-N-\delta})\] 
and 
\[|\mathscr{L}_V \mathscr{L}_V (g-\bar{g}-r^{2-N} \, Q)|_{\bar{g}} \leq O(r^{2-N-\delta}).\] 
Moreover, since $Q$ is a smooth tensor on $T^{n-1}$, we know that 
\[|\mathscr{L}_V (\bar{g}+r^{2-N} \, Q)|_{\bar{g}} = r^{2-N} \, |\mathscr{L}_{\frac{\partial}{\partial \theta_{n-2}}} Q|_{\bar{g}} \leq O(r^{-N})\] 
and 
\[|\mathscr{L}_V \mathscr{L}_V (\bar{g}+r^{2-N} \, Q)|_{\bar{g}} = r^{2-N} \, |\mathscr{L}_{\frac{\partial}{\partial \theta_{n-2}}} \mathscr{L}_{\frac{\partial}{\partial \theta_{n-2}}} Q|_{\bar{g}} \leq O(r^{-N})\] 
outside a compact set. Putting these facts together, the assertion follows. This completes the proof of Lemma \ref{Lie.derivative.of.metric.along.angular.vector.field}. \\

\begin{lemma}
\label{derivative.of.rho.along.angular.vector.field}
Let $n$ be an integer with $2 \leq n \leq N$. Let $(M,g,\rho)$ be an $(N,n)$-dataset. Let $V$ be a smooth vector field on $M$ with the property that $V = \frac{\partial}{\partial \theta_{n-2}}$ outside a compact set. Then $|V(\rho)| \leq O(r^{1-n-\delta})$ and $|V(V(\rho))| \leq O(r^{2-n-\delta})$.
\end{lemma} 

\textbf{Proof.} 
Let us fix a large number $\bar{r}$ and a point $p$ on the level set $\{r=\bar{r}\}$. Let $\varphi_s$ denote the flow generated by the vector field $\bar{r}^{-1} \, V$. Since the function $P: T^{n-1} \to \mathbb{R}$ is H\"older continuous with exponent $2\delta$, it follows that 
\begin{equation} 
\label{oscillation.of.P}
\sup_{s \in [-2,2]} P(\varphi_s(p)) - \inf_{s \in [-2,2]} P(\varphi_s(p)) \leq C \, \bar{r}^{-2\delta}. 
\end{equation}
Using the estimate $|\rho - r^{N-n} - r^{-n} \, P| \leq C \, r^{-n-2\delta}$, we obtain 
\begin{equation} 
\label{estimate.for.rho}
\sup_{s \in [-2,2]} |\rho(\varphi_s(p)) - \bar{r}^{N-n} - \bar{r}^{-n} \, P(\varphi_s(p))| \leq C \, \bar{r}^{-n-2\delta}. 
\end{equation}
Combining (\ref{oscillation.of.P}) and (\ref{estimate.for.rho}), we deduce that 
\begin{equation} 
\label{bound.for.oscillation.of.rho}
\sup_{s \in [-2,2]} \rho(\varphi_s(p)) - \inf_{s \in [-2,2]} \rho(\varphi_s(p)) \leq C \, \bar{r}^{-n-2\delta}. 
\end{equation}
By assumption, $|\bar{D}^m(\rho - r^{N-n})|_{\bar{g}} \leq O(r^{-n})$ for every nonnegative integer $m$. This implies 
\[|\underbrace{V \cdots V}_{\text{\rm $m$ times}}(\rho - r^{N-n})| \leq O(r^{m-n})\] 
for every positive integer $m$. Consequently, 
\begin{equation} 
\label{bound.for.higher.derivatives.of.rho}
\sup_{s \in [-2,2]} \Big | \frac{d^m}{ds^m} \rho(\varphi_s(p)) \Big | \leq C(m) \, \bar{r}^{-n} 
\end{equation}
for every positive integer $m$. Using (\ref{bound.for.oscillation.of.rho}), (\ref{bound.for.higher.derivatives.of.rho}), and standard interpolation inequalities, we conclude that 
\[\sup_{s \in [-1,1]} \Big | \frac{d}{ds} \rho(\varphi_s(p)) \Big | \leq C \, \bar{r}^{-n-\delta}\] 
and 
\[\sup_{s \in [-1,1]} \Big | \frac{d^2}{ds^2} \rho(\varphi_s(p)) \Big | \leq C \, \bar{r}^{-n-\delta}.\] 
Putting $s=0$, it follows that $|V(\rho)| \leq C \, \bar{r}^{1-n-\delta}$ and $|V(V(\rho))| \leq C \, \bar{r}^{2-n-\delta}$ at the point $p$. This completes the proof of Lemma \ref{derivative.of.rho.along.angular.vector.field}. \\

\begin{lemma}
\label{covariant.derivative.angular.vector.field}
Let $n$ be an integer with $2 \leq n \leq N$. Let $(M,g,\rho)$ be an $(N,n)$-dataset. Then 
\[\Big | D_{\frac{\partial}{\partial \theta_{n-2}}} \frac{\partial}{\partial \theta_{n-2}} + \Big ( b_{n-2}^2 \, r - \frac{N-2}{2} \, r^{1-N} \, Q(\frac{\partial}{\partial \theta_{n-2}},\frac{\partial}{\partial \theta_{n-2}}) \Big ) \, \nabla r \Big | \leq o(r^{2-N}).\] 
\end{lemma} 

\textbf{Proof.} 
For abbreviation, we define a metric $\hat{g}$ by 
\[\hat{g} = \bar{g} + r^{2-N} \, Q = r^{-2} \, dr \otimes dr + r^2 \, \gamma + r^{2-N} \, Q.\] 
Let $\hat{D}$ denote the Levi-Civita connection with respect to the metric $\hat{g}$. We compute 
\begin{align*} 
\hat{g}(\hat{D}_{\frac{\partial}{\partial \theta_{n-2}}} \frac{\partial}{\partial \theta_{n-2}},\frac{\partial}{\partial r}) 
&= \frac{\partial}{\partial \theta_{n-2}} \hat{g}(\frac{\partial}{\partial \theta_{n-2}},\frac{\partial}{\partial r}) - \frac{1}{2} \, \frac{\partial}{\partial r} \hat{g}(\frac{\partial}{\partial \theta_{n-2}},\frac{\partial}{\partial \theta_{n-2}}) \\ 
&= -b_{n-2}^2 \, r + \frac{N-2}{2} \, r^{1-N} \, Q(\frac{\partial}{\partial \theta_{n-2}},\frac{\partial}{\partial \theta_{n-2}}) 
\end{align*}
and 
\begin{align*} 
\hat{g}(\hat{D}_{\frac{\partial}{\partial \theta_{n-2}}} \frac{\partial}{\partial \theta_{n-2}},\frac{\partial}{\partial \theta_k}) 
&= \frac{\partial}{\partial \theta_{n-2}} \hat{g}(\frac{\partial}{\partial \theta_{n-2}},\frac{\partial}{\partial \theta_k}) - \frac{1}{2} \, \frac{\partial}{\partial \theta_k} \hat{g}(\frac{\partial}{\partial \theta_{n-2}},\frac{\partial}{\partial \theta_{n-2}}) \\ 
&= O(r^{2-N}) 
\end{align*} 
for $k=0,\hdots,n-2$. From this, we deduce that 
\[\Big | \hat{D}_{\frac{\partial}{\partial \theta_{n-2}}} \frac{\partial}{\partial \theta_{n-2}} + \Big ( b_{n-2}^2 \, r - \frac{N-2}{2} \, r^{1-N} \, Q(\frac{\partial}{\partial \theta_{n-2}},\frac{\partial}{\partial \theta_{n-2}}) \Big ) \, \hat{\nabla} r \Big | \leq O(r^{1-N}).\] 
Finally, $|g-\hat{g}|_{\bar{g}} \leq o(r^{-N})$ and $|\bar{D}(g-\hat{g})|_{\bar{g}} \leq o(r^{-N})$ by Lemma \ref{derivatives.of.g}. Putting these facts together, the assertion follows. This completes the proof of Lemma \ref{covariant.derivative.angular.vector.field}. \\

\begin{definition} 
\label{definition.tame}
Let $n$ be an integer with $3 \leq n \leq N$. Let $(M,g,\rho)$ be an $(N,n)$-dataset. Let $\Sigma$ be a properly embedded, connected, orientable hypersurface in $M$, and let $t_* \in S^1$. We say that $\Sigma$ is $t_*$-tame if there exists a large number $r_*$, and a function $f: [r_*,\infty) \times T^{n-2} \to S^1$ with the following properties: 
\begin{itemize}
\item $\Sigma \cap \{r \geq r_*\} = \{\theta_{n-2} = f(r,\theta_0,\hdots,\theta_{n-3})\}$.
\item We have $d_{S^1}(f(r,\theta_0,\hdots,\theta_{n-3}),t_*) \leq O(r^{-N})$, where $d_{S^1}$ denotes the Riemannian distance on $S^1$.
\item The higher order covariant derivatives of the function $f: [r_*,\infty) \times T^{n-2} \to S^1$ with respect to the hyperbolic metric $r^{-2} \, dr \otimes dr + \sum_{k=0}^{n-3} b_k^2 \, r^2 \, d\theta_k \otimes d\theta_k$ are bounded by $O(r^{-N})$. 
\end{itemize}
Finally, we say that $\Sigma$ is tame if $\Sigma$ is $t_*$-tame for some element $t_* \in S^1$.
\end{definition}

\begin{definition}
Let $n$ be an integer with $3 \leq n \leq N$. Let $(M,g,\rho)$ be an $(N,n)$-dataset, and let $\Sigma$ be a compact, orientable hypersurface in $M$. The $(g,\rho)$-area of $\Sigma$ is defined as $\int_\Sigma \rho \, d\text{\rm vol}_g$.
\end{definition}

\begin{definition}
\label{definition.stationary} 
Let $n$ be an integer with $3 \leq n \leq N$. Let $(M,g,\rho)$ be an $(N,n)$-dataset, and let $\Sigma$ be an orientable hypersurface in $M$. We say that $\Sigma$ is $(g,\rho)$-stationary if $H_\Sigma + \rho^{-1} \, \langle \nabla \rho,\nu_\Sigma \rangle = 0$ at each point on $\Sigma$. Here, $\nu_\Sigma$ denotes the unit normal vector field to $\Sigma$ and $H_\Sigma$ denotes the mean curvature of $\Sigma$. 
\end{definition}

\begin{definition}
\label{definition.stability}
Let $n$ be an integer with $3 \leq n \leq N$. Let $(M,g,\rho)$ be an $(N,n)$-dataset. Let $\Sigma$ be a properly embedded, connected, orientable hypersurface in $M$ which is $t_*$-tame for some $t_* \in S^1$. Suppose that $\Sigma$ is $(g,\rho)$-stationary. Moreover, suppose that $u: T^{n-1} \to \mathbb{R}$ is twice continuously differentiable. We say that $\Sigma$ is $(g,\rho,u)$-stable if the following holds. If $a$ is a real number and $V$ is a smooth vector field on $M$ with the property that $V = a \, \frac{\partial}{\partial \theta_{n-2}}$ outside a compact set, then 
\begin{align*} 
&\frac{1}{2} \int_\Sigma \rho \, \sum_{k=1}^{n-1} (\mathscr{L}_V \mathscr{L}_V g)(e_k,e_k) + \int_\Sigma V(V(\rho)) \\ 
&- \frac{1}{2} \int_\Sigma \rho \, \sum_{k,l=1}^{n-1} (\mathscr{L}_V g)(e_k,e_l) \, (\mathscr{L}_V g)(e_k,e_l) \\ 
&+ \frac{1}{4} \int_\Sigma \rho \, \sum_{k,l=1}^{n-1} (\mathscr{L}_V g)(e_k,e_k) \, (\mathscr{L}_V g)(e_l,e_l) \\ 
&+ \int_\Sigma V(\rho) \, \sum_{k=1}^{n-1} (\mathscr{L}_V g)(e_k,e_k) \\ 
&\geq -a^2 \int_{T^{n-2} \times \{t_*\}} \frac{\partial^2 u}{\partial \theta_{n-2}^2} \, d\text{\rm vol}_\gamma. 
\end{align*} 
Here, $\{e_1,\hdots,e_{n-1}\}$ denotes a local orthonormal frame on $\Sigma$. The integrals on the left hand side are well-defined in view of Lemma \ref{Lie.derivative.of.metric.along.angular.vector.field} and Lemma \ref{derivative.of.rho.along.angular.vector.field} and the assumption that $\Sigma$ is tame. 
\end{definition}

\begin{definition} 
\label{definition.weighted.Jacobi.operator}
Let $n$ be an integer with $3 \leq n \leq N$. Let $(M,g,\rho)$ be an $(N,n)$-dataset, and let $\Sigma$ be an orientable hypersurface in $M$. Given a smooth function $v: \Sigma \to \mathbb{R}$, we define 
\begin{align*} 
\mathbb{L}_\Sigma v 
&= -\text{\rm div}_\Sigma(\rho \, \nabla^\Sigma v) - \rho \, (\text{\rm Ric}(\nu_\Sigma,\nu_\Sigma) + |h_\Sigma|^2) \, v \\ 
&+ (D^2 \rho)(\nu_\Sigma,\nu_\Sigma) \, v - \rho^{-1} \, \langle \nabla \rho,\nu_\Sigma \rangle^2 \, v. 
\end{align*} 
Here, $\nu_\Sigma$ denotes the unit normal vector field to $\Sigma$ and $h_\Sigma$ denotes the second fundamental form of $\Sigma$. The operator $\mathbb{L}_\Sigma$ is referred to as the weighted Jacobi operator of $\Sigma$.
\end{definition}

\section{Proof of property $(\star_{N,2})$ for each $N \geq 3$}

\label{2D}

\begin{theorem}
\label{2D.case}
Let $N \geq 3$ be an integer. Then property $(\star_{N,2})$ is satisfied.
\end{theorem}

In the remainder of this section, we will describe the proof of Theorem \ref{2D.case}. Let $(M,g,\rho)$ be a $(N,2)$-dataset satisfying 
\[-2 \, \Delta \log \rho - \frac{N-1}{N-2} \, |\nabla \log \rho|^2 + 2K + N(N-1) \geq 0\] 
at each point in $M$. Let $\psi = \log \rho$. Then 
\[-2 \, \Delta \psi - \frac{N-1}{N-2} \, |\nabla \psi|^2 + 2K + N(N-1) \geq 0\] 
at each point in $M$. Moreover, the function $\psi$ satisfies 
\[|\bar{D}^m(\psi - (N-2) \, \log r)|_{\bar{g}} \leq O(r^{-N})\] 
for every nonnegative integer $m$, and 
\[|\psi - (N-2) \, \log r - r^{-N} \, P| \leq O(r^{-N-2\delta}).\] 
In particular, the function $\psi$ is proper, and the set of critical points of $\psi$ is compact.

\begin{lemma}
\label{level.sets.of.psi}
If $s$ is sufficiently large, then the level set $\{\psi=s\}$ is a one-dimensional submanifold diffeomorphic to $S^1$.
\end{lemma}

\textbf{Proof.} 
By assumption, $|d(\psi - (N-2) \, \log r)| \leq O(r^{-N})$. This implies $|\frac{\partial}{\partial r}(\psi - (N-2) \, \log r)| \leq O(r^{-N-1})$. Consequently, we can find a large number $r_*$ such that $\frac{\partial}{\partial r} \psi > 0$ on the set $\{r \geq r_*\}$. If $s$ is sufficiently large, then the level set $\{\psi=s\}$ is contained in the region $\{r \geq r_*\}$. Moreover, for each $t \in S^1$, the curve $[r_*,\infty) \times \{t\}$ intersects the level set $\{\psi=s\}$ exactly once and the intersection is transversal. Therefore, if $s$ is sufficiently large, then the level set $\{\psi=s\}$ is diffeomorphic to $S^1$. This completes the proof of Lemma \ref{level.sets.of.psi}. \\

\begin{lemma} 
\label{choice.of.r_j}
We can find a sequence $r_j \to \infty$ such that 
\[\liminf_{j \to \infty} r_j^{N-1} \int_{\{r=r_j\}} \langle \nabla \psi - (N-2) \, \nabla \log r,\eta \rangle \geq -\int_{S^1} N \, P \, d\text{\rm vol}_\gamma.\] 
Here, $\eta = \frac{\nabla r}{|\nabla r|}$ denotes the outward-pointing unit normal vector field along the level set $\{r=r_j\}$.
\end{lemma}

\textbf{Proof.} 
Using the estimates 
\[|\text{\rm div}(r^{N-1} \, \nabla r) - (N+1) \, r^N| \leq C\] 
and 
\[|\psi - (N-2) \, \log r| \leq C \, r^{-N},\] 
we obtain 
\begin{align*} 
&r^{N-1} \, \langle \nabla \psi - (N-2) \, \nabla \log r,\nabla r \rangle \\ 
&+ (N+1) \, r^N \, (\psi - (N-2) \, \log r) \\ 
&- \text{\rm div} \big ( r^{N-1} \, (\psi - (N-2) \, \log r) \, \nabla r \big ) \\ 
&= -\big ( \text{\rm div}(r^{N-1} \, \nabla r) - (N+1) \, r^N \big ) \, (\psi - (N-2) \, \log r) \\ 
&\geq - C \, r^{-N}. 
\end{align*} 
In the next step, we integrate this inequality over the domain $\{2r_0 \leq r \leq \bar{r}\}$, where $\bar{r} > 2r_0$. Using the divergence theorem, it follows that  
\begin{align*} 
&\int_{\{2r_0 \leq r \leq \bar{r}\}} r^{N-1} \, \langle \nabla \psi - (N-2) \, \nabla \log r,\nabla r \rangle \\ 
&+ \int_{\{2r_0 \leq r \leq \bar{r}\}} (N+1) \, r^N \, (\psi - (N-2) \, \log r) \\ 
&- \int_{\{r=\bar{r}\}} r^{N-1} \, (\psi - (N-2) \, \log r) \, |\nabla r| \geq -C 
\end{align*}
for $\bar{r} > 2r_0$. To estimate the second and third term on the left hand side, we use the inequality 
\[|\psi - (N-2) \, \log r - r^{-N} \, P| \leq C \, r^{-N-2\delta}.\] 
This gives 
\begin{align*} 
&\int_{\{2r_0 \leq r \leq \bar{r}\}} r^{N-1} \, \langle \nabla \psi - (N-2) \, \nabla \log r,\nabla r \rangle \\ 
&+ \bar{r} \int_{S^1} N \, P \, d\text{\rm vol}_\gamma \geq -C \, \bar{r}^{1-\delta} 
\end{align*} 
for $\bar{r} > 2r_0$. Using the co-area formula, we conclude that 
\[\limsup_{\bar{r} \to \infty} \int_{\{r=\bar{r}\}} r^{N-1} \, \Big \langle \nabla \psi - (N-2) \, \nabla \log r,\frac{\nabla r}{|\nabla r|} \Big \rangle \geq -\int_{S^1} N \, P \, d\text{\rm vol}_\gamma.\] 
Since $\eta = \frac{\nabla r}{|\nabla r|}$, the assertion follows. This completes the proof of Lemma \ref{choice.of.r_j}. \\

In the following, we assume that the sequence $r_j$ is chosen as in Lemma \ref{choice.of.r_j}. We define $M^{(j)} = M \setminus \{r>r_j\}$. We denote by $\kappa$ the geodesic curvature of the boundary $\partial M^{(j)}$. 

\begin{proposition} 
\label{choice.of.r_j.2}
We have 
\[\liminf_{j \to \infty} 2 \, |\partial M^{(j)}|^{N-1} \, \int_{\partial M^{(j)}} (\langle \nabla \psi,\eta \rangle + \kappa - (N-1)) \geq \Big ( \frac{4\pi}{N} \Big )^N.\] 
\end{proposition}

\textbf{Proof.} 
It follows from Lemma \ref{choice.of.r_j} that 
\begin{equation} 
\label{limit.a} 
\liminf_{j \to \infty} r_j^{N-1} \int_{\partial M^{(j)}} \langle \nabla \psi - (N-2) \, \nabla \log r,\eta \rangle \geq -\int_{S^1} N \, P \, d\text{\rm vol}_\gamma. 
\end{equation}
The estimate $|g-\bar{g} - r^{2-N} \, Q|_{\bar{g}} \leq O(r^{-N-2\delta})$ implies $|\nabla \log r| = 1 + O(r^{-N-\delta})$. Consequently, 
\begin{equation} 
\label{limit.b}
\lim_{j \to \infty} r_j^{N-1} \int_{\partial M^{(j)}} ((N-2) \, \langle \nabla \log r,\eta \rangle - (N-2)) = 0. 
\end{equation}
Using Lemma \ref{covariant.derivative.angular.vector.field}, we compute 
\[\kappa-1 = -\frac{N}{2} \, \text{\rm tr}_\gamma(Q) \, r_j^{-N} + o(r_j^{-N}),\] 
where $\text{\rm tr}_\gamma(Q) = b_0^{-2} \, Q(\frac{\partial}{\partial \theta_0},\frac{\partial}{\partial \theta_0})$. Integrating this identity over $\partial M^{(j)}$ gives 
\begin{equation} 
\label{limit.c} 
\lim_{j \to \infty} r_j^{N-1} \int_{\partial M^{(j)}} (\kappa - 1) = -\int_{S^1} \frac{N}{2} \, \text{\rm tr}_\gamma(Q) \, d\text{\rm vol}_\gamma. 
\end{equation}
Adding (\ref{limit.a}), (\ref{limit.b}), and (\ref{limit.c}), we obtain 
\[\liminf_{j \to \infty} r_j^{N-1} \int_{\partial M^{(j)}} (\langle \nabla \psi,\eta \rangle + \kappa - (N-1)) \geq -\int_{S^1} \Big ( \frac{N}{2} \, \text{\rm tr}_\gamma(Q) + N \, P \Big ) \, d\text{\rm vol}_\gamma.\] 
On the other hand, 
\[\int_{S^1} \Big ( N \, \text{\rm tr}_\gamma(Q) + 2N \, P + \Big ( \frac{2}{Nb_0} \Big )^N \Big ) \, d\text{\rm vol}_\gamma \leq 0\] 
by definition of a $(N,2)$-dataset. Note that $(S^1,\gamma)$ has length $2\pi b_0$. Consequently, 
\[\liminf_{j \to \infty} 2 \, r_j^{N-1} \int_{\partial M^{(j)}} (\langle \nabla \psi,\eta \rangle + \kappa - (N-1)) \geq 2\pi b_0 \, \Big ( \frac{2}{Nb_0} \Big )^N.\] 
Since $\lim_{j \to \infty} r_j^{-1} \, |\partial M^{(j)}| = 2\pi b_0$, we conclude that 
\[\liminf_{j \to \infty} 2 \, |\partial M^{(j)}|^{N-1} \int_{\partial M^{(j)}} (\langle \nabla \psi,\eta \rangle + \kappa - (N-1)) \geq \Big ( \frac{4\pi}{N} \Big )^N.\] 
This completes the proof of Proposition \ref{choice.of.r_j.2}. \\

For each $j$, we denote by $w^{(j)}: M^{(j)} \to [0,\infty)$ the distance function from $\partial M^{(j)}$. Note that the function $w^{(j)}$ is Lipschitz continuous with Lipschitz constant $1$. For each $j$, we define $l^{(j)} = \sup_{M^{(j)}} w^{(j)}$.

\begin{lemma}
\label{estimate.for.w_j}
We have $\sup_{M^{(j)} \cap \{r>2r_0\}} |w^{(j)} - \log r_j + \log r| \leq C$, where $C$ is a constant that does not depend on $j$. In particular, $|l^{(j)} - \log r_j| \leq C$, where $C$ is a constant that does not depend on $j$. 
\end{lemma}

\textbf{Proof.} 
This follows from our assumptions on the metric $g$. \\

We now continue the proof of Theorem \ref{2D.case}. For each $s \in [0,l^{(j)})$, we define $\Omega_s^{(j)} = \{w^{(j)} > s\}$. For almost every $s \in [0,l^{(j)})$, the boundary of $\Omega_s^{(j)}$ is a piecewise smooth curve. We denote the length of $\partial \Omega_s^{(j)}$ by $L^{(j)}(s)$.

As in Section 2 of \cite{Brendle-Hung}, we define a function $F: (0,\infty) \to (0,1)$ by 
\[F(s) = \tanh \Big ( \frac{Ns}{2} \Big ) = \frac{\sinh(Ns)}{1+\cosh(Ns)}\] 
for each $s \in (0,\infty)$. We define a function $G: (0,\infty) \to (1,\infty)$ by 
\[G(s) = \Big [ \cosh \Big ( \frac{Ns}{2} \Big ) \Big ]^{\frac{2(N-1)}{N}} = \Big [ \frac{1+\cosh(Ns)}{2} \Big ]^{\frac{N-1}{N}}\] 
for each $s \in (0,\infty)$. Moreover, we define 
\[I^{(j)}(s) = 2\pi \chi(M^{(j)}) - (N-1) \, F(l^{(j)}-s) \, L^{(j)}(s) + \int_{\Omega_s^{(j)}} (\Delta \psi - K)\] 
and 
\[J^{(j)}(s) = G(l^{(j)}-s) \, I^{(j)}(s)\] 
for $s \in [0,l^{(j)})$. 

\begin{proposition}
\label{initial.value.of.J}
We have $\liminf_{j \to \infty} J^{(j)}(0) \geq 2\pi$.
\end{proposition}

\textbf{Proof.} 
It follows from Young's inequality that $Nab \leq a^N + (N-1) \, b^{\frac{N}{N-1}}$ for all $a,b \geq 0$. Putting $a = \frac{4\pi}{N}$ and $b = G(l^{(j)})^{-1} \, |\partial M^{(j)}|^{N-1}$ gives 
\[4\pi \, G(l^{(j)})^{-1} \, |\partial M^{(j)}|^{N-1} \leq \Big ( \frac{4\pi}{N} \Big )^N + (N-1) \, G(l^{(j)})^{-\frac{N}{N-1}} \, |\partial M^{(j)}|^N.\] 
Using the estimate 
\[G(l^{(j)})^{-\frac{N}{N-1}} = 1-F(l^{(j)})^2 \leq 2 \, (1-F(l^{(j)})),\] 
we obtain 
\[4\pi \, G(l^{(j)})^{-1} \, |\partial M^{(j)}|^{N-1} \leq \Big ( \frac{4\pi}{N} \Big )^N + 2(N-1) \, (1-F(l^{(j)})) \, |\partial M^{(j)}|^N.\] 
On the other hand, the Gauss-Bonnet theorem gives 
\[I^{(j)}(0) = -(N-1) \, F(l^{(j)}) \, |\partial M^{(j)}| + \int_{\partial M^{(j)}} (\langle \nabla \psi,\eta \rangle + \kappa).\] 
Putting these facts together, we conclude that 
\begin{align*} 
&2 \, G(l^{(j)})^{-1} \, |\partial M^{(j)}|^{N-1} \, (2\pi - J^{(j)}(0)) \\ 
&= 4\pi \, G(l^{(j)})^{-1} \, |\partial M^{(j)}|^{N-1} + 2(N-1) \, F(l^{(j)}) \, |\partial M^{(j)}|^N \\ 
&- 2 \, |\partial M^{(j)}|^{N-1} \, \int_{\partial M^{(j)}} (\langle \nabla \psi,\eta \rangle + \kappa) \\ 
&\leq \Big ( \frac{4\pi}{N} \Big )^N - 2 \, |\partial M^{(j)}|^{N-1} \, \int_{\partial M^{(j)}} (\langle \nabla \psi,\eta \rangle + \kappa - (N-1)). 
\end{align*}
Note that $|l^{(j)} - \log r_j| \leq C$ by Lemma \ref{estimate.for.w_j}. Consequently, the sequence $G(l^{(j)})^{-1} \, |\partial M^{(j)}|^{N-1}$ is uniformly bounded from above and below by positive constants. Using Proposition \ref{choice.of.r_j.2}, we conclude that $\limsup_{j \to \infty} (2\pi - J^{(j)}(0)) \leq 0$, as claimed. \\

\begin{proposition}
\label{terminal.value.of.J}
For each $j$, we have $\limsup_{s \nearrow l^{(j)}} J^{(j)}(s) \leq 2\pi$.
\end{proposition} 

\textbf{Proof.} 
Since $M^{(j)}$ is connected, it follows that $\chi(M^{(j)}) \leq 1$. Consequently, 
\[\limsup_{s \nearrow l^{(j)}} I^{(j)}(s) \leq 2\pi \chi(M^{(j)}) \leq 2\pi.\] 
Since $\lim_{s \nearrow l^{(j)}} G(l^{(j)}-s) = 1$, we conclude that 
\[\limsup_{s \nearrow l^{(j)}} J^{(j)}(s) \leq 2\pi.\] 
This completes the proof of Proposition \ref{terminal.value.of.J}. \\

\begin{proposition}
\label{consequence.of.monotonicity}
For each $j$, we have 
\begin{align*} 
&\int_{M^{(j)}} G(l^{(j)}-w^{(j)}) \, \Big ( -2 \, \Delta \psi - \frac{N-1}{N-2} \, |\nabla \psi|^2 + 2K + N(N-1) \Big ) \\ 
&+ \int_{M^{(j)}} \frac{N-1}{N-2} \, G(l^{(j)}-w^{(j)}) \, |(N-2) \, F(l^{(j)}-w^{(j)}) \, \nabla w^{(j)} + \nabla \psi|^2 \\ 
&\leq 2 \, (2\pi - J^{(j)}(0))
\end{align*}
\end{proposition} 

\textbf{Proof.} 
For each $j$, we can find a large constant $C_j$ with the property that the function $s \mapsto J^{(j)}(s) + C_j s$ is monotone increasing (see \cite{Brendle-Hung}, Section 2). Moreover, 
\begin{align*} 
&\int_{\partial \Omega_s^{(j)}} G(l^{(j)}-s) \, \Big ( -2 \, \Delta \psi - \frac{N-1}{N-2} \, |\nabla \psi|^2 + 2K + N(N-1) \Big ) \\ 
&+ \int_{\partial \Omega_s^{(j)}} \frac{N-1}{N-2} \, G(l^{(j)}-s) \, |(N-2) \, F(l^{(j)}-s) \, \nabla w^{(j)} + \nabla \psi|^2 \\ 
&\leq 2 \, \frac{d}{ds} J^{(j)}(s) 
\end{align*}
for almost every $s \in (0,l^{(j)})$ (see \cite{Brendle-Hung}, Section 2). We integrate this inequality over the interval $(0,l^{(j)})$. Using Proposition \ref{terminal.value.of.J}, the assertion follows. This completes the proof of Proposition \ref{consequence.of.monotonicity}. \\

It follows from Lemma \ref{estimate.for.w_j} that $\sup_{M^{(j)} \cap \{r>2r_0\}} |l^{(j)} - w^{(j)} - \log r| \leq C$, where $C$ is a constant that does not depend on $j$. After passing to a subsequence, we may assume that the functions $l^{(j)}-w^{(j)}$ converge in $C_{\text{\rm loc}}^0$ to some limiting function $w: M \to \mathbb{R}$. The function $w$ is Lipschitz continuous with Lipschitz constant $1$. Moreover, $\sup_{(2r_0,\infty) \times S^1} |w - \log r| \leq C$. In particular, the function $w$ is proper. Finally, since $\inf_{M^{(j)}} (l^{(j)}-w^{(j)}) = 0$ for each $j$, we conclude that $\inf_M w = 0$.

\begin{lemma}
\label{zero.level.set.has.empty.interior}
The level set $\{w=0\}$ has empty interior.  
\end{lemma}

\textbf{Proof.} 
Let us fix a point $y \in M$ with $w(y) = 0$. Let $\delta$ be an arbitrary positive real number. If $j$ is sufficiently large, we can find a point $y^{(j)} \in M$ such that $d_{(M,g)}(y,y^{(j)}) = \delta$ and $w^{(j)}(y) - w^{(j)}(y^{(j)}) = \delta$. After passing to a subsequence, we may assume that the sequence $y^{(j)}$ converges to a point $\tilde{y} \in M$. Then $d_{(M,g)}(y,\tilde{y}) = \lim_{j \to \infty} d_{(M,g)}(y,y^{(j)}) = \delta$ and $w(\tilde{y}) - w(y) = \lim_{j \to \infty} (w^{(j)}(y) - w^{(j)}(y^{(j)})) = \delta$. Since $\delta>0$ is arbitrary, the assertion follows. This completes the proof of Lemma \ref{zero.level.set.has.empty.interior}. \\

\begin{proposition} 
\label{relation.between.w.and.psi}
The function $w$ satisfies 
\[\frac{(N-2) \, \sinh(Nw)}{1+\cosh(Nw)} = |\nabla \psi|\] 
and 
\[\frac{N-2}{N} \, \log(1+\cosh(Nw)) = \psi + c\] 
at each point on $M$, where $c$ is a constant. 
\end{proposition}

\textbf{Proof.} 
Using Proposition \ref{initial.value.of.J} and Proposition \ref{consequence.of.monotonicity}, we obtain 
\begin{align*} 
&\int_{M^{(j)}} G(l^{(j)}-w^{(j)}) \, \big | (N-2) \, F(l^{(j)}-w^{(j)}) - |\nabla \psi| \big |^2 \\ 
&\leq \int_{M^{(j)}} G(l^{(j)}-w^{(j)}) \, |(N-2) \, F(l^{(j)}-w^{(j)}) \, \nabla w^{(j)} + \nabla \psi|^2 \to 0 
\end{align*}
and 
\begin{align*} 
&\int_{M^{(j)}} G(l^{(j)}-w^{(j)}) \, \Big | \nabla \Big ( \frac{N-2}{N} \, \log \big ( 1+\cosh(N(l^{(j)}-w^{(j)})) \big ) - \psi \Big ) \Big |^2 \\ 
&= \int_{M^{(j)}} G(l^{(j)}-w^{(j)}) \, |(N-2) \, F(l^{(j)}-w^{(j)}) \, \nabla w^{(j)} + \nabla \psi|^2 \to 0 
\end{align*}
as $j \to \infty$. From this, the assertion follows. \\

\begin{proposition}
\label{pde.for.psi}
The function $\psi$ satisfies 
\[-2 \, \Delta \psi - \frac{N-1}{N-2} \, |\nabla \psi|^2 + 2K + N(N-1) = 0\] 
at each point in $M$. 
\end{proposition} 

\textbf{Proof.} 
Using Proposition \ref{initial.value.of.J} and Proposition \ref{consequence.of.monotonicity}, we obtain 
\[\int_{M^{(j)}} G(l^{(j)}-w^{(j)}) \, \Big ( -2 \, \Delta \psi - \frac{N-1}{N-2} \, |\nabla \psi|^2 + 2K + N(N-1) \Big ) \to 0\] 
as $j \to \infty$. From this, the assertion follows. \\

\begin{proposition}
\label{z}
The function $z := w^2$ is a smooth function on $M$ and $|\nabla z|^2 = 4z$ at each point in $M$. For each $s>0$, the level set $\{z=s^2\}$ is a smooth submanifold of dimension $1$ which is diffeomorphic to $S^1$.
\end{proposition}

\textbf{Proof.} 
Since $\psi$ is smooth, Proposition \ref{relation.between.w.and.psi} implies that the function $\cosh(Nw)$ is smooth. From this, we deduce that the function $w^2$ is smooth. Moreover, it follows from Proposition \ref{relation.between.w.and.psi} that $|\nabla w|^2 = 1$ in the region $\{w>0\}$. This implies $|\nabla z|^2 = 4z$ in the region $\{w>0\}$. By Lemma \ref{zero.level.set.has.empty.interior}, the set $\{w>0\}$ is dense. Since $z$ is a smooth function, it follows that $|\nabla z|^2 = 4z$ at each point in $M$. 

For each $s>0$, the level set $\{z=s^2\}$ is a smooth submanifold of dimension $1$. Finally, it follows from Lemma \ref{level.sets.of.psi} and Proposition \ref{relation.between.w.and.psi} that the level set $\{z=s^2\}$ is diffeomorphic to $S^1$ if $s$ is sufficiently large. By Morse theory, the level set $\{z=s^2\}$ is diffeomorphic to $S^1$ for each $s>0$. This completes the proof of Proposition \ref{z}. \\

For each $s \in [0,\infty)$, we define $\Gamma_s = \{w=s\} = \{z=s^2\}$. In the next step, we define a one-parameter family of smooth maps $\varphi_s: \Gamma_1 \to M$, $s \in (0,\infty)$. For each point $x \in \Gamma_1$, we define the path $\{\varphi_s(x): s>0\}$ to be the solution of the ODE 
\[\frac{\partial}{\partial s} \varphi_s(x) = \nabla w \big |_{\varphi_s(x)}\] 
with initial condition $\varphi_1(x) = x$. It is easy to see that the solution of this ODE is defined for all $s \in (0,\infty)$. Note that $\varphi_s(\Gamma_1) = \Gamma_s$ for each $s>0$. 

\begin{lemma}
\label{normal.exponential.map}
We have $\varphi_s(x) = \exp_x((s-1) \, \nabla w(x))$ for each point $x \in \Gamma_1$ and each $s>0$.
\end{lemma}

\textbf{Proof.} 
Recall that $|\nabla w|^2 = 1$ in the region $\{w>0\}$. Differentiating this identity gives $D_{\nabla w} \nabla w = 0$ in the region $\{w>0\}$. Consequently, for each point $x \in \Gamma_1$, the path $\{\varphi_s(x): s>0\}$ is a geodesic. This completes the proof of Lemma \ref{normal.exponential.map}. \\

We define a smooth map $\varphi_0: \Gamma_1 \to M$ by $\varphi_0(x) = \exp_x(-\nabla w(x))$ for each point $x \in \Gamma_1$.

\begin{lemma}
\label{Gamma_0}
We have $\varphi_0(\Gamma_1) = \Gamma_0$. 
\end{lemma}

\textbf{Proof.} 
We first consider a point $x \in \Gamma_1$. Then $\varphi_s(x) \in \Gamma_s$ for each $s>0$. Sending $s \searrow 0$, it follows that $\varphi_0(x) \in \Gamma_0$. Thus, $\varphi_0(\Gamma_1) \subset \Gamma_0$.

We next consider an arbitrary point $y \in \Gamma_0$. By Lemma \ref{zero.level.set.has.empty.interior}, we can find a sequence of positive real numbers $s_j \to 0$ and a sequence of points $y^{(j)} \in \Gamma_{s_j}$ such that $y^{(j)} \to y$. For each $j$, we can find a point $x^{(j)} \in \Gamma_1$ such that $\varphi_{s_j}(x^{(j)}) = y^{(j)}$. After passing to a subsequence, we may assume that the sequence $x^{(j)}$ converges to a point $x \in \Gamma_1$. Then $\varphi_0(x) = \lim_{j \to \infty} \varphi_{s_j}(x^{(j)}) = \lim_{j \to \infty} y^{(j)} = y$. Thus, $\Gamma_0 \subset \varphi_0(\Gamma_1)$. This completes the proof of Lemma \ref{Gamma_0}. \\

For each $s>0$, we denote by $\kappa$ the geodesic curvature of $\Gamma_s$. Since $|\nabla w| = 1$, we know that $\Delta w = \kappa$ at each point on $\Gamma_s$.

\begin{lemma}
\label{ode.for.kappa}
We have 
\begin{align*} 
\frac{\partial}{\partial s} \kappa(\varphi_s(x)) 
&= -\kappa(\varphi_s(x))^2 - \frac{(N-2) \, \sinh(Ns)}{1+\cosh(Ns)} \, \kappa(\varphi_s(x)) \\ 
&- \frac{N-2}{1+\cosh(Ns)} + (N-1) 
\end{align*} 
for each point $x \in \Gamma_1$ and each $s>0$. 
\end{lemma}

\textbf{Proof.} 
Proposition \ref{relation.between.w.and.psi} implies that 
\[|\nabla \psi|^2 = (N-2)^2 \, \Big ( 1 - \frac{2}{1+\cosh(Ns)} \Big )\] 
at each point on $\Gamma_s$. Moreover, using Proposition \ref{relation.between.w.and.psi}, we obtain 
\begin{align*} 
\Delta \psi 
&= (N-2) \, \Big ( \frac{\sinh(Ns)}{1+\cosh(Ns)} \, \Delta w + \frac{N}{1+\cosh(Ns)} \, |\nabla w|^2 \Big ) \\ 
&= (N-2) \, \Big ( \frac{\sinh(Ns)}{1+\cosh(Ns)} \, \kappa + \frac{N}{1+\cosh(Ns)} \Big ) 
\end{align*}
at each point on $\Gamma_s$. Using these identities together with Proposition \ref{pde.for.psi}, we conclude that 
\begin{align*} 
0 &= -2 \, \Delta \psi - \frac{N-1}{N-2} \, |\nabla \psi|^2 + 2K + N(N-1) \\ 
&= -\frac{2(N-2) \, \sinh(Ns)}{1+\cosh(Ns)} \, \kappa - \frac{2(N-2)}{1+\cosh(Ns)} + 2K + 2(N-1) 
\end{align*}
at each point on $\Gamma_s$. The assertion now follows from the fact that 
\[\frac{\partial}{\partial s} \kappa(\varphi_s(x)) = -\kappa(\varphi_s(x))^2 - K(\varphi_s(x))\] 
for each point $x \in \Gamma_1$ and each $s>0$. This completes the proof of Lemma \ref{ode.for.kappa}. \\

\begin{lemma} 
\label{first.case}
Suppose that $x$ is a point in $\Gamma_1$ with $\kappa(x) \neq \frac{\sinh N}{1+\cosh N} + \frac{N}{\sinh N}$. Let us write $\kappa(x) = \frac{\sinh N}{1+\cosh N} + \frac{Na}{1+\cosh N + a \sinh N}$ for some real number $a \neq -\frac{1+\cosh N}{\sinh N}$. Then $a \geq -1$ and 
\[\kappa(\varphi_s(x)) = \frac{\sinh(Ns)}{1+\cosh(Ns)} + \frac{Na}{1+\cosh(Ns) + a \sinh(Ns)}\] 
for all $s>0$. Moreover, the differential $(D\varphi_0)_x: T_x \Gamma_1 \to T_{\varphi_0(x)} M$ is non-zero.
\end{lemma}

\textbf{Proof.} 
Note that $1+\cosh N + a \sinh N \neq 0$ by definition of $a$. Let $\mathcal{S}$ denote the connected component of the set $\{s \in (0,\infty): 1+\cosh(Ns) + a \sinh(Ns) \neq 0\}$ with $1 \in \mathcal{S}$. Using Lemma \ref{ode.for.kappa} together with standard uniqueness results for ODE, we obtain 
\[\kappa(\varphi_s(x)) = \frac{\sinh(Ns)}{1+\cosh(Ns)} + \frac{Na}{1+\cosh(Ns) + a \sinh(Ns)}\] 
for all $s \in \mathcal{S}$. Since the function $s \mapsto \kappa(\varphi_s(x))$ is a smooth function defined for all $s \in (0,\infty)$, it follows that $\mathcal{S} = (0,\infty)$ and $a \geq -1$. This proves the first statement. The second statement follows from the fact that $\kappa(\varphi_s(x))$ is bounded as $s \searrow 0$. \\

\begin{lemma} 
\label{second.case}
Suppose that $x$ is a point in $\Gamma_1$ with $\kappa(x) = \frac{\sinh N}{1+\cosh N} + \frac{N}{\sinh N}$. Then 
\[\kappa(\varphi_s(x)) = \frac{\sinh(Ns)}{1+\cosh(Ns)} + \frac{N}{\sinh(Ns)}\] 
for all $s>0$. Moreover, the differential $(D\varphi_0)_x: T_x \Gamma_1 \to T_{\varphi_0(x)} M$ vanishes. 
\end{lemma} 

\textbf{Proof.} 
The first statement again follows from Lemma \ref{ode.for.kappa} together with standard uniqueness results for ODE. The second statement follows from the fact that $\kappa(\varphi_s(x)) - \frac{1}{s}$ is bounded as $s \searrow 0$. \\

\subsection{The case when $\kappa(x) \neq \frac{\sinh N}{1+\cosh N} + \frac{N}{\sinh N}$ for each point $x \in \Gamma_1$} 

In this case, Lemma \ref{first.case} implies that the map $\varphi_0: \Gamma_1 \to M$ is a smooth immersion. 

\begin{lemma}
\label{Hessian.of.z.case.1}
For each point $y \in \Gamma_0$, the Hessian of the function $z$ has eigenvalues $2$ and $0$.
\end{lemma}

\textbf{Proof.} 
Let us fix an arbitrary point $y \in \Gamma_0$. By Lemma \ref{Gamma_0}, we can find a point $x \in \Gamma_1$ such that $y = \varphi_0(x)$. For each $s>0$, the Hessian of the function $w$ at the point $\varphi_s(x)$ has eigenvalues $0$ and $\kappa(\varphi_s(x))$. Moreover, for each $s>0$, the Hessian of the function $z$ at the point $\varphi_s(x)$ has eigenvalues $2$ and $2s \, \kappa(\varphi_s(x))$. By Lemma \ref{first.case}, $2s \, \kappa(\varphi_s(x)) \to 0$ as $s \searrow 0$. Therefore, the Hessian of the function $z$ at the point $\varphi_0(x)$ has eigenvalues $2$ and $0$. This completes the proof of Lemma \ref{Hessian.of.z.case.1}. \\

In the next step, we show that the function $z$ is a Morse-Bott function (see \cite{Banyaga-Hurtubise} for a definition). 

\begin{lemma}
\label{morse.bott.function} 
The set $\Gamma_0$ is a smooth submanifold of dimension $1$ and the function $z$ is a Morse-Bott function. 
\end{lemma} 

\textbf{Proof.} 
Let us fix an arbitrary point $y \in \Gamma_0$. Let $\xi_1,\xi_2 \in T_y M$ denote the eigenvectors of the Hessian of the function $z$ at the point $y$. We assume that $\xi_1$ is an eigenvector with eigenvalue $2$, and $\xi_2$ is an eigenvector with eigenvalue $0$. We extend $\xi_1$ and $\xi_2$ to smooth vector fields on $M$. By the implicit function theorem, we can find an open neighborhood $U$ of $y$ with the property that the set $U \cap \{\langle \nabla z,\xi_1 \rangle = 0\}$ is contained in a smooth submanifold $Z$ of dimension $1$. This implies $U \cap \Gamma_0 \subset U \cap \{\nabla z=0\} \subset Z$. On the other hand, we know that $\Gamma_1$ is a smooth submanifold and $\Gamma_0$ is the image of $\Gamma_1$ under a smooth immersion. Consequently, we can find a smooth curve $\gamma: [-1,1] \to M$ such that $\gamma(0) = y$, $\gamma'(0) \neq 0$, and $\gamma([-1,1]) \subset \Gamma_0$. Thus, we can find an open neighborhood $\tilde{U}$ of $y$ such that $\tilde{U} \subset U$ and $\tilde{U} \cap \Gamma_0 = \tilde{U} \cap Z$. This shows that $\Gamma_0$ is a smooth submanifold of dimension $1$. In view of Lemma \ref{Hessian.of.z.case.1}, it follows that the function $z$ is a Morse-Bott function. This completes the proof of Lemma \ref{morse.bott.function}. \\

Since $M$ is orientable, the submanifold $\Gamma_0$ is two-sided. It follows from Lemma \ref{morse.bott.function} that the set $\{z \leq s^2\}$ is diffeomorphic to $\Gamma_0 \times [-1,1]$ if $s>0$ is sufficiently small. In particular, if $s>0$ is sufficiently small, then the boundary $\{z=s^2\}$ is disconnected. This contradicts Proposition \ref{z}. Therefore, this case cannot occur.

\subsection{The case when $\kappa(x) = \frac{\sinh N}{1+\cosh N} + \frac{N}{\sinh N}$ for some point $x \in \Gamma_1$}

In this case, we fix a point $x_0 \in \Gamma_1$ such that $\kappa(x_0) = \frac{\sinh N}{1+\cosh N} + \frac{N}{\sinh N}$. Let $y_0 = \varphi_0(x_0) \in \Gamma_0$.

\begin{lemma}
\label{Hessian.of.z.case.2}
The Hessian of the function $z$ at the point $y_0$ has a single eigenvalue $2$ of multiplicity $2$. 
\end{lemma} 

\textbf{Proof.} 
For each $s>0$, the Hessian of the function $z$ at the point $\varphi_s(x_0)$ has eigenvalues $2$ and $2s \, \kappa(\varphi_s(x_0))$. By Lemma \ref{second.case}, $2s \, \kappa(\varphi_s(x_0)) \to 2$ as $s \searrow 0$. Therefore, the Hessian of the function $z$ at the point $\varphi_0(x_0)$ has a single eigenvalue $2$ of multiplicity $2$. This completes the proof of Lemma \ref{Hessian.of.z.case.2}. \\

In view of Lemma \ref{Hessian.of.z.case.2}, we can find an open neighborhood $U$ of $y_0$ such that $U \cap \Gamma_0 = \{y_0\}$. Consequently, the set $\{x \in \Gamma_1: \varphi_0(x) = y_0\}$ is both open and closed as a subset of $\Gamma_1$. Moreover, the set $\{x \in \Gamma_1: \varphi_0(x) = y_0\}$ contains the point $x_0$. Since $\Gamma_1$ is connected by Proposition \ref{z}, we conclude that $\varphi_0(x)=y_0$ for each point $x \in \Gamma_1$. Since $\varphi_0(\Gamma_1) = \Gamma_0$ by Lemma \ref{Gamma_0}, it follows that the set $\Gamma_0$ consists of a single point. Moreover, for each point $x \in \Gamma_1$, the differential $(D\varphi_0)_x: T_x \Gamma_1 \to T_{\varphi_0(x)} M$ vanishes. Using Lemma \ref{first.case}, it follows that $\kappa(x) = \frac{\sinh N}{1+\cosh N} + \frac{N}{\sinh N}$ for each point $x \in \Gamma_1$. Using Lemma \ref{second.case}, we deduce that 
\[\kappa(\varphi_s(x)) = \frac{\sinh(Ns)}{1+\cosh(Ns)} + \frac{N}{\sinh(Ns)}\] 
for each point $x \in \Gamma_1$ and each $s>0$. Therefore, if we define 
\begin{align*} 
L(s) 
&= \frac{1}{2} \, \Big [ \frac{1+\cosh(Ns)}{2} \Big ]^{-\frac{N-1}{N}} \, \sinh(Ns) \\ 
&= \Big [ \cosh \Big ( \frac{Ns}{2} \Big ) \Big ]^{-\frac{N-2}{N}} \, \sinh \Big ( \frac{Ns}{2} \Big ), 
\end{align*}
then $\kappa(\varphi_s(x)) = \frac{d}{ds} \log L(s)$ for each $x \in \Gamma_1$ and each $s>0$. From this, it is easy to see that $(M,g)$ is locally isometric to $(\mathbb{R}^2,g_{\text{\rm HM},N,2})$. By Proposition \ref{relation.between.w.and.psi}, the function 
\[\frac{N-2}{N} \, \log \Big ( \frac{1+\cosh(Nw)}{2} \Big ) - \log \rho = \frac{2(N-2)}{N} \, \log \cosh \Big ( \frac{Nw}{2} \Big ) - \log \rho\] 
is constant. From this, we deduce that $(M,g,\rho)$ is a model $(N,2)$-dataset.

\section{Properties of $(g,\rho)$-stationary hypersurfaces which are $(g,\rho,u)$-stable in the sense of Definition \ref{definition.stability}}

\label{properties.of.stable.hypersurfaces}

Throughout this section, we assume that $N$ and $n$ are integers satisfying $3 \leq n \leq N$ and $(M,g,\rho)$ is an $(N,n)$-dataset. Let us fix a function $u: T^{n-1} \to \mathbb{R}$ such that 
\[\Delta_\gamma u + \frac{N}{2} \, \text{\rm tr}_\gamma(Q) + N \, P = \text{\rm constant}.\] 
The function $u$ is twice continuously differentiable with H\"older continuous second derivatives. Note that 
\[\int_{T^{n-1}} \Big ( N \, \text{\rm tr}_\gamma(Q) + 2N \, P + \Big ( \frac{2}{Nb_0} \Big )^N \Big ) \, d\text{\rm vol}_\gamma \leq 0\] 
by definition of an $(N,n)$-dataset. This implies 
\begin{equation} 
\label{pde.for.u}
\Delta_\gamma u + \frac{N}{2} \, \text{\rm tr}_\gamma(Q) + N \, P + \frac{1}{2} \, \Big ( \frac{2}{Nb_0} \Big )^N \leq 0. 
\end{equation}
at each point on $T^{n-1}$.

Throughout this section, we assume that $\Sigma$ is a properly embedded, connected, orientable hypersurface in $M$ which is $t_*$-tame for some $t_* \in S^1$ (see Definition \ref{definition.tame}). Let $r_*$ be chosen as in Definition \ref{definition.tame}. We further assume that $\Sigma$ is $(g,\rho)$-stationary and $(g,\rho,u)$-stable in the sense of Definition \ref{definition.stability}. We denote by $\mathbb{L}_\Sigma$ the weighted Jacobi operator of $\Sigma$ (see Definition \ref{definition.weighted.Jacobi.operator}).

\begin{definition}
\label{definition.hyperbolic.metric.on.Sigma}
Consider the map $\pi: \Sigma \cap \{r \geq r_*\} \to [r_*,\infty) \times T^{n-2}$ which maps a point $(r,\theta_0,\hdots,\theta_{n-3},\theta_{n-2}) \in \Sigma \cap \{r \geq r_*\}$ to $(r,\theta_0,\hdots,\theta_{n-3}) \in [r_*,\infty) \times T^{n-2}$. We denote by $g_{\text{\rm hyp}}$ the pull-back of the hyperbolic metric $r^{-2} \, dr \otimes dr + \sum_{k=0}^{n-3} b_k^2 \, r^2 \, d\theta_k \otimes d\theta_k$ on $[r_*,\infty) \times T^{n-2}$ under the map $\pi$. Note that $g_{\text{\rm hyp}}$ is a hyperbolic metric on $\Sigma \cap \{r \geq r_*\}$. The metric $g_{\text{\rm hyp}}$ is obtained by restricting the $(0,2)$-tensor $\bar{g} - b_{n-2}^2 \, r^2 \, d\theta_{n-2} \otimes d\theta_{n-2}$ on $[r_*,\infty) \times T^{n-1}$ to $\Sigma \cap \{r \geq r_*\}$.
\end{definition}

In the following, we assume that the unit normal vector field along $\Sigma$ is chosen so that $\langle \frac{\partial}{\partial \theta_{n-2}},\nu_\Sigma \rangle > 0$ outside a compact set. Moreover, we fix a positive smooth function $\bar{v}: \Sigma \to \mathbb{R}$ with the property that $\bar{v} = \langle \frac{\partial}{\partial \theta_{n-2}},\nu_\Sigma \rangle$ outside a compact set. 

\begin{lemma} 
\label{higher.derivative.bound.for.bar.v}
Let $m$ be a nonnegative integer. Then $|D_{\text{\rm hyp}}^{\Sigma,m} \bar{v}|_{g_{\text{\rm hyp}}} \leq O(r)$, where $D_{\text{\rm hyp}}^{\Sigma,m}$ denotes the covariant derivative of order $m$ with respect to the metric $g_{\text{\rm hyp}}$.
\end{lemma} 

\textbf{Proof.} 
This follows directly from the assumption that $\Sigma$ is tame. \\

\begin{lemma}
\label{asymptotics.for.bar.v}
Let $m$ be a nonnegative integer. Then 
\[\Big | D_{\text{\rm hyp}}^{\Sigma,m} \Big ( \bar{v} - b_{n-2} \, r -\frac{1}{2} \, b_{n-2}^{-1} \, r^{1-N} \, Q(\frac{\partial}{\partial \theta_{n-2}},\frac{\partial}{\partial \theta_{n-2}}) \Big ) \Big |_{g_{\text{\rm hyp}}} \leq O(r^{1-N-\delta}),\] 
where $D_{\text{\rm hyp}}^{\Sigma,m}$ denotes the covariant derivative of order $m$ with respect to the metric $g_{\text{\rm hyp}}$.
\end{lemma} 

\textbf{Proof.} 
Since $\Sigma$ is tame, we obtain $|(\frac{\partial}{\partial \theta_{n-2}})^{\text{\rm tan}}| \leq O(r^{2-N})$ along $\Sigma$. This implies 
\[\Big | \langle \frac{\partial}{\partial \theta_{n-2}},\frac{\partial}{\partial \theta_{n-2}} \rangle - \bar{v}^2 \Big | = \Big | \Big ( \frac{\partial}{\partial \theta_{n-2}} \Big )^{\text{\rm tan}} \Big |^2 \leq O(r^{4-2N})\] 
outside a compact set. The asymptotic expansion of the metric $g$ gives 
\[\Big | \langle \frac{\partial}{\partial \theta_{n-2}},\frac{\partial}{\partial \theta_{n-2}} \rangle - b_{n-2}^2 \, r^2 - r^{2-N} \, Q(\frac{\partial}{\partial \theta_{n-2}},\frac{\partial}{\partial \theta_{n-2}}) \Big | \leq O(r^{2-N-2\delta}).\] 
Putting these facts together, we obtain 
\[\Big | \bar{v}^2 - b_{n-2}^2 \, r^2 - r^{2-N} \, Q(\frac{\partial}{\partial \theta_{n-2}},\frac{\partial}{\partial \theta_{n-2}}) \Big | \leq O(r^{2-N-2\delta}).\] 
Since $\bar{v}$ is a positive function, it follows that 
\[\Big | \bar{v} - b_{n-2} \, r -\frac{1}{2} \, b_{n-2}^{-1} \, r^{1-N} \, Q(\frac{\partial}{\partial \theta_{n-2}},\frac{\partial}{\partial \theta_{n-2}}) \Big | \leq O(r^{1-N-2\delta}).\] 
Finally, using Lemma \ref{higher.derivative.bound.for.bar.v}, we obtain 
\[\Big | D_{\text{\rm hyp}}^{\Sigma,m} \Big ( \bar{v} - b_{n-2} \, r -\frac{1}{2} \, b_{n-2}^{-1} \, r^{1-N} \, Q(\frac{\partial}{\partial \theta_{n-2}},\frac{\partial}{\partial \theta_{n-2}}) \Big ) \Big |_{g_{\text{\rm hyp}}} \leq O(r)\] 
for every nonnegative integer $m$. The assertion now follows from standard interpolation inequalities. \\

\begin{lemma}
\label{Jacobi.operator.of.bar.v}
We have $|\mathbb{L}_\Sigma \bar{v}| \leq O(r^{1-n-\delta})$. Moreover, $|D_{\text{\rm hyp}}^{\Sigma,m} \mathbb{L}_\Sigma \bar{v}|_{g_{\text{\rm hyp}}} \leq O(r^{1-n})$ for every nonnegative integer $m$.
\end{lemma}

\textbf{Proof.} 
Let $V$ be a smooth vector field on $M$ with the property that $V = \frac{\partial}{\partial \theta_{n-2}}$ outside a compact set. It follows from Proposition \ref{formula.weighted.Jacobi.operator} that 
\begin{align*} 
\mathbb{L}_\Sigma \bar{v} 
&= -\rho \sum_{k=1}^{n-1} (D_{e_k} (\mathscr{L}_V g))(e_k,\nu_\Sigma) + \frac{1}{2} \, \rho \sum_{k=1}^{n-1} (D_{\nu_\Sigma} (\mathscr{L}_V g))(e_k,e_k) \\ 
&- \rho \sum_{k,l=1}^{n-1} h_\Sigma(e_k,e_l) \, (\mathscr{L}_V g)(e_k,e_l) - (\mathscr{L}_V g)(\nabla \rho,\nu_\Sigma) + \rho \, \big \langle \nabla(V(\log \rho)),\nu_\Sigma \big \rangle 
\end{align*} 
outside a compact set. Using Lemma \ref{Lie.derivative.of.metric.along.angular.vector.field} and Lemma \ref{derivative.of.rho.along.angular.vector.field}, we obtain $|\mathbb{L}_\Sigma \bar{v}| \leq O(r^{1-n-\delta})$. On the other hand, Lemma \ref{higher.derivative.bound.for.bar.v} implies $|D_{\text{\rm hyp}}^{\Sigma,m} \mathbb{L}_\Sigma \bar{v}|_{g_{\text{\rm hyp}}} \leq O(r^{N-n+1})$ for every nonnegative integer $m$. Using standard interpolation inequalities, we conclude that $|D_{\text{\rm hyp}}^{\Sigma,m} \mathbb{L}_\Sigma \bar{v}|_{g_{\text{\rm hyp}}} \leq O(r^{1-n})$ for every nonnegative integer $m$. This completes the proof of Lemma \ref{Jacobi.operator.of.bar.v}. \\

\begin{lemma} 
\label{Jacobi.operator.of.test.function}
We have 
\[|\mathbb{L}_\Sigma(r^{-N} \, \bar{v})| \leq o(r^{1-n-\delta})\] 
and 
\[|\mathbb{L}_\Sigma(r^{-N-\delta} \, \bar{v}) + b_{n-2} \, \delta(N+\delta) \, r^{1-n-\delta}| \leq o(r^{1-n-\delta}).\] 
\end{lemma} 

\textbf{Proof.} 
We compute 
\[\mathbb{L}_\Sigma(r^{-N} \, \bar{v}) = r^{-N} \, \mathbb{L}_\Sigma \bar{v} - 2\rho \, \langle \nabla^\Sigma (r^{-N}),\nabla^\Sigma \bar{v} \rangle - \text{\rm div}_\Sigma(\rho \, \nabla^\Sigma (r^{-N})) \, \bar{v}\] 
and 
\[\mathbb{L}_\Sigma(r^{-N-\delta} \, \bar{v}) = r^{-N-\delta} \, \mathbb{L}_\Sigma \bar{v} - 2\rho \, \langle \nabla^\Sigma (r^{-N-\delta}),\nabla^\Sigma \bar{v} \rangle - \text{\rm div}_\Sigma(\rho \, \nabla^\Sigma (r^{-N-\delta})) \, \bar{v}.\] 
It is easy to see that 
\[|2\rho \, \langle \nabla^\Sigma (r^{-N}),\nabla^\Sigma \bar{v} \rangle + \text{\rm div}_\Sigma(\rho \, \nabla^\Sigma (r^{-N})) \, \bar{v}| \leq o(r^{1-n-\delta})\] 
and 
\begin{align*} 
&|2\rho \, \langle \nabla^\Sigma (r^{-N-\delta}),\nabla^\Sigma \bar{v} \rangle + \text{\rm div}_\Sigma(\rho \, \nabla^\Sigma (r^{-N-\delta})) \, \bar{v} \\ 
&- b_{n-2} \, \delta(N+\delta) \, r^{1-n-\delta}| \leq o(r^{1-n-\delta}). 
\end{align*} 
The assertion now follows from Lemma \ref{Jacobi.operator.of.bar.v}. This completes the proof of Lemma \ref{Jacobi.operator.of.test.function}. \\

In the following, we consider a sequence $r_j \to \infty$. For each $j$, we define $\Sigma^{(j)} = \Sigma \setminus \{r>r_j\}$. Since $\Sigma$ is connected and tame, it follows that $\Sigma^{(j)}$ is connected if $j$ is sufficiently large. For each $j$, we define 
\begin{align} 
\label{definition.Lambda_j}
\Lambda_j 
&= -\int_{\partial \Sigma^{(j)}} \rho \, \Big \langle \Big ( D_{\frac{\partial}{\partial \theta_{n-2}}} \frac{\partial}{\partial \theta_{n-2}} \Big )^{\text{\rm tan}},\eta \Big \rangle \notag \\ 
&- \int_{T^{n-2} \times \{t_*\}} \frac{\partial^2 u}{\partial \theta_{n-2}^2} \, d\text{\rm vol}_\gamma - r_j^{-\frac{\delta}{2}}. 
\end{align} 
Here, $\eta = \frac{\nabla^\Sigma r}{|\nabla^\Sigma r|}$ denotes the outward-pointing unit normal vector to $\partial \Sigma^{(j)}$ in $\Sigma$. It follows from Lemma \ref{covariant.derivative.angular.vector.field} that the sequence $r_j^{-N} \, \Lambda_j$ converges to a positive real number as $j \to \infty$.

\begin{lemma} 
\label{asymptotics.for.Lambda_j} 
We have 
\[\Lambda_j - \int_{\partial \Sigma^{(j)}} \rho \, \bar{v} \, \langle \nabla^\Sigma \bar{v},\eta \rangle \to -\int_{T^{n-2} \times \{t_*\}} \frac{\partial^2 u}{\partial \theta_{n-2}^2} \, d\text{\rm vol}_\gamma\] 
as $j \to \infty$. 
\end{lemma}

\textbf{Proof.} 
Using Lemma \ref{covariant.derivative.angular.vector.field}, we obtain 
\begin{align} 
\label{boundary.integral.1}
&\int_{\partial \Sigma^{(j)}} \rho \, \Big \langle \Big ( D_{\frac{\partial}{\partial \theta_{n-2}}} \frac{\partial}{\partial \theta_{n-2}} + b_{n-2}^2 \, r \, \nabla r \Big )^{\text{\rm tan}},\eta \Big \rangle \notag \\ 
&\to \int_{T^{n-2} \times \{t_*\}} \frac{N-2}{2} \, Q(\frac{\partial}{\partial \theta_{n-2}},\frac{\partial}{\partial \theta_{n-2}}) \, d\text{\rm vol}_\gamma 
\end{align} 
as $j \to \infty$. Moreover, Lemma \ref{asymptotics.for.bar.v} implies 
\begin{align} 
\label{boundary.integral.2}
&\int_{\partial \Sigma^{(j)}} \rho \, (\bar{v} \, \langle \nabla^\Sigma \bar{v},\eta \rangle - b_{n-2}^2 \, r \, \langle \nabla^\Sigma r,\eta \rangle) \notag \\ 
&\to -\int_{T^{n-2} \times \{t_*\}} \frac{N-2}{2} \, Q(\frac{\partial}{\partial \theta_{n-2}},\frac{\partial}{\partial \theta_{n-2}}) \, d\text{\rm vol}_\gamma 
\end{align} 
as $j \to \infty$. Adding (\ref{boundary.integral.1}) and (\ref{boundary.integral.2}), we conclude that 
\begin{equation} 
\label{sum.of.boundary.integrals}
\int_{\partial \Sigma^{(j)}} \rho \, \Big \langle \Big ( D_{\frac{\partial}{\partial \theta_{n-2}}} \frac{\partial}{\partial \theta_{n-2}} \Big )^{\text{\rm tan}},\eta \Big \rangle + \int_{\partial \Sigma^{(j)}} \rho \, \bar{v} \, \langle \nabla^\Sigma \bar{v},\eta \rangle \to 0 
\end{equation}
as $j \to \infty$. The assertion follows by combining (\ref{definition.Lambda_j}) and (\ref{sum.of.boundary.integrals}). This completes the proof of Lemma \ref{asymptotics.for.Lambda_j}. \\

\begin{proposition}
\label{consequence.of.stability.inequality.1} 
Let $a$ be a real number and let $V$ be a smooth vector field on $M$ with the property that $V = a \, \frac{\partial}{\partial \theta_{n-2}}$ in a neighborhood of the set $\{r=r_j\}$. Then 
\begin{align*} 
&\frac{1}{2} \int_{\Sigma^{(j)}} \rho \, \sum_{k=1}^{n-1} (\mathscr{L}_V \mathscr{L}_V g)(e_k,e_k) + \int_{\Sigma^{(j)}} V(V(\rho)) \\ 
&- \frac{1}{2} \int_{\Sigma^{(j)}} \rho \, \sum_{k,l=1}^{n-1} (\mathscr{L}_V g)(e_k,e_l) \, (\mathscr{L}_V g)(e_k,e_l) \\ 
&+ \frac{1}{4} \int_{\Sigma^{(j)}} \rho \, \sum_{k,l=1}^{n-1} (\mathscr{L}_V g)(e_k,e_k) \, (\mathscr{L}_V g)(e_l,e_l) \\ 
&+ \int_{\Sigma^{(j)}} V(\rho) \, \sum_{k=1}^{n-1} (\mathscr{L}_V g)(e_k,e_k) \\ 
&\geq -a^2 \int_{T^{n-2} \times \{t_*\}} \frac{\partial^2 u}{\partial \theta_{n-2}^2} \, d\text{\rm vol}_\gamma - C \, a^2 \, r_j^{-\delta}. 
\end{align*} 
Here, $\{e_1,\hdots,e_{n-1}\}$ denotes a local orthonormal frame on $\Sigma$, and $C$ is a positive constant which is independent of $j$.
\end{proposition}

\textbf{Proof.} 
We may assume that $V = a \, \frac{\partial}{\partial \theta_{n-2}}$ in the region $\{r>r_j\}$. Using Lemma \ref{Lie.derivative.of.metric.along.angular.vector.field}, we obtain 
\[|\mathscr{L}_V g| \leq C \, |a| \, r^{1-N-\delta}, \quad |\mathscr{L}_V \mathscr{L}_V g| \leq C \, |a|^2 \, r^{2-N-\delta}\] 
in the region $\{r>r_j\}$. Moreover, Lemma \ref{derivative.of.rho.along.angular.vector.field} gives 
\[|V(\rho)| \leq C \, |a| \, r^{1-n-\delta}, \quad |V(V(\rho))| \leq C \, |a|^2 \, r^{2-n-\delta}\] 
in the region $\{r>r_j\}$. Putting these facts together, we conclude that 
\begin{align*} 
&\frac{1}{2} \int_{\Sigma \cap \{r>r_j\}} \rho \, \sum_{k=1}^{n-1} (\mathscr{L}_V \mathscr{L}_V g)(e_k,e_k) + \int_{\Sigma \cap \{r>r_j\}} V(V(\rho)) \\ 
&- \frac{1}{2} \int_{\Sigma \cap \{r>r_j\}} \rho \, \sum_{k,l=1}^{n-1} (\mathscr{L}_V g)(e_k,e_l) \, (\mathscr{L}_V g)(e_k,e_l) \\ 
&+ \frac{1}{4} \int_{\Sigma \cap \{r>r_j\}} \rho \, \sum_{k,l=1}^{n-1} (\mathscr{L}_V g)(e_k,e_k) \, (\mathscr{L}_V g)(e_l,e_l) \\ 
&+ \int_{\Sigma \cap \{r>r_j\}} V(\rho) \, \sum_{k=1}^{n-1} (\mathscr{L}_V g)(e_k,e_k) \\ 
&\leq C \, a^2 \, r_j^{-\delta}. 
\end{align*} 
On the other hand, $\Sigma$ is $(g,\rho,u)$-stable in the sense of Definition \ref{definition.stability}. Putting these facts together, the assertion follows. \\

\begin{proposition} 
\label{consequence.of.stability.inequality.2}
If $j$ is sufficiently large, then the following holds. Let $a$ be a real number and let $v$ be a smooth function on $\Sigma^{(j)}$ such that $v = a \, \bar{v}$ in a neighborhood of the boundary $\partial \Sigma^{(j)}$. Then  
\begin{align*} 
&-\Lambda_j \, a^2 + \int_{\Sigma^{(j)}} \rho \, |\nabla^\Sigma v|^2 - \int_{\Sigma^{(j)}} \rho \, (\text{\rm Ric}(\nu_\Sigma,\nu_\Sigma) + |h_\Sigma|^2) \, v^2 \\ 
&+ \int_{\Sigma^{(j)}} (D^2 \rho)(\nu_\Sigma,\nu_\Sigma) \, v^2 - \int_{\Sigma^{(j)}} \rho^{-1} \, \langle \nabla \rho,\nu_\Sigma \rangle^2 \, v^2 \geq 0. 
\end{align*} 
\end{proposition} 

\textbf{Proof.} 
We can find a smooth vector field $V$ on $M$ such that $v = \langle V,\nu_\Sigma \rangle$ at each point on $\Sigma^{(j)}$ and $V = a \, \frac{\partial}{\partial \theta_{n-2}}$ in a neighborhood of the set $\{r=r_j\}$. Let $W = D_V V$. Then $W = a^2 \, D_{\frac{\partial}{\partial \theta_{n-2}}} \frac{\partial}{\partial \theta_{n-2}}$ in a neighborhood of the set $\{r=r_j\}$. We define a tangential vector field $Z$ along $\Sigma$ by 
\[Z = D_{V^{\text{\rm tan}}}^\Sigma (V^{\text{\rm tan}}) - \text{\rm div}_\Sigma(V^{\text{\rm tan}}) \, V^{\text{\rm tan}} + 2 \sum_{k=1}^{n-1} h_\Sigma(V^{\text{\rm tan}},e_k) \, \langle V,\nu_\Sigma \rangle \, e_k,\] 
where $\{e_1,\hdots,e_{n-1}\}$ denotes a local orthonormal frame on $\Sigma$. Using Proposition \ref{consequence.of.stability.inequality.1} and Proposition \ref{formula.second.variation}, we obtain 
\begin{align} 
\label{stability.inequality.for.v}
&\int_{\Sigma^{(j)}} \rho \, |\nabla^\Sigma v|^2 - \int_{\Sigma^{(j)}} \rho \, (\text{\rm Ric}(\nu_\Sigma,\nu_\Sigma) + |h_\Sigma|^2) \, v^2 \notag \\ 
&+ \int_{\Sigma^{(j)}} (D^2 \rho)(\nu_\Sigma,\nu_\Sigma) \, v^2 - \int_{\Sigma^{(j)}} \rho^{-1} \, \langle \nabla \rho,\nu_\Sigma \rangle^2 \, v^2 \notag \\ 
&+ \int_{\Sigma^{(j)}} \big ( \text{\rm div}_\Sigma(\rho \, W^{\text{\rm tan}}) - \text{\rm div}_\Sigma(\rho \, Z) + \text{\rm div}_\Sigma(\langle V^{\text{\rm tan}},\nabla^\Sigma \rho \rangle \, V^{\text{\rm tan}}) \big ) \\ 
&\geq -a^2 \int_{T^{n-2} \times \{t_*\}} \frac{\partial^2 u}{\partial \theta_{n-2}^2} \, d\text{\rm vol}_\gamma - C \, a^2 \, r_j^{-\delta}. \notag 
\end{align}
By assumption, $\Sigma$ is tame and $V = a \, \frac{\partial}{\partial \theta_{n-2}}$ in a neighborhood of the set $\{r=r_j\}$. This implies $|h_\Sigma| \leq C \, r_j^{1-N}$, $|V| \leq C \, |a| \, r_j$, and $|V^{\text{\rm tan}}| + |D^\Sigma(V^{\text{\rm tan}})| \leq C \, |a| \, r_j^{2-N}$ at each point on $\partial \Sigma^{(j)}$. Thus, $|\rho \, Z| \leq C \, a^2 \, r_j^{4-N-n}$ and $|\langle V^{\text{\rm tan}},\nabla^\Sigma \rho \rangle \, V^{\text{\rm tan}}| \leq C \, a^2 \, r_j^{4-N-n}$ at each point on $\partial \Sigma^{(j)}$. Using the divergence theorem, we obtain 
\begin{align} 
\label{boundary.term}
&\int_{\Sigma^{(j)}} \big ( \text{\rm div}_\Sigma(\rho \, W^{\text{\rm tan}}) - \text{\rm div}_\Sigma(\rho \, Z) + \text{\rm div}_\Sigma(\langle V^{\text{\rm tan}},\nabla^\Sigma \rho \rangle \, V^{\text{\rm tan}}) \big ) \notag \\ 
&= \int_{\partial \Sigma^{(j)}} \big ( \rho \, \langle W^{\text{\rm tan}},\eta \rangle - \rho \, \langle Z,\eta \rangle + \langle V^{\text{\rm tan}},\nabla^\Sigma \rho \rangle \, \langle V^{\text{\rm tan}},\eta \rangle \big ) \\ 
&\leq a^2 \int_{\partial \Sigma^{(j)}} \rho \, \Big \langle \Big ( D_{\frac{\partial}{\partial \theta_{n-2}}} \frac{\partial}{\partial \theta_{n-2}} \Big )^{\text{\rm tan}},\eta \Big \rangle + C \, a^2 \, r_j^{2-N}. \notag
\end{align} 
The assertion follows by combining (\ref{definition.Lambda_j}), (\ref{stability.inequality.for.v}), and (\ref{boundary.term}). This completes the proof of Proposition \ref{consequence.of.stability.inequality.2}. \\

In the following, we assume that $j$ is chosen sufficiently large so that $r_j > 8r_*$ and the conclusion of Proposition \ref{consequence.of.stability.inequality.2} holds. Let us fix a nonnegative smooth function $\mu: \mathbb{R} \to \mathbb{R}$ which is supported in the interval $(2r_*,5r_*)$ and is strictly positive on the interval $[3r_*,4r_*]$. By scaling, we may assume that $\int_{\Sigma \cap \{2r_* < r < 5r_*\}} \mu(r)^3 = 1$. 

\begin{definition} 
\label{definition.lambda_j}
For each $j$, we denote by $\lambda_j$ the infimum of the functional 
\begin{align*} 
&-\Lambda_j \, a^2 + \int_{\Sigma^{(j)}} \rho \, |\nabla^\Sigma v|^2 - \int_{\Sigma^{(j)}} \rho \, (\text{\rm Ric}(\nu_\Sigma,\nu_\Sigma) + |h_\Sigma|^2) \, v^2 \\ 
&+ \int_{\Sigma^{(j)}} (D^2 \rho)(\nu_\Sigma,\nu_\Sigma) \, v^2 - \int_{\Sigma^{(j)}} \rho^{-1} \, \langle \nabla \rho,\nu_\Sigma \rangle^2 \, v^2 
\end{align*} 
over all pairs $(v,a) \in H^1(\Sigma^{(j)}) \times \mathbb{R}$ with the property that $v - a \, \bar{v} \in H_0^1(\Sigma^{(j)})$ and $\int_{\Sigma \cap \{2r_* < r < 5r_*\}} \mu(r) \, v^2 = 1$. 
\end{definition}

\begin{proposition}
\label{bounds.for.lambda_j}
If $j$ is sufficiently large, then $0 \leq \lambda_j \leq C$. Here, $C$ is a positive constant that does not depend on $j$.
\end{proposition}

\textbf{Proof.} 
It follows from Proposition \ref{consequence.of.stability.inequality.2} that $\lambda_j$ is nonnegative. To show that $\lambda_j$ is bounded from above, we recall that $\int_{\Sigma \cap \{2r_* < r < 5r_*\}} \mu(r)^3 = 1$. Using the function $\mu(r)$ as a test function in Definition \ref{definition.lambda_j}, we conclude that $\lambda_j \leq C$, where $C$ is independent of $j$. This completes the proof of Proposition \ref{bounds.for.lambda_j}. \\

After passing to a subsequence, we may assume that the sequence $\lambda_j$ converges to a nonnegative real number $\lambda_\infty$.

\begin{proposition} 
\label{existence.of.eigenfunction}
For each $j$, we can find a pair $(v^{(j)},a^{(j)}) \in H^1(\Sigma^{(j)}) \times \mathbb{R}$ such that $v^{(j)} - a^{(j)} \, \bar{v} \in H_0^1(\Sigma^{(j)})$, $\int_{\Sigma \cap \{2r_* < r < 5r_*\}} \mu(r) \, (v^{(j)})^2 = 1$, and  
\begin{align*} 
&-\Lambda_j \, (a^{(j)})^2 + \int_{\Sigma^{(j)}} \rho \, |\nabla^\Sigma v^{(j)}|^2 - \int_{\Sigma^{(j)}} \rho \, (\text{\rm Ric}(\nu_\Sigma,\nu_\Sigma) + |h_\Sigma|^2) \, (v^{(j)})^2 \\ 
&+ \int_{\Sigma^{(j)}} (D^2 \rho)(\nu_\Sigma,\nu_\Sigma) \, (v^{(j)})^2 - \int_{\Sigma^{(j)}} \rho^{-1} \, \langle \nabla \rho,\nu_\Sigma \rangle^2 \, (v^{(j)})^2 = \lambda_j.
\end{align*} 
\end{proposition}

\textbf{Proof.} 
In the following, we fix an integer $j$. Given any positive integer $l$, we can find a pair $(v_{j,l},a_{j,l}) \in H^1(\Sigma^{(j)}) \times \mathbb{R}$ such that $v_{j,l} - aa_{j,l} \, \bar{v} \in H_0^1(\Sigma^{(j)})$, $\int_{\Sigma \cap \{2r_* < r < 5r_*\}} \mu(r) \, v_{j,l}^2 = 1$, and  
\begin{align} 
\label{minimizing.sequence}
&-\Lambda_j \, a_{j,l}^2 + \int_{\Sigma^{(j)}} \rho \, |\nabla^\Sigma v_{j,l}|^2 - \int_{\Sigma^{(j)}} \rho \, (\text{\rm Ric}(\nu_\Sigma,\nu_\Sigma) + |h_\Sigma|^2) \, v_{j,l}^2 \notag \\ 
&+ \int_{\Sigma^{(j)}} (D^2 \rho)(\nu_\Sigma,\nu_\Sigma) \, v_{j,l}^2 - \int_{\Sigma^{(j)}} \rho^{-1} \, \langle \nabla \rho,\nu_\Sigma \rangle^2 \, v_{j,l}^2 \leq \lambda_j + \frac{1}{l}.
\end{align} 
We distinguish two cases: 

\textit{Case 1:} Suppose that $\sup_l \|(v_{j,l},a_{j,l})\|_{L^2(\Sigma^{(j)}) \times \mathbb{R}} < \infty$. Using (\ref{minimizing.sequence}), we deduce that $\sup_l \|(v_{j,l},a_{j,l})\|_{H^1(\Sigma^{(j)}) \times \mathbb{R}} < \infty$. After passing to a subsequence, the sequence $(v_{j,l},a_{j,l})$ converges strongly in $L^2(\Sigma^{(j)}) \times \mathbb{R}$ and weakly in $H^1(\Sigma^{(j)}) \times \mathbb{R}$. Let $(v^{(j)},a^{(j)}) \in H^1(\Sigma^{(j)}) \times \mathbb{R}$ denote the limit. Then $v^{(j)} - a^{(j)} \, \bar{v} \in H_0^1(\Sigma^{(j)})$ and $\int_{\Sigma \cap \{2r_* < r < 5r_*\}} \mu(r) \, (v^{(j)})^2 = 1$. Using (\ref{minimizing.sequence}) and the lower semicontinuity of the Dirichlet energy, we obtain 
\begin{align*} 
&-\Lambda_j \, (a^{(j)})^2 + \int_{\Sigma^{(j)}} \rho \, |\nabla^\Sigma v^{(j)}|^2 - \int_{\Sigma^{(j)}} \rho \, (\text{\rm Ric}(\nu_\Sigma,\nu_\Sigma) + |h_\Sigma|^2) \, (v^{(j)})^2 \\ 
&+ \int_{\Sigma^{(j)}} (D^2 \rho)(\nu_\Sigma,\nu_\Sigma) \, (v^{(j)})^2 - \int_{\Sigma^{(j)}} \rho^{-1} \, \langle \nabla \rho,\nu_\Sigma \rangle^2 \, (v^{(j)})^2 \leq \lambda_j.
\end{align*} 
By definition of $\lambda_j$, equality holds in the preceding inequality. Thus, $(v^{(j)},a^{(j)})$ has all the required properties.

\textit{Case 2:} Suppose that $\sup_l \|(v_{j,l},a_{j,l})\|_{L^2(\Sigma^{(j)}) \times \mathbb{R}} = \infty$. After passing to a subsequence, we may assume that $\|(v_{j,l},a_{j,l})\|_{L^2(\Sigma^{(j)}) \times \mathbb{R}} \to \infty$ as $l \to \infty$. For each $l$, we define 
\[(\tilde{v}_{j,l},\tilde{a}_{j,l}) = \frac{(v_{j,l},a_{j,l})}{\|(v_{j,l},a_{j,l})\|_{L^2(\Sigma^{(j)}) \times \mathbb{R}}} \in H^1(\Sigma^{(j)}) \times \mathbb{R}.\]
Clearly, $\|(\tilde{v}_{j,l},\tilde{a}_{j,l})\|_{L^2(\Sigma^{(j)}) \times \mathbb{R}} = 1$ for each $l$. Using (\ref{minimizing.sequence}), we deduce that $\sup_l \|(\tilde{v}_{j,l},\tilde{a}_{j,l})\|_{H^1(\Sigma^{(j)}) \times \mathbb{R}} < \infty$. After passing to a subsequence, the sequence $(\tilde{v}_{j,l},\tilde{a}_{j,l})$ converges strongly in $L^2(\Sigma^{(j)}) \times \mathbb{R}$ and weakly in $H^1(\Sigma^{(j)}) \times \mathbb{R}$. The limit $(\tilde{v}^{(j)},\tilde{a}^{(j)}) \in H^1(\Sigma^{(j)}) \times \mathbb{R}$ satisfies $\|(\tilde{v}^{(j)},\tilde{a}^{(j)})\|_{L^2(\Sigma^{(j)}) \times \mathbb{R}} = 1$. Moreover, $\tilde{v}^{(j)} - \tilde{a}^{(j)} \, \bar{v} \in H_0^1(\Sigma^{(j)})$ and $\int_{\Sigma \cap \{2r_* < r < 5r_*\}} \mu(r) \, (\tilde{v}^{(j)})^2 = 0$. Using (\ref{minimizing.sequence}) and the lower semicontinuity of the Dirichlet energy, we obtain 
\begin{align*} 
&-\Lambda_j \, (\tilde{a}^{(j)})^2 + \int_{\Sigma^{(j)}} \rho \, |\nabla^\Sigma \tilde{v}^{(j)}|^2 - \int_{\Sigma^{(j)}} \rho \, (\text{\rm Ric}(\nu_\Sigma,\nu_\Sigma) + |h_\Sigma|^2) \, (\tilde{v}^{(j)})^2 \\ 
&+ \int_{\Sigma^{(j)}} (D^2 \rho)(\nu_\Sigma,\nu_\Sigma) \, (\tilde{v}^{(j)})^2 - \int_{\Sigma^{(j)}} \rho^{-1} \, \langle \nabla \rho,\nu_\Sigma \rangle^2 \, (\tilde{v}^{(j)})^2 \leq 0.
\end{align*} 
By Proposition \ref{consequence.of.stability.inequality.2}, equality holds in the preceding inequality. Consequently, $\tilde{v}^{(j)}$ is a weak solution of the PDE $\mathbb{L}_\Sigma \tilde{v}^{(j)} = 0$ on $\Sigma^{(j)}$. Moreover, $\tilde{v}^{(j)} - \tilde{a}^{(j)} \, \bar{v} \in H_0^1(\Sigma^{(j)})$. Since $\bar{v}$ is smooth, it follows from standard elliptic regularity theory that $\tilde{v}^{(j)}$ is a smooth solution of the PDE $\mathbb{L}_\Sigma \tilde{v}^{(j)} = 0$ on $\Sigma^{(j)}$ with Dirichlet boundary condition $\tilde{v}^{(j)} = \tilde{a}^{(j)} \, \bar{v}$ on $\partial \Sigma^{(j)}$. Since $\int_{\Sigma \cap \{2r_* < r < 5r_*\}} \mu(r) \, (\tilde{v}^{(j)})^2 = 0$, the function $\tilde{v}^{(j)}$ vanishes on the set $\Sigma \cap \{3r_* \leq r \leq 4r_*\} \subset \Sigma^{(j)}$. Since $\Sigma^{(j)}$ is connected, standard unique continuation theorems for elliptic PDE (see e.g. \cite{Aronszajn}) imply that $\tilde{v}^{(j)}$ vanishes identically. Using the boundary condition, we conclude that $\tilde{a}^{(j)}=0$. This contradicts the fact that $\|(\tilde{v}^{(j)},\tilde{a}^{(j)})\|_{L^2(\Sigma^{(j)}) \times \mathbb{R}} = 1$. This completes the proof of Proposition \ref{existence.of.eigenfunction}. \\

Let $(v^{(j)},a^{(j)}) \in H^1(\Sigma^{(j)}) \times \mathbb{R}$ denote the minimizer constructed in Proposition \ref{existence.of.eigenfunction}. After replacing the pair $(v^{(j)},a^{(j)}) \in H^1(\Sigma^{(j)}) \times \mathbb{R}$ by the pair $(|v^{(j)}|,|a^{(j)}|) \in H^1(\Sigma^{(j)}) \times \mathbb{R}$, we may assume that $v^{(j)}$ is a nonnegative function on $\Sigma^{(j)}$ and $a^{(j)}$ is a nonnegative real number. 

\begin{proposition}
\label{properties.of.minimizer}
For each $j$, the function $v^{(j)}$ is a smooth solution of the PDE 
\begin{equation} 
\label{Jacobi.operator.of.v_j}
\mathbb{L}_\Sigma v^{(j)} = \begin{cases} \lambda_j \, \mu(r) \, v^{(j)} & \text{\rm on $\Sigma \cap \{2r_* < r < 5r_*\}$} \\ 0 & \text{\rm on $\Sigma^{(j)} \setminus \{2r_* < r < 5r_*\}$} \end{cases}
\end{equation} 
with Dirichlet boundary condition $v^{(j)} = a^{(j)} \, \bar{v}$ on $\partial \Sigma^{(j)}$. Moreover, 
\begin{equation} 
\label{normal.derivative.of.v_j}
\int_{\partial \Sigma^{(j)}} \rho \, \bar{v} \, \langle \nabla^\Sigma v^{(j)},\eta \rangle = \Lambda_j \, a^{(j)}. 
\end{equation}
\end{proposition}

\textbf{Proof.} 
The minimization property of $(v^{(j)},a^{(j)})$ implies that 
\begin{align*} 
&-\Lambda_j \, a^{(j)} \, a + \int_{\Sigma^{(j)}} \rho \, \langle \nabla^\Sigma v^{(j)},\nabla^\Sigma v \rangle - \int_{\Sigma^{(j)}} \rho \, (\text{\rm Ric}(\nu_\Sigma,\nu_\Sigma) + |h_\Sigma|^2) \, v^{(j)} \, v \\ 
&+ \int_{\Sigma^{(j)}} (D^2 \rho)(\nu_\Sigma,\nu_\Sigma) \, v^{(j)} \, v - \int_{\Sigma^{(j)}} \rho^{-1} \, \langle \nabla \rho,\nu_\Sigma \rangle^2 \, v^{(j)} \, v \\ 
&= \lambda_j \int_{\Sigma \cap \{2r_* < r < 5r_*\}} \mu(r) \, v^{(j)} \, v 
\end{align*} 
for all pairs $(v,a) \in H^1(\Sigma^{(j)}) \times \mathbb{R}$ satisfying $v - a \, \bar{v} \in H_0^1(\Sigma^{(j)})$. Putting $a=0$, we conclude that $v^{(j)}$ is a weak solution of (\ref{Jacobi.operator.of.v_j}). Moreover, $v^{(j)} - a^{(j)} \, \bar{v} \in H_0^1(\Sigma^{(j)})$. Since $\bar{v}$ is smooth, it follows from standard elliptic regularity theory that $v^{(j)}$ is a smooth solution of (\ref{Jacobi.operator.of.v_j}) with Dirichlet boundary condition $v^{(j)} = a^{(j)} \, \bar{v}$ on $\partial \Sigma^{(j)}$. Using the identity above and integration by parts, we obtain 
\[-\Lambda_j \, a^{(j)} \, a + \int_{\partial \Sigma^{(j)}} \rho \, \langle \nabla^\Sigma v^{(j)},\eta \rangle \, v = 0\] 
for all pairs $(v,a) \in H^1(\Sigma^{(j)}) \times \mathbb{R}$ satisfying $v - a \, \bar{v} \in H_0^1(\Sigma^{(j)})$. From this, the identity (\ref{normal.derivative.of.v_j}) follows. This completes the proof of Proposition \ref{properties.of.minimizer}. \\

\begin{lemma}
\label{strict.positivity.of.v_j}
The function $v^{(j)}$ is strictly positive at each point in the interior of $\Sigma^{(j)}$.
\end{lemma} 

\textbf{Proof.} 
Note that the function $v^{(j)}$ is nonnegative. Therefore, the assertion follows from the strict maximum principle for elliptic PDE. \\

\begin{lemma} 
\label{strict.positivity.of.a_j}
The number $a^{(j)}$ is strictly positive. Consequently, the function $v^{(j)}$ is strictly positive at each point on the boundary $\partial \Sigma^{(j)}$.
\end{lemma}

\textbf{Proof.} 
We argue by contradiction. Suppose that $a^{(j)} = 0$. Then the function $v^{(j)}$ vanishes along the boundary $\partial \Sigma^{(j)}$. Using Lemma \ref{strict.positivity.of.v_j} and the Hopf boundary point lemma (see \cite{Gilbarg-Trudinger}, Lemma 3.4), we conclude that $\langle \nabla^\Sigma v^{(j)},\eta \rangle < 0$ at each point on the boundary $\partial \Sigma^{(j)}$. This contradicts (\ref{normal.derivative.of.v_j}). This completes the proof of Lemma \ref{strict.positivity.of.a_j}. \\

For each $j$, we define a smooth function $w^{(j)}: \Sigma^{(j)} \to \mathbb{R}$ by 
\begin{equation} 
\label{definition.of.w_j}
w^{(j)} = \frac{v^{(j)}}{a^{(j)}} - \bar{v}. 
\end{equation} 
Note that $w^{(j)}$ is well-defined by Lemma \ref{strict.positivity.of.a_j}. Moreover, $w^{(j)} = 0$ on $\partial \Sigma^{(j)}$. 

\begin{lemma} 
\label{normal.derivative.of.w_j}
We have 
\[\bigg | \int_{\partial \Sigma^{(j)}} \rho \, \bar{v} \, \langle \nabla^\Sigma w^{(j)},\eta \rangle \bigg | \leq C.\] 
\end{lemma}

\textbf{Proof.} 
Using (\ref{normal.derivative.of.v_j}), we obtain 
\[\int_{\partial \Sigma^{(j)}} \rho \, \bar{v} \, \langle \nabla^\Sigma w^{(j)},\eta \rangle = \Lambda_j - \int_{\partial \Sigma^{(j)}} \rho \, \bar{v} \, \langle \nabla^\Sigma \bar{v},\eta \rangle.\] 
Therefore, the assertion follows from Lemma \ref{asymptotics.for.Lambda_j}. \\

\begin{lemma}
\label{normalization} 
We have $\sup_{\Sigma \cap \{2r_* < r < 5r_*\}} v^{(j)} \geq \frac{1}{C}$ and $\inf_{\Sigma \cap \{2r_* < r < 5r_*\}} v^{(j)} \leq C$ for some uniform constant $C$. 
\end{lemma}

\textbf{Proof.} 
This follows from the fact that $\int_{\Sigma \cap \{2r_* < r < 5r_*\}} \mu(r) \, (v^{(j)})^2 = 1$ (see Proposition \ref{existence.of.eigenfunction}). \\

\begin{lemma}
\label{lower.bound.for.a_j}
The sequence $a^{(j)}$ is bounded from below by a positive constant.
\end{lemma}

\textbf{Proof.} 
Suppose that the assertion is false. After passing to a subsequence, we may assume that $a^{(j)} \to 0$. By Lemma \ref{Jacobi.operator.of.bar.v}, we can find a large constant $C_0$ such that 
\[|\mathbb{L}_\Sigma \bar{v}| \leq C_0 \, r^{1-n-\delta}\]
on $\Sigma \cap \{r \geq 6r_*\}$. Using (\ref{Jacobi.operator.of.v_j}) and (\ref{definition.of.w_j}), we obtain 
\[\mathbb{L}_\Sigma w^{(j)} = -\mathbb{L}_\Sigma \bar{v} \geq -C_0 \, r^{1-n-\delta}\] 
on $\Sigma \cap \{6r_* \leq r \leq r_j\}$. On the other hand, it follows from Lemma \ref{Jacobi.operator.of.test.function} that we can find a large constant $\sigma \in [6r_*,\infty)$ with the following properties: 
\begin{itemize}
\item The function $-\rho \, (\text{\rm Ric}(\nu_\Sigma,\nu_\Sigma) + |h_\Sigma|^2) + (D^2 \rho)(\nu_\Sigma,\nu_\Sigma) - \rho^{-1} \, \langle \nabla \rho,\nu_\Sigma \rangle^2$ is positive on $\Sigma \cap \{r \geq \sigma\}$. 
\item $\mathbb{L}_\Sigma((r^{-N} + r^{-N-\delta}) \, \bar{v}) \leq -b_{n-2} \, \delta N \, r^{1-n-\delta}$ on $\Sigma \cap \{r \geq \sigma\}$. 
\end{itemize} 
It follows from Lemma \ref{normalization} and the Harnack inequality that $\inf_{\Sigma \cap \{r=\sigma\}} v^{(j)}$ is uniformly bounded from below by a positive constant that may depend on $\sigma$, but not on $j$. Since $a^{(j)} \to 0$, it follows that $\inf_{\Sigma \cap \{r=\sigma\}} a^{(j)} \, w^{(j)}$ is uniformly bounded from below by a positive constant that may depend on $\sigma$, but not on $j$. Therefore, if $j$ is sufficiently large, then 
\[(r^{-N} + r^{-N-\delta} - r_j^{-N} - r_j^{-N-\delta}) \, \bar{v} \leq K \, a^{(j)} \, w^{(j)}\] 
on $\Sigma \cap \{r=\sigma\}$, where $K$ is a constant that may depend on $\sigma$, but not on $j$. Moreover, $b_{n-2} \, \delta N - C_0 \, (r_j^N + r_j^{-N-\delta}) > C_0 \, K \, a^{(j)}$ if $j$ is sufficiently large. This implies 
\[\mathbb{L}_\Sigma((r^{-N} + r^{-N-\delta} - r_j^{-N} - r_j^{-N-\delta}) \, \bar{v}) < K \, a^{(j)} \, \mathbb{L}_\Sigma w^{(j)}\] 
on $\Sigma \cap \{\sigma \leq r \leq r_j\}$. We now apply a standard comparison principle (cf. Theorem 3.3 in \cite{Gilbarg-Trudinger}) to the operator $\mathbb{L}_\Sigma$ on $\Sigma \cap \{\sigma \leq r \leq r_j\}$. If $j$ is sufficiently large, we conclude that 
\[(r^{-N} + r^{-N-\delta} - r_j^{-N} - r_j^{-N-\delta}) \, \bar{v} \leq K \, a^{(j)} \, w^{(j)}\] 
on $\Sigma \cap \{\sigma \leq r \leq r_j\}$, and equality holds on the set $\Sigma \cap \{r=r_j\}$. In particular, if $j$ is sufficiently large, then 
\[-\big \langle \nabla^\Sigma ((r^{-N} + r^{-N-\delta} - r_j^{-N} - r_j^{-N-\delta}) \, \bar{v}),\eta \big \rangle \leq -K \, a^{(j)} \, \langle \nabla^\Sigma w^{(j)},\eta \rangle\] 
at each point on $\partial \Sigma^{(j)}$. This implies 
\[-\bar{v} \, \langle \nabla^\Sigma (r^{-N} + r^{-N-\delta}),\eta \rangle \leq -K \, a^{(j)} \, \langle \nabla^\Sigma w^{(j)},\eta \rangle\] 
at each point on $\partial \Sigma^{(j)}$. Thus, 
\[-\int_{\partial \Sigma^{(j)}} \rho \, \bar{v}^2 \, \langle \nabla^\Sigma (r^{-N} + r^{-N-\delta}),\eta \rangle \leq -K \, a^{(j)} \, \int_{\partial \Sigma^{(j)}} \rho \, \bar{v} \, \langle \nabla^\Sigma w^{(j)},\eta \rangle\] 
if $j$ is sufficiently large. Finally, we send $j \to \infty$. The expression on the left hand side is bounded from below by a positive constant, while the expression on the right hand side converges to $0$ by Lemma \ref{normal.derivative.of.w_j}. This is a contradiction. This completes the proof of Lemma \ref{lower.bound.for.a_j}. \\

\begin{lemma}
\label{lower.bound.for.w_j}
We can find a large constant $\sigma$ with the following property. If $j$ is sufficiently large, then $w^{(j)} \geq -2\sigma^N \, r^{-N} \, \bar{v}$ on $\Sigma \cap \{\sigma \leq r \leq r_j\}$. 
\end{lemma}

\textbf{Proof.} 
Using (\ref{Jacobi.operator.of.v_j}), (\ref{definition.of.w_j}), and Lemma \ref{Jacobi.operator.of.bar.v}, we obtain 
\[\mathbb{L}_\Sigma w^{(j)} = -\mathbb{L}_\Sigma \bar{v} \geq -C_0 \, r^{1-n-\delta}\] 
on $\Sigma \cap \{6r_* \leq r \leq r_j\}$, where $C_0$ is independent of $j$. On the other hand, it follows from Lemma \ref{Jacobi.operator.of.test.function} that we can find a large constant $\sigma \in [6r_*,\infty)$ with the following properties: 
\begin{itemize}
\item The function $-\rho \, (\text{\rm Ric}(\nu_\Sigma,\nu_\Sigma) + |h_\Sigma|^2) + (D^2 \rho)(\nu_\Sigma,\nu_\Sigma) - \rho^{-1} \, \langle \nabla \rho,\nu_\Sigma \rangle^2$ is positive on $\Sigma \cap \{r \geq \sigma\}$. 
\item $\mathbb{L}_\Sigma((r^{-N} - r^{-N-\delta}) \, \bar{v}) \geq b_{n-2} \, \delta N \, r^{1-n-\delta}$ on $\Sigma \cap \{r \geq \sigma\}$. 
\end{itemize} 
By increasing $\sigma$ if necessary, we may arrange that $C_0 < 2b_{n-2} \, \delta N \, \sigma^N$ and $2\sigma^{-\delta} < 1$. 

We next observe that 
\[w^{(j)} \geq -\bar{v} \geq -2\sigma^N \, (r^{-N} - r^{-N-\delta}) \, \bar{v}\] 
on $\Sigma \cap \{r=\sigma\}$. Moreover, using the inequality $C_0 < 2b_{n-2} \, \delta N \, \sigma^N$, we obtain 
\[\mathbb{L}_\Sigma w^{(j)} > -2\sigma^N \, \mathbb{L}_\Sigma((r^{-N} - r^{-N-\delta}) \, \bar{v})\] 
on $\Sigma \cap \{\sigma \leq r \leq r_j\}$. We now apply a standard comparison principle (cf. Theorem 3.3 in \cite{Gilbarg-Trudinger}) to the operator $\mathbb{L}_\Sigma$ on $\Sigma \cap \{\sigma \leq r \leq r_j\}$. If $j$ is sufficiently large, we conclude that 
\[w^{(j)} \geq -2\sigma^N \, (r^{-N} - r^{-N-\delta}) \, \bar{v}\] 
on $\Sigma \cap \{\sigma \leq r \leq r_j\}$. This completes the proof of Lemma \ref{lower.bound.for.w_j}. \\

\begin{lemma}
\label{upper.bound.for.a_j}
The sequence $a^{(j)}$ is bounded from above.
\end{lemma}

\textbf{Proof.} 
Let $\sigma$ denote the constant in Lemma \ref{lower.bound.for.w_j}. It follows from Lemma \ref{lower.bound.for.w_j} that $w^{(j)} \geq -2^{1-N} \, \bar{v}$ on $\Sigma \cap \{r=2\sigma\}$. This implies $v^{(j)} \geq (1-2^{1-N}) \, a^{(j)} \, \bar{v}$ on $\Sigma \cap \{r=2\sigma\}$. On the other hand, it follows from the Harnack inequality and Lemma \ref{normalization} that $\sup_{\Sigma \cap \{r=2\sigma\}} v^{(j)}$ is bounded from above by a constant that may depend on $\sigma$, but not on $j$. Putting these facts together, the assertion follows. \\

\begin{lemma} 
\label{upper.bound.for.w_j}
We can find a large constant $\sigma$ and a large constant $C$ with the following property. If $j$ is sufficiently large, then $w^{(j)} \leq C \, r^{-N} \, \bar{v}$ on $\Sigma \cap \{\sigma \leq r \leq r_j\}$. 
\end{lemma} 

\textbf{Proof.} 
Using (\ref{Jacobi.operator.of.v_j}), (\ref{definition.of.w_j}), and Lemma \ref{Jacobi.operator.of.bar.v}, we obtain 
\[\mathbb{L}_\Sigma w^{(j)} = -\mathbb{L}_\Sigma \bar{v} \leq C_0 \, r^{1-n-\delta}\] 
on $\Sigma \cap \{6r_* \leq r \leq r_j\}$, where $C_0$ is independent of $j$. On the other hand, it follows from Lemma \ref{Jacobi.operator.of.test.function} that we can find a large constant $\sigma \in [6r_*,\infty)$ with the following properties: 
\begin{itemize}
\item The function $-\rho \, (\text{\rm Ric}(\nu_\Sigma,\nu_\Sigma) + |h_\Sigma|^2) + (D^2 \rho)(\nu_\Sigma,\nu_\Sigma) - \rho^{-1} \, \langle \nabla \rho,\nu_\Sigma \rangle^2$ is positive on $\Sigma \cap \{r \geq \sigma\}$. 
\item $\mathbb{L}_\Sigma((r^{-N} - r^{-N-\delta}) \, \bar{v}) \geq b_{n-2} \, \delta N \, r^{1-n-\delta}$ on $\Sigma \cap \{r \geq \sigma\}$. 
\end{itemize}
It follows from the Harnack inequality and Lemma \ref{normalization} that $\sup_{\Sigma \cap \{r=\sigma\}} v^{(j)}$ is bounded from above by a positive constant that may depend on $\sigma$, but not on $j$. Moreover, Lemma \ref{lower.bound.for.a_j} implies that $a^{(j)}$ is bounded from below by a positive constant which is independent of $j$. Therefore, if $j$ is sufficiently large, then 
\[w^{(j)} \leq \frac{v^{(j)}}{a^{(j)}} \leq K \, (r^{-N} - r^{-N-\delta}) \, \bar{v}\] 
on $\Sigma \cap \{r=\sigma\}$, where $K$ is a large constant that may depend on $\sigma$, but not on $j$. By increasing $K$ if necessary, we may arrange that $C_0 < b_{n-2} \, \delta N \, K$. This implies 
\[\mathbb{L}_\Sigma w^{(j)} < K \, \mathbb{L}_\Sigma((r^{-N} - r^{-N-\delta}) \, \bar{v})\] 
on $\Sigma \cap \{\sigma \leq r \leq r_j\}$. We now apply a standard comparison principle (cf. Theorem 3.3 in \cite{Gilbarg-Trudinger}) to the operator $\mathbb{L}_\Sigma$ on $\Sigma \cap \{\sigma \leq r \leq r_j\}$. If $j$ is sufficiently large, we conclude that 
\[w^{(j)} \leq K \, (r^{-N} - r^{-N-\delta}) \, \bar{v}\] 
on $\Sigma \cap \{\sigma \leq r \leq r_j\}$. This completes the proof of Lemma \ref{upper.bound.for.w_j}. \\

\begin{proposition}
\label{existence.of.check.v}
After passing to a subsequence if necessary, the functions $\frac{v^{(j)}}{a^{(j)}}$ converge in $C_{\text{\rm loc}}^\infty$ to a positive smooth function $\check{v}: \Sigma \to \mathbb{R}$. The function $\check{v}$ satisfies the PDE 
\begin{equation} 
\label{Jacobi.operator.of.check.v}
\mathbb{L}_\Sigma \check{v} = \begin{cases} \lambda_\infty \, \mu(r) \, \check{v} & \text{\rm on $\Sigma \cap \{2r_* < r < 5r_*\}$} \\ 0 & \text{\rm on $\Sigma \setminus \{2r_* < r < 5r_*\}$,} \end{cases} 
\end{equation}
where $\lambda_\infty = \lim_{j \to \infty} \lambda_j$. Moreover, we can find a large constant $\sigma$ and a large constant $C$ such that $|\check{v}-\bar{v}| \leq C \, r^{-N} \, \bar{v}$ on $\Sigma \cap \{r \geq \sigma\}$. 
\end{proposition} 

\textbf{Proof.} 
Using Lemma \ref{normalization}, Lemma \ref{lower.bound.for.a_j}, and Lemma \ref{upper.bound.for.a_j}, we obtain $\sup_{\Sigma \cap \{2r_* < r < 5r_*\}} \frac{v^{(j)}}{a^{(j)}} \geq \frac{1}{C}$ and $\inf_{\Sigma \cap \{2r_* < r < 5r_*\}} \frac{v^{(j)}}{a^{(j)}} \leq C$ for some uniform constant $C$. The Harnack inequality gives uniform upper and lower bounds for $\frac{v^{(j)}}{a^{(j)}}$ on every compact subset of $\Sigma$. After passing to a subsequence if necessary, the functions $\frac{v^{(j)}}{a^{(j)}}$ converge in $C_{\text{\rm loc}}^\infty$ to a positive smooth function $\check{v}: \Sigma \to \mathbb{R}$. Finally, Lemma \ref{lower.bound.for.w_j} and Lemma \ref{upper.bound.for.w_j} imply that $|\check{v}-\bar{v}| \leq C \, r^{-N} \, \bar{v}$ on $\Sigma \cap \{r \geq \sigma\}$. \\

\begin{proposition}
\label{higher.derivative.bound.for.check.v.minus.bar.v}
Let $\check{v}: \Sigma \to \mathbb{R}$ be defined as in Proposition \ref{existence.of.check.v}. Then $|D_{\text{\rm hyp}}^{\Sigma,m} (\check{v}-\bar{v})|_{g_{\text{\rm hyp}}} \leq O(r^{1-N})$ for every nonnegative integer $m$.
\end{proposition}

\textbf{Proof.} 
Proposition \ref{existence.of.check.v} implies that $|\check{v}-\bar{v}| \leq O(r^{1-N})$. Moreover, using Lemma \ref{Jacobi.operator.of.bar.v} and the PDE (\ref{Jacobi.operator.of.check.v}), we obtain $|D_{\text{\rm hyp}}^{\Sigma,m} \mathbb{L}_\Sigma (\check{v}-\bar{v})|_{g_{\text{\rm hyp}}} \leq O(r^{1-n})$ for every nonnegative integer $m$. The assertion now follows from standard interior estimates for elliptic PDE. \\

\begin{proposition}
\label{asymptotic.expansion.for.check.v.minus.bar.v}
Let $\check{v}: \Sigma \to \mathbb{R}$ be defined as in Proposition \ref{existence.of.check.v}. Then we can find a function $A \in C^{\frac{\delta}{10}}(T^{n-2},\gamma)$ such that 
\[|\check{v} - \bar{v} - b_{n-2} \, r^{1-N} \, A(\theta_0,\hdots,\theta_{n-3})| \leq O(r^{1-N-\frac{\delta}{10}})\] 
and 
\[|\langle \nabla^\Sigma r,\nabla^\Sigma (\check{v} - \bar{v}) \rangle + (N-1) \, b_{n-2} \, r^{2-N} \, A(\theta_0,\hdots,\theta_{n-3})| \leq O(r^{2-N-\frac{\delta}{10}}).\] 
\end{proposition} 

\textbf{Proof.} 
It follows from (\ref{Jacobi.operator.of.v_j}), (\ref{definition.of.w_j}), and Lemma \ref{Jacobi.operator.of.bar.v} that 
\begin{equation} 
\label{Lw}
|\mathbb{L}_\Sigma w^{(j)}| = |\mathbb{L}_\Sigma \bar{v}| \leq C \, r^{1-n-\delta} 
\end{equation}
and 
\begin{equation} 
\label{derivatives.of.Lw}
\sum_{m=1}^2 |D_{\text{\rm hyp}}^{\Sigma,m} \mathbb{L}_\Sigma w^{(j)}|_{g_{\text{\rm hyp}}} = \sum_{m=1}^2 |D_{\text{\rm hyp}}^{\Sigma,m} \mathbb{L}_\Sigma \bar{v}|_{g_{\text{\rm hyp}}} \leq C \, r^{1-n} 
\end{equation}
on $\Sigma \cap \{6r_* \leq r \leq r_j\}$. Moreover, it follows from Lemma \ref{lower.bound.for.w_j} and Lemma \ref{upper.bound.for.w_j} that $|w^{(j)}| \leq C \, r^{1-N}$ on $\Sigma \cap \{6r_* \leq r \leq r_j\}$. Finally, we know that $w^{(j)} = 0$ on $\Sigma \cap \{r=r_j\}$. Using (\ref{Lw}) and (\ref{derivatives.of.Lw}) together with the standard regularity theory for elliptic PDE (see \cite{Gilbarg-Trudinger}, Theorem 6.6), we conclude that 
\begin{equation} 
\label{derivatives.of.w}
\sum_{m=0}^3 |D_{\text{\rm hyp}}^{\Sigma,m} w^{(j)}|_{g_{\text{\rm hyp}}} \leq C \, r^{1-N} 
\end{equation}
on $\Sigma \cap \{8r_* \leq r \leq r_j\}$. We define a function $\zeta^{(j)}$ on $\Sigma \cap \{8r_* \leq r \leq r_j\}$ by 
\[-\text{\rm div}_{g_{\text{\rm hyp}}}(r^{N-n} \, dw^{(j)}) + (N-1) \, r^{N-n} \, w^{(j)} = \zeta^{(j)}.\] 
Using (\ref{Lw}) and (\ref{derivatives.of.w}), we obtain $|\zeta^{(j)}| \leq C \, r^{1-n-\delta}$ on $\Sigma \cap \{8r_* \leq r \leq r_j\}$. Moreover, (\ref{derivatives.of.w}) implies $|d\zeta^{(j)}|_{g_{\text{\rm hyp}}} \leq C \, r^{1-n}$ on $\Sigma \cap \{8r_* \leq r \leq r_j\}$. If we apply Theorem \ref{asymptotics.for.solutions.of.linear.PDE} to the functions $w^{(j)} = \frac{v^{(j)}}{a^{(j)}} - \bar{v}$, the assertion follows. This completes the proof of Proposition \ref{asymptotic.expansion.for.check.v.minus.bar.v}. \\

\begin{proposition} 
\label{flux.1}
We have 
\[\int_{T^{n-2} \times \{t_*\}} \Big ( N \, b_{n-2}^2 \, A - \frac{\partial^2 u}{\partial \theta_{n-2}^2} \Big ) \, d\text{\rm vol}_\gamma = 0.\] 
\end{proposition} 

\textbf{Proof.} 
It follows from Lemma \ref{upper.bound.for.w_j} that $v^{(j)} \leq C \, a^{(j)} \, \bar{v}$ in the domain $\Sigma \cap \{8r_* \leq r \leq r_j\}$, where $C$ is independent of $j$. Using (\ref{Jacobi.operator.of.v_j}) and Lemma \ref{Jacobi.operator.of.bar.v}, we obtain 
\begin{align} 
\label{divergence.identity.1}
|\text{\rm div}_\Sigma(\rho \, \bar{v} \, \nabla^\Sigma v^{(j)} - \rho \, v^{(j)} \, \nabla^\Sigma \bar{v})| 
&= |v^{(j)} \, \mathbb{L}_\Sigma \bar{v} - \bar{v} \, \mathbb{L}_\Sigma v^{(j)}| \notag \\ 
&= v^{(j)} \, |\mathbb{L}_\Sigma \bar{v}| \\ 
&\leq C \, a^{(j)} \, r^{2-n-\delta} \notag
\end{align}
on $\Sigma \cap \{8r_* \leq r \leq r_j\}$, where $C$ is independent of $j$. We integrate both sides of (\ref{divergence.identity.1}) over $\Sigma \cap \{\bar{r} \leq r \leq r_j\}$, where $\bar{r} \in (8r_*,r_j)$. Using the divergence theorem, we deduce that 
\begin{align} 
\label{divergence.identity.2} 
&\bigg | \int_{\partial \Sigma^{(j)}} \rho \, \bar{v} \, \langle \nabla^\Sigma v^{(j)},\eta \rangle - \int_{\partial \Sigma^{(j)}} \rho \, v^{(j)} \, \langle \nabla^\Sigma \bar{v},\eta \rangle \notag \\ 
&- \int_{\Sigma \cap \{r=\bar{r}\}} \rho \, \bar{v} \, \Big \langle \nabla^\Sigma v^{(j)},\frac{\nabla^\Sigma r}{|\nabla^\Sigma r|} \Big \rangle + \int_{\Sigma \cap \{r=\bar{r}\}} \rho \, v^{(j)} \, \Big \langle \nabla^\Sigma \bar{v},\frac{\nabla^\Sigma r}{|\nabla^\Sigma r|} \Big \rangle \bigg | \\ 
&\leq C \, a^{(j)} \, \bar{r}^{-\delta} \notag 
\end{align}
for $\bar{r} \in (8r_*,r_j)$, where $C$ is independent of $\bar{r}$ and $j$. In the next step, we use the identity (\ref{normal.derivative.of.v_j}) and the fact that $v^{(j)} = a^{(j)} \, \bar{v}$ on $\partial \Sigma^{(j)}$. This gives 
\begin{align} 
\label{divergence.identity.3}
&\bigg | \bigg ( \Lambda_j - \int_{\partial \Sigma^{(j)}} \rho \, \bar{v} \, \langle \nabla^\Sigma \bar{v},\eta \rangle \bigg ) \, a^{(j)} \notag \\ 
&- \int_{\Sigma \cap \{r=\bar{r}\}} \rho \, \bar{v} \, \Big \langle \nabla^\Sigma v^{(j)},\frac{\nabla^\Sigma r}{|\nabla^\Sigma r|} \Big \rangle + \int_{\Sigma \cap \{r=\bar{r}\}} \rho \, v^{(j)} \, \Big \langle \nabla^\Sigma \bar{v},\frac{\nabla^\Sigma r}{|\nabla^\Sigma r|} \Big \rangle \bigg | \\ 
&\leq C \, a^{(j)} \, \bar{r}^{-\delta} \notag 
\end{align}
for $\bar{r} \in (8r_*,r_j)$, where $C$ is independent of $\bar{r}$ and $j$. We divide both sides of (\ref{divergence.identity.3}) by $a^{(j)}$ and send $j \to \infty$, while keeping $\bar{r}$ fixed. Using Lemma \ref{asymptotics.for.Lambda_j} and Proposition \ref{existence.of.check.v}, we conclude that 
\begin{align} 
\label{limit.1}
&\bigg | -\int_{T^{n-2} \times \{t_*\}} \frac{\partial^2 u}{\partial \theta_{n-2}^2} \, d\text{\rm vol}_\gamma \notag \\ 
&- \int_{\Sigma \cap \{r=\bar{r}\}} \rho \, \bar{v} \, \Big \langle \nabla^\Sigma \check{v},\frac{\nabla^\Sigma r}{|\nabla^\Sigma r|} \Big \rangle + \int_{\Sigma \cap \{r=\bar{r}\}} \rho \, \check{v} \, \Big \langle \nabla^\Sigma \bar{v},\frac{\nabla^\Sigma r}{|\nabla^\Sigma r|} \Big \rangle \bigg | \\ 
&\leq C \, \bar{r}^{-\delta} \notag
\end{align}
for $\bar{r} > 8r_*$, where $C$ is independent of $\bar{r}$. Finally, we send $\bar{r} \to \infty$. Using Lemma \ref{asymptotics.for.bar.v} and Proposition \ref{asymptotic.expansion.for.check.v.minus.bar.v}, we obtain 
\begin{equation} 
\label{limit.2}
\int_{\Sigma \cap \{r=\bar{r}\}} \rho \, (\check{v}-\bar{v}) \, \Big \langle \nabla^\Sigma \bar{v},\frac{\nabla^\Sigma r}{|\nabla^\Sigma r|} \Big \rangle \to \int_{T^{n-2} \times \{t_*\}} b_{n-2}^2 \, A \, d\text{\rm vol}_\gamma 
\end{equation} 
and 
\begin{equation} 
\label{limit.3}
\int_{\Sigma \cap \{r=\bar{r}\}} \rho \, \bar{v} \, \Big \langle \nabla^\Sigma (\check{v}-\bar{v}),\frac{\nabla^\Sigma r}{|\nabla^\Sigma r|} \Big \rangle \to -\int_{T^{n-2} \times \{t_*\}} (N-1) \, b_{n-2}^2 \, A \, d\text{\rm vol}_\gamma 
\end{equation} 
and 
as $\bar{r} \to \infty$. Subtracting (\ref{limit.3}) from (\ref{limit.2}) gives 
\begin{align} 
\label{limit.4} 
&-\int_{\Sigma \cap \{r=\bar{r}\}} \rho \, \bar{v} \, \Big \langle \nabla^\Sigma \check{v},\frac{\nabla^\Sigma r}{|\nabla^\Sigma r|} \Big \rangle + \int_{\Sigma \cap \{r=\bar{r}\}} \rho \, \check{v} \, \Big \langle \nabla^\Sigma \bar{v},\frac{\nabla^\Sigma r}{|\nabla^\Sigma r|} \Big \rangle \notag \\ 
&\to \int_{T^{n-2} \times \{t_*\}} N \, b_{n-2}^2 \, A \, d\text{\rm vol}_\gamma 
\end{align}
as $\bar{r} \to \infty$. If we combine (\ref{limit.1}) and (\ref{limit.4}), the assertion follows. This completes the proof of Proposition \ref{flux.1}. \\

\begin{corollary} 
\label{flux.2}
We have 
\[\int_{T^{n-2} \times \{t_*\}} \Big ( N \, \text{\rm tr}_\gamma(Q) + 2N \, (P+A) + \Big ( \frac{2}{Nb_0} \Big )^N \Big ) \, d\text{\rm vol}_\gamma \leq 0.\] 
\end{corollary} 

\textbf{Proof.}
Integrating the pointwise inequality (\ref{pde.for.u}) over $T^{n-2} \times \{t_*\}$ gives 
\[\int_{T^{n-2} \times \{t_*\}} \Big ( b_{n-2}^{-2} \, \frac{\partial^2 u}{\partial \theta_{n-2}^2} + \frac{N}{2} \, \text{\rm tr}_\gamma(Q) + N \, P + \frac{1}{2} \, \Big ( \frac{2}{Nb_0} \Big )^N \Big ) \, d\text{\rm vol}_\gamma \leq 0.\] 
On the other hand, 
\[\int_{T^{n-2} \times \{t_*\}} \Big ( NA - b_{n-2}^{-2} \, \frac{\partial^2 u}{\partial \theta_{n-2}^2} \Big ) \, d\text{\rm vol}_\gamma = 0\] 
by Proposition \ref{flux.1}. The assertion follows by adding these two inequalities. This completes the proof of Corollary \ref{flux.2}. \\

\begin{corollary} 
\label{mass.term.Sigma}
Let $\check{\gamma} = \sum_{k=0}^{n-3} b_k^2 \, d\theta_k \otimes d\theta_k$ denote the restriction of $\gamma$ to $T^{n-2} \times \{t_*\}$. Let $\check{Q}$ denote the restriction of $Q$ to $T^{n-2} \times \{t_*\}$, and let $\check{P}$ denote the restriction of the function $P+A + \frac{1}{2} \, b_{n-2}^{-2} \, Q(\frac{\partial}{\partial \theta_{n-2}},\frac{\partial}{\partial \theta_{n-2}})$ to $T^{n-2} \times \{t_*\}$. Then the function $\check{P}: T^{n-2} \times \{t_*\} \to \mathbb{R}$ is H\"older continuous and 
\[\int_{T^{n-2} \times \{t_*\}} \Big ( N \, \text{\rm tr}_{\check{\gamma}}(\check{Q}) + 2N \, \check{P} + \Big ( \frac{2}{Nb_0} \Big )^N \Big ) \, d\text{\rm vol}_\gamma \leq 0.\] 
\end{corollary} 

\textbf{Proof.}
By Proposition \ref{asymptotic.expansion.for.check.v.minus.bar.v}, the function $A$ is H\"older continuous. This implies that the function $\check{P}$ is H\"older continuous. Using Corollary \ref{flux.2} together with the identity $\text{\rm tr}_{\check{\gamma}}(\check{Q}) = \text{\rm tr}_\gamma(Q) - b_{n-2}^{-2} \, Q(\frac{\partial}{\partial \theta_{n-2}},\frac{\partial}{\partial \theta_{n-2}})$, we obtain 
\[\int_{T^{n-2} \times \{t_*\}} \Big ( N \, \text{\rm tr}_{\check{\gamma}}(\check{Q}) + 2N \, \check{P} + \Big ( \frac{2}{Nb_0} \Big )^N \Big ) \, d\text{\rm vol}_\gamma \leq 0.\] 
This completes the proof of Corollary \ref{mass.term.Sigma}. \\

\begin{proposition}
\label{asymptotics.of.check.g}
Let $\check{g}$ denote the induced metric on $\Sigma$, and let $g_{\text{\rm hyp}}$ be defined as in Definition \ref{definition.hyperbolic.metric.on.Sigma}. For every nonnegative integer $m$, we have 
\[|D_{\text{\rm hyp}}^{\Sigma,m} (\check{g}-g_{\text{\rm hyp}})|_{g_{\text{\rm hyp}}} \leq O(r^{-N}).\] 
Moreover, 
\[|\check{g}-g_{\text{\rm hyp}} - r^{2-N} \, \check{Q}|_{g_{\text{\rm hyp}}} \leq O(r^{-N-\delta}),\] 
where $\check{Q}$ denotes the restriction of $Q$ to $T^{n-2} \times \{t_*\}$. 
\end{proposition}

\textbf{Proof.} 
Note that $\check{g}-g_{\text{\rm hyp}}$ is a $(0,2)$-tensor on $\Sigma \cap \{r \geq r_*\}$. It is obtained by restricting the $(0,2)$-tensor $g-\bar{g} + b_{n-2}^2 \, r^2 \, d\theta_{n-2} \otimes d\theta_{n-2}$ on $[r_*,\infty) \times T^{n-1}$ to $\Sigma \cap \{r \geq r_*\}$. The one-form $d\theta_{n-2}$ on $[r_*,\infty) \times T^{n-1}$ restricts to a one-form on $\Sigma \cap \{r \geq r_*\}$. Since $\Sigma$ is tame, this one-form has norm at most $O(r^{-N})$ with respect to the metric $g_{\text{\rm hyp}}$, and its higher order covariant derivatives with respect to $g_{\text{\rm hyp}}$ are bounded by $O(r^{-N})$ as well. From this, the assertion follows. \\

\begin{proposition}
\label{asymptotics.of.check.rho}
Let $\check{v}: \Sigma \to \mathbb{R}$ be defined as in Proposition \ref{existence.of.check.v}. Let us define a positive function $\check{\rho}$ on $\Sigma$ by $\check{\rho} = b_{n-2}^{-1} \, \check{v} \, \rho$. For every nonnegative integer $m$, we have 
\[|D_{\text{\rm hyp}}^{\Sigma,m} (\check{\rho}-r^{N-n+1})|_{g_{\text{\rm hyp}}} \leq O(r^{1-n}).\] 
Moreover, 
\[|\check{\rho}-r^{N-n+1}-r^{1-n} \, \check{P}(\theta_0,\hdots,\theta_{n-3})| \leq O(r^{1-n-\frac{\delta}{10}}),\] 
where $\check{P}$ denotes the restriction of the function $P+A + \frac{1}{2} \, b_{n-2}^{-2} \, Q(\frac{\partial}{\partial \theta_{n-2}},\frac{\partial}{\partial \theta_{n-2}})$ to $T^{n-2} \times \{t_*\}$. 
\end{proposition}

\textbf{Proof.} 
This follows by combining Lemma \ref{asymptotics.for.bar.v}, Proposition \ref{higher.derivative.bound.for.check.v.minus.bar.v}, and Proposition \ref{asymptotic.expansion.for.check.v.minus.bar.v}. \\

Combining Proposition \ref{asymptotics.of.check.g}, Proposition \ref{asymptotics.of.check.rho}, and Corollary \ref{mass.term.Sigma}, we conclude that $(\Sigma,\check{g},\check{\rho})$ is an $(N,n-1)$-dataset.

\begin{proposition}
\label{scalar.curvature.of.hypersurface}
Let $\check{v}: \Sigma \to \mathbb{R}$ be defined as in Proposition \ref{existence.of.check.v}. Let us define a positive function $\check{\rho}$ on $\Sigma$ by $\check{\rho} = b_{n-2}^{-1} \, \check{v} \, \rho$. If $n=N$, we assume that $\rho = 1$ and $R + N(N-1) \geq 0$ at each point on $\Sigma$. If $n<N$, we assume that 
\[-2 \, \Delta \log \rho - \frac{N-n+1}{N-n} \, |\nabla \log \rho|^2 + R + N(N-1) \geq 0\] 
at each point on $\Sigma$. Then 
\[-2 \, \Delta_\Sigma \log \check{\rho} - \frac{N-n+2}{N-n+1} \, |\nabla^\Sigma \log \check{\rho}|^2 + R_\Sigma + N(N-1) \geq 0\] 
at each point on $\Sigma$.
\end{proposition}

\textbf{Proof.} 
Using the PDE (\ref{Jacobi.operator.of.check.v}) and the inequality $\lambda_\infty \geq 0$, we obtain $\mathbb{L}_\Sigma \check{v} \geq 0$ at each point on $\Sigma$. In the next step, we use a crucial formula which originates in the work of Schoen and Yau \cite{Schoen-Yau1},\cite{Schoen-Yau2} and is closely related to the toric symmetrization technique of Gromov and Lawson (see \cite{Gromov-Lawson}, Sections 11 and 12). This gives 
\begin{align} 
\label{slicing.identity}
&-2 \, \Delta_\Sigma \log \check{\rho} - |\nabla^\Sigma \log \check{\rho}|^2 + R_\Sigma \notag \\ 
&+ 2 \, \Delta \log \rho + |\nabla \log \rho|^2 - R - |\nabla^\Sigma \log \check{v}|^2 - |h_\Sigma|^2 \\ 
&= 2 \, \check{v}^{-1} \, \mathbb{L}_\Sigma \check{v} \geq 0 \notag 
\end{align}
at each point on $\Sigma$ (see also \cite{Brendle-Hung}, Section 4). We distinguish two cases: 

\textit{Case 1:} Suppose first that $n=N$. In this case, our assumption implies that $\rho = 1$ and $R + N(N-1) \geq 0$. Using (\ref{slicing.identity}) and the identity $\check{\rho} = b_{n-2}^{-1} \, \check{v}$, we obtain 
\[-2 \, \Delta_\Sigma \log \check{\rho} - 2 \, |\nabla^\Sigma \log \check{\rho}|^2 + R_\Sigma + N(N-1) \geq 0\] 
at each point on $\Sigma$.

\textit{Case 2:} Suppose now that $n<N$. In this case, our assumption implies that 
\[-2 \, \Delta \log \rho - \frac{N-n+1}{N-n} \, |\nabla \log \rho|^2 + R + N(N-1) \geq 0\] 
at each point on $\Sigma$. Using (\ref{slicing.identity}), we obtain 
\begin{align} 
\label{consequence.of.slicing.identity}
&-2 \, \Delta_\Sigma \log \check{\rho} - |\nabla^\Sigma \log \check{\rho}|^2 + R_\Sigma + N(N-1) \notag \\ 
&- \frac{1}{N-n} \, |\nabla \log \rho|^2 - |\nabla^\Sigma \log \check{v}|^2 \geq 0 
\end{align}
at each point on $\Sigma$. Moreover, using the identity $\check{\rho} = b_{n-2}^{-1} \, \check{v} \, \rho$, we compute 
\begin{align} 
\label{gradient.terms}
&\frac{1}{N-n} \, |\nabla \log \rho|^2 + |\nabla^\Sigma \log \check{v}|^2 - \frac{1}{N-n+1} \, |\nabla^\Sigma \log \check{\rho}|^2 \notag \\ 
&= \frac{N-n}{N-n+1} \, \Big | \frac{1}{N-n} \, \nabla^\Sigma \log \rho - \nabla^\Sigma \log \check{v} \Big |^2 \\ 
&+ \frac{1}{N-n} \, \langle \nabla \log \rho,\nu_\Sigma \rangle^2 \geq 0 \notag
\end{align} 
at each point on $\Sigma$. Adding (\ref{consequence.of.slicing.identity}) and (\ref{gradient.terms}), we conclude that 
\[-2 \, \Delta_\Sigma \log \check{\rho} - \frac{N-n+2}{N-n+1} \, |\nabla^\Sigma \log \check{\rho}|^2 + R_\Sigma + N(N-1) \geq 0\] 
at each point on $\Sigma$. This completes the proof of Proposition \ref{scalar.curvature.of.hypersurface}. \\

\section{The conformal compactification and a foliation near infinity}

\label{foliation}

Throughout this section, we fix integers $N$ and $n$ such that $3 \leq n \leq N$. We define a flat metric $\gamma$ on $T^{n-1}$ by $\gamma = \sum_{k=0}^{n-2} b_k^2 \, d\theta_k \otimes d\theta_k$. Given a positive real number $r_0$, we define a hyperbolic metric $\bar{g}$ on $(r_0,\infty) \times T^{n-1}$ by $\bar{g} = r^{-2} \, dr \otimes dr + r^2 \, \gamma$. 

Let $(M,g)$ be a complete, connected, orientable Riemannian manifold of dimension $n$. We assume that there exists a bounded open domain $E \subset M$ with smooth boundary such that the complement $M \setminus E$ is diffeomorphic to $[r_0,\infty) \times T^{n-1}$. For every nonnegative integer $m$, we assume that 
\[|\bar{D}^m (g-\bar{g})|_{\bar{g}} \leq O(r^{-N})\] 
on $M \setminus E$, where $\bar{D}^m$ denotes the covariant derivative of order $m$ with respect to the hyperbolic metric $\bar{g}$. We further assume that the metric $g$ satisfies 
\[|g-\bar{g} - r^{2-N} \, Q|_{\bar{g}} \leq O(r^{-N-2\delta}),\] 
where $Q$ is a smooth symmetric $(0,2)$-tensor on $T^{n-1}$. Arguing as in Lemma \ref{derivatives.of.g}, we conclude that 
\[|\bar{D}^m(g-\bar{g} - r^{2-N} \, Q)|_{\bar{g}} \leq O(r^{-N-\delta})\] 
for every nonnegative integer $m$.

It is convenient to perform a change of variables and put $z = r^{-1}$. In the new coordinates, the hyperbolic metric takes the form $\bar{g} = z^{-2} \, (dz \otimes dz + \gamma)$. For abbrevation, we denote by $g_{\text{\rm flat}}$ the flat metric $dz \otimes dz + \gamma$. This gives 
\[|D_{\text{\rm flat}}^m(z^2 \, g - g_{\text{\rm flat}} - z^N \, Q)|_{g_{\text{\rm flat}}} \leq O(z^{N-m+\delta}),\] 
where $D_{\text{\rm flat}}^m$ denotes the covariant derivative of order $m$ with respect to the flat metric $g_{\text{\rm flat}}$. From this, it is easy to see that the conformal metric $\tilde{g} = z^2 \, g$ extends to a metric of class $C^N$ on a compact manifold $\tilde{M}$ with boundary. The manifold $\tilde{M}$ is referred to as the conformal compactification of $M$. The functions $z,\theta_0,\hdots,\theta_{n-2}$ extend smoothly to $\tilde{M}$. Moreover, $z=0$ on the boundary $\partial \tilde{M}$.

We next consider an interval $I \subset \mathbb{R}$ and a curve $\alpha: I \to M \setminus E$ satisfying 
\begin{equation} 
\label{geod.g}
D_s \dot{\alpha}(s) = -z^{-3} \, |dz|_g^2 \, \dot{\alpha}(s). 
\end{equation}
Every curve $\alpha$ satisfying (\ref{geod.g}) is a reparametrization of a geodesic. Let us consider the conformal metric $\tilde{g} = z^2 \, g$, and let $\tilde{D}$ denote the Levi-Civita connection with respect to the metric $\tilde{g}$. The equation (\ref{geod.g}) is equivalent to 
\begin{align} 
\label{geod.tilde.g}
\tilde{D}_s \dot{\alpha}(s) 
&= -z^{-1} \, |\dot{\alpha}(s)|_{\tilde{g}}^2 \, \tilde{\nabla} z \big |_{\alpha(s)} \notag \\ 
&+ 2 \, z^{-1} \, \big \langle \tilde{\nabla} z \big |_{\alpha(s)},\dot{\alpha}(s) \big \rangle_{\tilde{g}} \, \dot{\alpha}(s) \\ 
&- z^{-1} \, |dz|_{\tilde{g}}^2 \, \dot{\alpha}(s). \notag
\end{align}
Here, $\tilde{\nabla} z$ denotes the gradient of the function $z$ with respect to the metric $\tilde{g}$. In the next step, we define a vector field $\zeta$ along $\alpha$ by 
\[\zeta(s) = z^{-1} \, \big ( \dot{\alpha}(s) - \tilde{\nabla} z \big |_{\alpha(s)} \big ).\]
Then 
\begin{equation} 
\label{ode.1}
\dot{\alpha}(s) = \tilde{\nabla} z \big |_{\alpha(s)} + z \, \zeta(s). 
\end{equation}
Moreover, the identity (\ref{geod.tilde.g}) is equivalent to 
\begin{align} 
\label{ode.2}
\tilde{D}_s \zeta(s) 
&= -z^{-1} \, \sum_{k=1}^n (\tilde{D}^2 z)_{\alpha(s)} \big ( \tilde{\nabla} z \big |_{\alpha(s)},\tilde{e}_k \big ) \, \tilde{e}_k \notag \\ 
&- \sum_{k=1}^n (\tilde{D}^2 z)_{\alpha(s)}(\zeta(s),\tilde{e}_k) \, \tilde{e}_k \\ 
&- |\zeta(s)|_{\tilde{g}}^2 \; \tilde{\nabla} z \big |_{\alpha(s)} + \big \langle \tilde{\nabla} z \big |_{\alpha(s)},\zeta(s) \big \rangle_{\tilde{g}} \, \zeta(s), \notag
\end{align}
where $\{\tilde{e}_1,\hdots,\tilde{e}_n\}$ denotes a local orthonormal frame with respect to the metric $\tilde{g}$. The equations (\ref{ode.1}) and (\ref{ode.2}) give a system of two ODEs of first order for the pair $(\alpha,\zeta)$.

Recall that the metric $\tilde{g}$ is of class $C^N$ up to the boundary. Since the Hessian $\tilde{D}^2 z$ involves first derivatives of the metric, it is of class $C^{N-1}$ up to the boundary. Moreover, the Hessian $\tilde{D}^2 z$ vanishes along the boundary $\partial \tilde{M}$. Consequently, the vector field $\sum_{k=1}^n (\tilde{D}^2 z)(\tilde{\nabla} z,\tilde{e}_k) \, \tilde{e}_k$ is of class $C^{N-1}$ up to the boundary. Moreover, this vector field vanishes along the boundary $\partial \tilde{M}$. Therefore, the vector field $z^{-1} \sum_{k=1}^n (\tilde{D}^2 z)(\tilde{\nabla} z,\tilde{e}_k) \, \tilde{e}_k$ is of class $C^{N-2}$ up to the boundary.

\begin{proposition}
\label{well.posedness}
We can find small positive constants $z_*$ and $s_*$ with the following properties: 
\begin{itemize} 
\item Suppose that $q$ is a point in $\tilde{M} \cap \{z \leq z_*\}$ and $\xi \in T_q \tilde{M}$ is a tangent vector with $|\xi|_{\tilde{g}} \leq 1$. Then there is a unique solution $(\alpha(s),\zeta(s))$, $s \in [0,s_*]$, of the system (\ref{ode.1})--(\ref{ode.2}) with initial conditions $\alpha(0) = q$ and $\zeta(0) = \xi$. We define $\Psi_s(q,\xi) = \alpha(s)$ for each $s \in [0,s_*]$. 
\item For each $s \in [0,s_*]$, the map $\Psi_s$ is of class $C^{N-2}$ up to the boundary.
\end{itemize}
\end{proposition}

\textbf{Proof.} 
This follows from standard local existence theory for ODEs. \\

In the following, we assume that $s_* > 0$ has been chosen sufficiently small so that the map 
\[\partial \tilde{M} \times [0,s_*] \to \tilde{M}, \, (q,s) \mapsto \Psi_s(q,0)\]
is a diffeomorphism of class $C^{N-2}$. Consequently, we can find a small positive number $z_{\text{\rm fol}}$ and a map $\Xi: [0,z_{\text{\rm fol}}] \times T^{n-1} \to S^1$ of class $C^{N-2}$ with the following properties: 
\begin{itemize} 
\item We have $\Xi(0,\theta_0,\hdots,\theta_{n-2}) = \theta_{n-2}$ and $\frac{\partial}{\partial z} \Xi(0,\theta_0,\hdots,\theta_{n-2}) = 0$ for all points $(\theta_0,\hdots,\theta_{n-2}) \in T^{n-1}$. 
\item We have $\frac{\partial}{\partial \theta_{n-2}} \Xi(z,\theta_0,\hdots,\theta_{n-2}) \neq 0$ for each point $(z,\theta_0,\hdots,\theta_{n-2}) \in [0,z_{\text{\rm fol}}] \times T^{n-1}$. 
\item For each $t \in S^1$, the set $\{\Xi=t\} \subset [0,z_{\text{\rm fol}}] \times T^{n-1}$ can be written as a graph $\{\theta_{n-2} = G_t(z,\theta_0,\hdots,\theta_{n-3})\}$. The map 
\[[0,z_{\text{\rm fol}}] \times T^{n-2} \times S^1 \to S^1, \quad (z,\theta_0,\hdots,\theta_{n-3},t) \mapsto G_t(z,\theta_0,\hdots,\theta_{n-3})\]
is of class $C^{N-2}$. Moreover, $G_t(0,\theta_0,\hdots,\theta_{n-3}) = t$ for all points $(\theta_0,\hdots,\theta_{n-3}) \in T^{n-2}$ and all $t \in S^1$. 
\item For each $t \in S^1$ and each point $p \in [0,z_{\text{\rm fol}}] \times T^{n-1}$ with $\Xi=t$, there exists a point $q \in \partial \tilde{M} \cap \{\theta_{n-2}=t\}$ and a real number $s \in [0,s_*]$ such that $\Psi_s(q,0) = p$.
\end{itemize}
In the following, we put $r_{\text{\rm fol}} = z_{\text{\rm fol}}^{-1}$. By choosing $z_{\text{\rm fol}}$ sufficiently small, we can further arrange that the Hessian of the function $r$ with respect to the metric $g$ is positive definite in the region $\{r > r_{\text{\rm fol}}\}$.

For each $t \in S^1$, we denote by $\mathcal{Z}_t$ the set of all points in $[r_{\text{\rm fol}},\infty) \times T^{n-1}$ with $\Xi=t$. For each $t \in S^1$, $\mathcal{Z}_t$ is a hypersurface of class $C^{N-2}$. Moreover, $[r_{\text{\rm fol}},\infty) \times T^{n-1} = \bigcup_{t \in S^1} \mathcal{Z}_t$, and the sets $\mathcal{Z}_t$ are pairwise disjoint.

\begin{proposition} 
\label{tame.and.totally.geodesic}
Assume that $\Sigma$ is a properly embedded hypersurface in $M$ which is $t_*$-tame for some element $t_* \in S^1$. If $\Sigma$ is totally geodesic, then $\Sigma \cap \{r > r_{\text{\rm fol}}\} = \mathcal{Z}_{t_*}$. 
\end{proposition}

\textbf{Proof.} 
The proof consists of three steps.

\textit{Step 1:} We claim that $\mathcal{Z}_{t_*} \subset \Sigma$. To prove this, we fix an arbitrary point $q \in \partial \tilde{M} \cap \{\theta_{n-2}=t_*\}$. Since $\Sigma$ is $t_*$-tame, we can find a sequence of points $q^{(j)} \in \Sigma$ such that $d_{\tilde{g}}(q^{(j)},q) \to 0$ as $j \to \infty$. Let $z^{(j)} > 0$ denote the value of the function $z$ at the point $q^{(j)}$. Since $q \in \partial \tilde{M}$, it follows that $z^{(j)} \to 0$ as $j \to \infty$. Since $\Sigma$ is $t_*$-tame, we can find a sequence of vectors $\xi^{(j)} \in T_{q^{(j)}} M$ such that $\tilde{\nabla} z \big |_{q^{(j)}} + z^{(j)} \, \xi^{(j)} \in T_{q^{(j)}} \Sigma$ for each $j$ and $|\xi^{(j)}|_{\tilde{g}} \to 0$ as $j \to \infty$. We define $\alpha^{(j)}(s) = \Psi_s(q^{(j)},\xi^{(j)})$ for all $j$ and all $s \in [0,s_*]$. Since the map $\Psi_s$ is continuous up to the boundary, it follows that $\alpha^{(j)}(s) = \Psi_s(q^{(j)},\xi^{(j)}) \to \Psi_s(q,0)$ for each $s \in (0,s_*]$. For each $j$, the path $\alpha^{(j)}$ is a solution of the ODE (\ref{geod.g}) with initial conditions $\alpha^{(j)}(0) = q^{(j)} \in \Sigma$ and $\dot{\alpha}^{(j)}(0) = \tilde{\nabla} z \big |_{q^{(j)}} + z^{(j)} \, \xi^{(j)} \in T_{q^{(j)}} \Sigma$. Since $\Sigma$ is totally geodesic, it follows that $\alpha^{(j)}(s) \in \Sigma$ for all $j$ and all $s \in [0,s_*]$. Passing to the limit as $j \to \infty$, we conclude that $\Psi_s(q,0) \in \Sigma$ for each $s \in (0,s_*]$. Since $\mathcal{Z}_{t_*} \subset \big \{ \Psi_s(q,0): q \in \partial \tilde{M} \cap \{\theta_{n-2}=t\}, \, s \in (0,s_*] \big \}$, we conclude that $\mathcal{Z}_{t_*} \subset \Sigma$.

\textit{Step 2:} We claim that $\Sigma \cap \{r > r_{\text{\rm fol}}\}$ is connected. Since $\Sigma$ is tame, the set $\Sigma \cap \{r > r_{\text{\rm fol}}\}$ has at exactly one unbounded connected component. If the set $\Sigma \cap \{r > r_{\text{\rm fol}}\}$ has a bounded connected component, we consider a point on that connected component where the function $r$ attains its maximum. Since $\Sigma$ is totally geodesic and the Hessian of the function $r$ is positive definite in the region $\{r > r_{\text{\rm fol}}\}$, this leads to a contradiction. Thus, the set $\Sigma \cap \{r > r_{\text{\rm fol}}\}$ has exactly one unbounded connected component and no bounded connected components. 

\textit{Step 3:} Finally, we claim that $\mathcal{Z}_{t_*} = \Sigma \cap \{r > r_{\text{\rm fol}}\}$. In view of Step 1, the set $\mathcal{Z}_{t_*}$ is contained in $\Sigma \cap \{r > r_{\text{\rm fol}}\}$. It is easy to see that the set $\mathcal{Z}_{t_*}$ is both open and closed as a subset of $\Sigma \cap \{r > r_{\text{\rm fol}}\}$. Since the set $\Sigma \cap \{r > r_{\text{\rm fol}}\}$ is connected by Step 2, it follows that $\mathcal{Z}_{t_*} = \Sigma \cap \{r > r_{\text{\rm fol}}\}$. This completes the proof of Proposition \ref{tame.and.totally.geodesic}. \\

\begin{proposition}
\label{uniqueness} 
For each $\bar{t} \in S^1$, there is at most one properly embedded, connected, orientable hypersurface $\Sigma$ with the property that $\Sigma$ is totally geodesic and $\Sigma \cap \{r > 2r_{\text{\rm fol}}\} = \mathcal{Z}_{\bar{t}} \cap \{r > 2r_{\text{\rm fol}}\}$.
\end{proposition}

\textbf{Proof.} 
Suppose that $\Sigma$ and $\tilde{\Sigma}$ are two hypersurfaces with the required properties. Let $A$ denote the set of all points $p \in \Sigma$ with the property that $p \in \tilde{\Sigma}$ and $T_p \Sigma = T_p \tilde{\Sigma}$. Clearly, $A$ is a closed subset of $\Sigma$. Since $\Sigma$ and $\tilde{\Sigma}$ are totally geodesic, it is easy to see that $A$ is open as a subset of $\Sigma$. Finally, our assumptions imply that $A$ is non-empty. Since $\Sigma$ is connected, it follows that $A = \Sigma$. Thus, $\Sigma \subset \tilde{\Sigma}$. An analogous argument shows that $\tilde{\Sigma} \subset \Sigma$. This completes the proof of Proposition \ref{uniqueness}. \\

\section{Construction of barriers}

\label{barrier}

Throughout this section, we assume that $N$ and $n$ are integers satisfying $3 \leq n \leq N$ and $(M,g,\rho)$ is an $(N,n)$-dataset. Our goal is to construct a family of domains that are mean concave with respect to the conformal metric $\rho^{\frac{2}{n-1}} \, g$.

By the intermediate value theorem, we can find a real number $\hat{s}_N \in (2,\infty)$ such that 
\[(1-\hat{s}_N^{-2})^{-\frac{1}{2}} = 2^N \, (\hat{s}_N^{2-N} - \hat{s}_N^{1-N}).\] 

\begin{definition}
\label{definition.of.psi}
We define a function $\psi: (1,\infty) \to \mathbb{R}$ by 
\[\psi(s) = (1 - s^{-2})^{\frac{1}{2}}\] 
for $s \in (1,\hat{s}_N]$ and 
\[\psi(s) = (1-\hat{s}_N^{-2})^{\frac{1}{2}} + \frac{2^N}{N} \, (\hat{s}_N^{-N} - s^{-N}) - \frac{2^N}{N+1} \, (\hat{s}_N^{-N-1}-s^{-N-1})\] 
for $s \in [\hat{s}_N,\infty)$.
\end{definition}

\begin{lemma} 
\label{properties.of.psi}
The function $\psi$ is continuously differentiable and strictly monotone increasing. Moreover, $\psi$ takes values in the interval $(0,1+\frac{1}{N})$.
\end{lemma} 

\textbf{Proof.} 
It follows from our choice of $\hat{s}_N$ that the function $\psi$ is continuously differentiable. The derivative of $\psi$ is given by  
\[\psi'(s) = (1 - s^{-2})^{-\frac{1}{2}} \, s^{-3}\] 
for $s \in (1,\hat{s}_N]$ and 
\[\psi'(s) = 2^N\, (s^{-N-1} - s^{-N-2})\] 
for $s \in [\hat{s}_N,\infty)$. In particular, the function $\psi$ is strictly monotone increasing. We next observe that $\lim_{s \to 1} \psi(s) = 0$. Moreover, since $\hat{s}_N \in (2,\infty)$, we obtain 
\[\lim_{s \to \infty} \psi(s) = (1-\hat{s}_N^{-2})^{\frac{1}{2}} + \frac{2^N}{N} \, \hat{s}_N^{-N} - \frac{2^N}{N+1} \, \hat{s}_N^{-N-1} \leq 1+\frac{1}{N}.\] 
Therefore, the function $\psi$ takes values in the interval $(0,1+\frac{1}{N})$. This completes the proof of Lemma \ref{properties.of.psi}. \\

\begin{definition}
\label{definition.of.chi}
We define a function $\chi: (1,\infty) \setminus \{\hat{s}_N\} \to \mathbb{R}$ by 
\[\chi(s) = s^{2-N} \, \frac{d}{ds} \bigg ( \frac{s^N \, \psi'(s)}{(s^{-2} + s^2 \, \psi'(s)^2)^{\frac{1}{2}}} \bigg )\] 
for $s \in (1,\infty) \setminus \{\hat{s}_N\}$. 
\end{definition}

\begin{lemma}
\label{estimate.chi} 
We have $\chi(s) \geq s^{-N}$ for all $s \in (1,\infty) \setminus \{\hat{s}_N\}$. 
\end{lemma} 

\textbf{Proof.} 
We compute  
\[\frac{s^N \, \psi'(s)}{(s^{-2} + s^2 \, \psi'(s)^2)^{\frac{1}{2}}} = s^{N-2}\] 
for $s \in (1,\hat{s}_N)$ and 
\[\frac{s^N \, \psi'(s)}{(s^{-2} + s^2 \, \psi'(s)^2)^{\frac{1}{2}}} = \big( 2^{-2N} \, (1 - s^{-1})^{-2} + s^{2-2N} \big )^{-\frac{1}{2}}\] 
for $s \in (\hat{s}_N,\infty)$. This implies 
\[\chi(s) = (N-2) \, s^{-1}\] 
for $s \in (1,\hat{s}_N)$ and 
\begin{align*} 
\chi(s) 
&= \big( 2^{-2N} \, (1 - s^{-1})^{-2} + s^{2-2N} \big )^{-\frac{3}{2}} \\ 
&\cdot \big ( 2^{-2N} \, (1 - s^{-1})^{-3} + (N-1) \, s^{3-2N} \big ) \, s^{-N}
\end{align*} 
for $s \in (\hat{s}_N,\infty)$. Since $\hat{s}_N \in (2,\infty)$, we obtain 
\begin{align*} 
2^{-2N} \, (1 - s^{-1})^{-2} + s^{2-2N} 
&\leq 2 \, \big ( 2^{-3N} \, (1-s^{-1})^{-3} + s^{3-3N} \big )^{\frac{2}{3}} \\ 
&\leq 2 \, \big ( 2^{-3N} \, (1-s^{-1})^{-3} + 2^{-N} \, s^{3-2N} \big )^{\frac{2}{3}} 
\end{align*} 
for all $s \in (\hat{s}_N,\infty)$. From this, we deduce that $\chi(s) \geq 2^{N-\frac{3}{2}} \, s^{-N}$ for all $s \in (\hat{s}_N,\infty)$. This completes the proof of Lemma \ref{estimate.chi}. \\

\begin{definition} 
Let $\sigma \in (r_0,\infty)$ be sufficiently large and let $\bar{t} \in S^1$. We define a domain $\Omega_{\sigma,\bar{t}} \subset [\sigma,\infty) \times T^{n-1}$ by 
\[\Omega_{\sigma,\bar{t}} = \{b_{n-2} \, \sigma \, d_{S^1}(\theta_{n-2},\bar{t}) < \psi(\sigma^{-1} r)\},\] 
where $d_{S^1}$ denotes the Riemannian distance on $S^1$.
\end{definition} 

It follows from Lemma \ref{properties.of.psi} that  
\[\Omega_{\sigma,\bar{t}} \subset \{b_{n-2} \, \sigma \, d_{S^1}(\theta_{n-2},\bar{t}) < 1+\frac{1}{N}\}.\] 
Moreover, the boundary $\partial \Omega_{\sigma,\bar{t}}$ is a hypersurface of class $C^1$.

\begin{proposition} 
\label{nested.family.of.sets}
Let us fix an element $\bar{t} \in S^1$. Then the sets $\Omega_{\sigma,\bar{t}}$, $\sigma \in (r_0,\infty)$, form a decreasing family of sets.
\end{proposition} 

\textbf{Proof.} 
This follows immediately from Lemma \ref{properties.of.psi}. \\

\begin{proposition} 
\label{mean.curvature.of.Omega.bar.g}
Let $\bar{\nu}$ denote the outward-pointing unit normal vector field along $\partial \Omega_{\sigma,\bar{t}}$ with respect to the hyperbolic metric $\bar{g}$, and let $\bar{H}$ denote the mean curvature of $\partial \Omega_{\sigma,\bar{t}}$ with respect to the hyperbolic metric $\bar{g}$. Then $\bar{H} + (N-n) \, r^{-1} \, \langle \bar{\nabla} r,\bar{\nu} \rangle_{\bar{g}} = -\chi(\sigma^{-1} r)$ for $r \in (\sigma,\infty) \setminus \{\hat{s}_N \sigma\}$. 
\end{proposition}

\textbf{Proof.} 
The Hessian of the function $r$ with respect to the hyperbolic metric $\bar{g}$ is given by $\bar{D}^2 r = r \, \bar{g}$. We define a function $F: (r_0,\infty) \times T^{n-1} \to \mathbb{R}$ by $F = d_{S^1}(\theta_{n-2},\bar{t})$. Note that $F$ is smooth for $0 < F < \pi$. The Hessian of the function $F$ with respect to the hyperbolic metric $\bar{g}$ satisfies 
\[\bar{D}^2 F + r^{-1} \, (dr \otimes dF + dF \otimes dr) = 0\] 
for $0 < F < \pi$. In particular, $\bar{\Delta} F = 0$ for $0 < F < \pi$, where $\bar{\Delta}$ denotes the Laplacian with respect to the hyperbolic metric $\bar{g}$. We next observe that 
\[|b_{n-2} \, \sigma \, \bar{\nabla} F - \sigma^{-1} \, \psi'(\sigma^{-1} r) \, \bar{\nabla} r|_{\bar{g}}^2 = \sigma^2 r^{-2} + \sigma^{-2} r^2 \, \psi'(\sigma^{-1} r)^2,\] 
provided that $0 < F < \pi$ and $r \in (\sigma,\infty) \setminus \{\hat{s}_N \sigma\}$. Consequently, the outward-pointing unit normal vector field along $\partial \Omega_{\sigma,\bar{t}}$ is given by 
\[\bar{\nu} = \frac{b_{n-2} \, \sigma \, \bar{\nabla} F - \sigma^{-1} \, \psi'(\sigma^{-1} r) \, \bar{\nabla} r}{|b_{n-2} \, \sigma \, \bar{\nabla} F - \sigma^{-1} \, \psi'(\sigma^{-1} r) \, \bar{\nabla} r|_{\bar{g}}} = \frac{b_{n-2} \, \sigma \, \bar{\nabla} F - \sigma^{-1} \, \psi'(\sigma^{-1} r) \, \bar{\nabla} r}{(\sigma^2 r^{-2} + \sigma^{-2} r^2 \, \psi'(\sigma^{-1} r)^2)^{\frac{1}{2}}}\] 
for $r \in (\sigma,\infty) \setminus \{\hat{s}_N \sigma\}$. In particular, 
\[r^{-1} \, \langle \bar{\nabla} r,\bar{\nu} \rangle_{\bar{g}} = -\frac{\sigma^{-1} r \, \psi'(\sigma^{-1} r)}{(\sigma^2 r^{-2} + \sigma^{-2} r^2 \, \psi'(\sigma^{-1} r)^2)^{\frac{1}{2}}}\] 
and 
\[b_{n-2} \, \sigma \, \langle \bar{\nabla} F,\bar{\nu} \rangle_{\bar{g}} = \frac{\sigma^2 r^{-2}}{(\sigma^2 r^{-2} + \sigma^{-2} r^2 \, \psi'(\sigma^{-1} r)^2)^{\frac{1}{2}}}\] 
for $r \in (\sigma,\infty) \setminus \{\hat{s}_N \sigma\}$. The mean curvature of $\partial \Omega_{\sigma,\bar{t}}$ with respect to the hyperbolic metric $\bar{g}$ satisfies 
\begin{align*} 
&|b_{n-2} \, \sigma \, \bar{\nabla} F - \psi'(\sigma^{-1} r) \, \bar{\nabla} r|_{\bar{g}} \, \bar{H} \\ 
&= b_{n-2} \, \sigma \, \text{\rm tr}_{\partial \Omega_{\sigma,\bar{t}}}(\bar{D}^2 F) - \sigma^{-1} \, \psi'(\sigma^{-1} r) \, \text{\rm tr}_{\partial \Omega_{\sigma,\bar{t}}}(\bar{D}^2 r) \\ 
&- \sigma^{-2} \, \psi''(\sigma^{-1} r) \, \text{\rm tr}_{\partial \Omega_{\sigma,\bar{t}}}(dr \otimes dr) 
\end{align*} 
for $r \in (\sigma,\infty) \setminus \{\hat{s}_N \sigma\}$. Since $\bar{\Delta} F = 0$ for $0 < F < \pi$, it follows that 
\begin{align*} 
&(\sigma^2 r^{-2} + \sigma^{-2} r^2 \, \psi'(\sigma^{-1} r)^2)^{\frac{1}{2}} \, \bar{H} \\ 
&= -b_{n-2} \, \sigma \, (\bar{D}^2 F)(\bar{\nu},\bar{\nu}) - (n-1) \, \sigma^{-1} r \, \psi'(\sigma^{-1} r) \\ 
&- \sigma^{-2} \, \psi''(\sigma^{-1} r) \, (|\bar{\nabla} r|_{\bar{g}}^2 - \langle \bar{\nabla} r,\bar{\nu} \rangle_{\bar{g}}^2) \\ 
&= 2b_{n-2} \, \sigma r^{-1} \, \langle \bar{\nabla} r,\bar{\nu} \rangle_{\bar{g}} \, \langle \bar{\nabla} F,\bar{\nu} \rangle_{\bar{g}} - (n-1) \, \sigma^{-1} r \, \psi'(\sigma^{-1} r) \\ 
&- \sigma^{-2} r^2 \, \psi''(\sigma^{-1} r) \, (1 - r^{-2} \, \langle \bar{\nabla} r,\bar{\nu} \rangle_{\bar{g}}^2) \\ 
&= -\frac{2\sigma r^{-1} \, \psi'(\sigma^{-1} r)}{\sigma^2 r^{-2} + \sigma^{-2} r^2 \, \psi'(\sigma^{-1} r)^2} - (n-1) \, \sigma^{-1} r \, \psi'(\sigma^{-1} r) \\ 
&- \frac{\psi''(\sigma^{-1} r)}{\sigma^2 r^{-2} + \sigma^{-2} r^2 \, \psi'(\sigma^{-1} r)^2} 
\end{align*} 
for $r \in (\sigma,\infty) \setminus \{\hat{s}_N \sigma\}$. Consequently, 
\begin{align*} 
&(\sigma^2 r^{-2} + \sigma^{-2} r^2 \, \psi'(\sigma^{-1} r)^2)^{\frac{1}{2}} \, (\bar{H} + (N-n) \, r^{-1} \, \langle \bar{\nabla} r,\bar{\nu} \rangle_{\bar{g}}) \\ 
&= -\frac{2\sigma r^{-1} \, \psi'(\sigma^{-1} r)}{\sigma^2 r^{-2} + \sigma^{-2} r^2 \, \psi'(\sigma^{-1} r)^2} - (N-1) \, \sigma^{-1} r \, \psi'(\sigma^{-1} r) \\ 
&- \frac{\psi''(\sigma^{-1} r)}{\sigma^2 r^{-2} + \sigma^{-2} r^2 \, \psi'(\sigma^{-1} r)^2}  
\end{align*}
for $r \in (\sigma,\infty) \setminus \{\hat{s}_N \sigma\}$. On the other hand, it follows from Definition \ref{definition.of.chi} that 
\begin{align*} 
(s^{-2} + s^2 \, \psi'(s)^2)^{\frac{1}{2}} \, \chi(s) 
&= \frac{2s^{-1} \, \psi'(s)}{s^{-2} + s^2 \, \psi'(s)^2} + (N-1) \, s \, \psi'(s) \\ 
&+ \frac{\psi''(s)}{s^{-2} + s^2 \, \psi'(s)^2} 
\end{align*} 
for $s \in (1,\infty) \setminus \{\hat{s}_N\}$. This completes the proof of Proposition \ref{mean.curvature.of.Omega.bar.g}. \\

\begin{corollary} 
\label{mean.curvature.of.Omega.g}
We can find a large number $r_{\text{\rm barrier}} \in (r_0,\infty)$ with the following property. Assume that $\sigma \in [r_{\text{\rm barrier}},\infty)$ and $\bar{t} \in S^1$. Let $\nu$ denote the outward-pointing unit normal vector field along $\partial \Omega_{\sigma,\bar{t}}$ with respect to the metric $g$, and let $H$ denote the mean curvature of $\partial \Omega_{\sigma,\bar{t}}$ with respect to the metric $g$. Then $H + \rho^{-1} \, \langle \nabla \rho,\nu \rangle < 0$ at each point on $\partial \Omega_{\sigma,\bar{t}} \setminus \{r=\hat{s}_N \sigma\}$. 
\end{corollary} 

\textbf{Proof.} 
Let $\bar{\nu}$ denote the outward-pointing unit normal vector field along $\partial \Omega_{\sigma,\bar{t}}$ with respect to the hyperbolic metric $\bar{g}$, and let $\bar{H}$ denote the mean curvature of $\partial \Omega_{\sigma,\bar{t}}$ with respect to the hyperbolic metric $\bar{g}$. It follows from Lemma \ref{estimate.chi} and Proposition \ref{mean.curvature.of.Omega.bar.g} that 
\[\bar{H} + (N-n) \, r^{-1} \, \langle \bar{\nabla} r,\bar{\nu} \rangle_{\bar{g}} \leq -\sigma^N \, r^{-N}\] 
at each point on $\partial \Omega_{\sigma,\bar{t}} \setminus \{r=\hat{s}_N \sigma\}$. The second fundamental form of $\partial \Omega_{\sigma,\bar{t}}$ with respect to the hyperbolic metric $\bar{g}$ is uniformly bounded, and the higher order covariant derivatives of the second fundamental form with respect to $\bar{g}$ are bounded as well. Since $|g-\bar{g}|_{\bar{g}} \leq O(r^{-N})$ and $|\bar{D}(g-\bar{g})| \leq O(r^{-N})$, it follows that 
\[|H-\bar{H}| \leq C \, r^{-N}\] 
and 
\[|r^{-1} \, \langle \nabla r,\nu \rangle - r^{-1} \, \langle \bar{\nabla} r,\bar{\nu} \rangle_{\bar{g}}| \leq C \, r^{-N}\] 
at each point on $\partial \Omega_{\sigma,\bar{t}} \setminus \{r=\hat{s}_N \sigma\}$, where $C$ is independent of $\sigma$. Finally, 
\[|\rho^{-1} \, \langle \nabla \rho,\nu \rangle - (N-n) \, r^{-1} \, \langle \nabla r,\nu \rangle| \leq C \, r^{-N}\] 
at each point on $\partial \Omega_{\sigma,\bar{t}} \setminus \{r=\hat{s}_N \sigma\}$, where $C$ is independent of $\sigma$. Putting these facts together, we conclude that 
\begin{align*} 
H + \rho^{-1} \, \langle \nabla \rho,\nu \rangle 
&\leq \bar{H} + (N-n) \, r^{-1} \, \langle \bar{\nabla} r,\bar{\nu} \rangle_{\bar{g}} + C \, r^{-N} \\ 
&\leq -\sigma^N \, r^{-N} + C \, r^{-N} 
\end{align*}
at each point on $\partial \Omega_{\sigma,\bar{t}} \setminus \{r=\hat{s}_N \sigma\}$, where $C$ is independent of $\sigma$. This completes the proof of Corollary \ref{mean.curvature.of.Omega.g}. \\

\section{Existence of $(g,\rho)$-stationary hypersurfaces which are $(g,\rho,u)$-stable in the sense of Definition \ref{definition.stability}}

\label{existence.theory}

Throughout this section, we assume that $N$ and $n$ are integers satisfying $3 \leq n \leq N \leq 7$ and $(M,g,\rho)$ is an $(N,n)$-dataset. If $n=N$, we assume that $\rho=1$ and $R + N(N-1) \geq 0$ at each point in $M$. If $n<N$, we assume that  
\[-2 \, \Delta \log \rho - \frac{N-n+1}{N-n} \, |\nabla \log \rho|^2 + R + N(N-1) \geq 0\] 
at each point in $M$. 

As in Section \ref{properties.of.stable.hypersurfaces}, we assume that $u: T^{n-1} \to \mathbb{R}$ is a solution of the PDE 
\[\Delta_\gamma u + \frac{N}{2} \, \text{\rm tr}_\gamma(Q) + N \, P = \text{\rm constant}.\] 
The function $u$ is twice continuously differentiable with H\"older continuous second derivatives. 

Let $r_{\text{\rm fol}}$ denote the constant introduced in Section \ref{foliation}. We define an open domain $U$ by $U = M \setminus \{r \geq 8r_{\text{\rm fol}}\}$. Throughout this section, we fix an arbitrary point $p_* \in M \setminus \{r \geq 2r_{\text{\rm fol}}\}$. 

Our goal is to construct an orientable hypersurface passing through $p_*$ which is $(g,\rho)$-stationary and is $(g,\rho,u)$-stable in the sense of Definition \ref{definition.stability}. Our arguments are inspired by the work of Gang Liu \cite{Liu}. 

\begin{proposition}
\label{perturbed.metrics}
We can find a sequence of positive real numbers $\varepsilon_i \to 0$ and a sequence of Riemannian metrics $g^{(i)}$ with the following properties:
\begin{itemize}
\item $\frac{1}{2} \, g \leq g^{(i)} \leq 2g$ at each point in $M$.
\item $g^{(i)} = g$ at each point in $M \setminus U$.
\item $g^{(i)} \to g$ in $C^\infty(\bar{U})$.
\item If $n=N$, then $R_{g^{(i)}} + N(N-1) > 0$ at each point on $U \setminus B_{(M,g)}(p_*,\varepsilon_i)$. If $n<N$, then
\[-2 \, \Delta_{g^{(i)}} \log \rho - \frac{N-n+1}{N-n} \, |d\log \rho|_{g^{(i)}}^2 + R_{g^{(i)}} + N(N-1) > 0\]
at each point in $U \setminus B_{(M,g)}(p_*,\varepsilon_i)$.
\end{itemize}
\end{proposition}

\textbf{Proof.}
Let us fix a sequence of positive real numbers $\varepsilon_i \to 0$. In the following, we assume that $i$ is chosen sufficiently large. For each $i$, we can find a smooth function $\varphi_i: \bar{U} \to \mathbb{R}$ such that 
\begin{equation}
\label{pde.for.varphi_i.in.a.small.ball}
-(n-1) \, \Delta_g \varphi_i - (n-2) \, \langle d\log \rho,d\varphi_i \rangle_g + \varphi_i = \exp \Big ( -\frac{1}{\varepsilon_i^2 - d_{(M,g)}(p_*,x)^2} \Big ) 
\end{equation} 
at each point in $B_{(M,g)}(p_*,\varepsilon_i)$, 
\begin{equation}
\label{pde.for.varphi_i.outside.a.small.ball}
-(n-1) \, \Delta_g \varphi_i - (n-2) \, \langle d\log \rho,d\varphi_i \rangle_g + \varphi_i = 0 
\end{equation} 
at each point in $U \setminus B_{(M,g)}(p_*,\varepsilon_i)$, and $\varphi_i=0$ on $\partial U$. Note that $\varphi_i \to 0$ in $C^\infty(\bar{U})$. It follows from the strict maximum principle that $\varphi_i > 0$ at each point in $U$. For each $i$, we define a smooth function $\omega_i: M \to \mathbb{R}$ by
\[\omega_i = \begin{cases} \exp(-\varphi_i^{-1}) & \text{\rm on $U$} \\ 0 & \text{\rm on $M \setminus U$.} \end{cases}\]
Note that $\omega_i \to 0$ in $C^\infty(\bar{U})$. If $i$ is sufficiently large, then $|d\omega_i|_g^2 \leq \omega_i$ and $\Delta_g \omega_i \geq \omega_i \, \varphi_i^{-2} \, \Delta_g \varphi_i$ at each point in $U$. For each $i$, we define a conformal metric $g^{(i)}$ on $M$ by
\[g^{(i)} = (1+\omega_i)^{-1} \, g.\] 
Using the standard formula for the change of the scalar curvature under a conformal change of the metric (see \cite{Besse}, Theorem 1.159), we obtain 
\begin{align} 
\label{scal}
R_{g^{(i)}} 
&= (1+\omega_i) \, R_g + (n-1) \, \Delta_g \omega_i \notag \\ 
&- \frac{(n-1)(n+2)}{4} \, (1+\omega_i)^{-1} \, |d\omega_i|_g^2 \\ 
&\geq R_g + (n-1) \, \omega_i \, \varphi_i^{-2} \, \Delta_g \varphi_i - C \, \omega_i \notag
\end{align} 
at each point in $U$. Moreover, 
\begin{align} 
\label{Laplacian.rho}
-\Delta_{g^{(i)}} \log \rho 
&= -(1+\omega_i) \, \Delta_g \log \rho + \frac{n-2}{2} \, \langle d\log \rho,d\omega_i \rangle_g \notag \\ 
&\geq -\Delta_g \log \rho + \frac{n-2}{2} \, \omega_i \, \varphi_i^{-2} \, \langle d\log \rho,d\varphi_i \rangle_g - C \, \omega_i 
\end{align} 
at each point in $U$. We distinguish two cases: 

\textit{Case 1:} Suppose first that $n=N$. By assumption, $\rho=1$ and $R_g + N(N-1) \geq 0$. Moreover, (\ref{pde.for.varphi_i.outside.a.small.ball}) gives $-(N-1) \, \Delta_g \varphi_i + \varphi_i = 0$ at each point in $U \setminus B_{(M,g)}(p_*,\varepsilon_i)$. Using (\ref{scal}), we obtain 
\begin{align*} 
R_{g^{(i)}} + N(N-1) 
&\geq R_g + N(N-1) + (N-1) \, \omega_i \, \varphi_i^{-2} \, \Delta_g \varphi_i - C \, \omega_i \\ 
&\geq \omega_i \, \varphi_i^{-1} - C \, \omega_i 
\end{align*}
at each point in $U \setminus B_{(M,g)}(p_*,\varepsilon_i)$. If $i$ is sufficiently large, then the expression on the right hand side is strictly positive at each point in $U \setminus B_{(M,g)}(p_*,\varepsilon_i)$. 

\textit{Case 2:} Suppose now that $n<N$. By assumption, 
\begin{equation} 
\label{assumption}
-2 \, \Delta_g \log \rho - \frac{N-n+1}{N-n} \, |d\log \rho|_g^2 + R_g + N(N-1) \geq 0. 
\end{equation}
Using (\ref{pde.for.varphi_i.outside.a.small.ball}), (\ref{scal}), (\ref{Laplacian.rho}), and (\ref{assumption}), we obtain 
\begin{align*} 
&-2 \, \Delta_{g^{(i)}} \log \rho - \frac{N-n+1}{N-n} \, |d\log \rho|_{g^{(i)}}^2 + R_{g^{(i)}} + N(N-1) \\ 
&\geq -2 \, \Delta_g \log \rho - \frac{N-n+1}{N-n} \, |d\log \rho|_g^2 + R_g + N(N-1) \\ 
&+ (n-1) \, \omega_i \, \varphi_i^{-2} \, \Delta_g \varphi_i + (n-2) \, \omega_i \, \varphi_i^{-2} \, \langle d\log \rho,d\varphi_i \rangle_g - C \, \omega_i \\ 
&\geq \omega_i \, \varphi_i^{-1} - C \, \omega_i 
\end{align*}
at each point in $U \setminus B_{(M,g)}(p_*,\varepsilon_i)$. If $i$ is sufficiently large, then the expression on the right hand side is strictly positive at each point in $U \setminus B_{(M,g)}(p_*,\varepsilon_i)$. This completes the proof of Proposition \ref{perturbed.metrics}. \\

Let us consider an arbitrary sequence $r_j \to \infty$. For each $j$, we define $M^{(j)} = M \setminus \{r > r_j\}$. For each $j$ and each $t \in S^1$, we define $\Gamma_t^{(j)} = \partial M^{(j)} \cap \{\theta_{n-2}=t\}$. 

\begin{lemma}
\label{Gamma.bounds.a.hypersurface}
For each $j$ and each $t \in S^1$, $\Gamma_t^{(j)}$ bounds a compact, orientable hypersurface.
\end{lemma}

\textbf{Proof.}
By definition of an $(N,n)$-dataset, $\theta_{n-2}$ extends to a smooth map from $M$ to $S^1$. If $t \in S^1$ is a regular value of the map $\theta_{n-2}: M \to S^1$, then $M^{(j)} \cap \{\theta_{n-2}=t\}$ is a compact, orientable hypersurface with boundary $\Gamma_t^{(j)}$. This proves the assertion in the special case when $t$ is a regular value of the map $\theta_{n-2}: M \to S^1$. Since the set of regular values is dense, the assertion is true for each $t \in S^1$. This completes the proof of Lemma \ref{Gamma.bounds.a.hypersurface}. \\

Given positive integers $i,j$ and an element $t \in S^1$, we minimize the $(g^{(i)},\rho)$-area over all compact, orientable hypersurfaces $\Sigma \subset M^{(j)}$ with boundary $\Gamma_t^{(j)}$. Let $\mathcal{A}^{(i,j)}(t)$ denote the infimum of the $(g^{(i)},\rho)$-area in this class of hypersurfaces. It is easy to see that the function $t \mapsto \mathcal{A}^{(i,j)}(t)$ is continuous. Given positive integers $i,j$, we minimize the function 
\begin{equation} 
\label{modified.area}
t \mapsto \mathcal{A}^{(i,j)}(t) + \int_{T^{n-2} \times \{t\}} u \, d\text{\rm vol}_\gamma 
\end{equation}
over all $t \in S^1$. Given positive integers $i,j$, we can find an element $t_{i,j} \in S^1$ where the function (\ref{modified.area}) attains its minimum. Moreover, given positive integers $i,j$, we can find a compact, orientable hypersurface $\Sigma^{(i,j)}$ with boundary $\Gamma_{t_{i,j}}^{(j)}$ such that the $(g^{(i)},\rho)$-area of $\Sigma^{(i,j)}$ is equal to $\mathcal{A}^{(i,j)}(t_{i,j})$. Note that $\Sigma^{(i,j)}$ is connected. 

The minimization problem above can be viewed as a hybrid between a Plateau problem and a free boundary problem. It is inspired by the work of Schoen and Yau on the positive mass theorem (see \cite{Schoen}, Section 4, and \cite{Eichmair-Huang-Lee-Schoen}). 

\begin{proposition}
\label{upper.bound.for.area.1}
The $(g^{(i)},\rho)$-area of $\Sigma^{(i,j)}$ is bounded from above by 
\[(2\pi)^{n-2} \, \Big ( \prod_{k=0}^{n-3} b_k \Big ) \, \frac{r_j^{N-2}}{N-2} + C.\] 
The constant $C$ is independent of $i$ and $j$.
\end{proposition}

\textbf{Proof.} 
By definition of an $(N,n)$-dataset, $\theta_{n-2}$ extends to a smooth map from $M$ to $S^1$. Let us fix an element $\bar{t} \in S^1$ so that $\bar{t}$ is a regular value of the map $\theta_{n-2}: M \to S^1$. Then $M^{(j)} \cap \{\theta_{n-2}=\bar{t}\}$ is a compact, orientable hypersurface with boundary $\Gamma_{\bar{t}}^{(j)}$. This implies 
\[\mathcal{A}^{(i,j)}(\bar{t}) \leq \int_{M^{(j)} \cap \{\theta_{n-2}=\bar{t}\}} \rho \, d\text{\rm vol}_{g^{(i)}} \leq (2\pi)^{n-2} \, \Big ( \prod_{k=0}^{n-3} b_k \Big ) \, \frac{r_j^{N-2}}{N-2} + C,\] 
where $C$ is independent of $i$ and $j$. On the other hand, since the function (\ref{modified.area}) attains its minimum at $t_{i,j}$, we know that 
\[\int_{\Sigma^{(i,j)}} \rho \, d\text{\rm vol}_{g^{(i)}} = \mathcal{A}^{(i,j)}(t_{i,j}) \leq \mathcal{A}^{(i,j)}(\bar{t}) + C,\] 
where $C$ is independent of $i$ and $j$. Putting these facts together, we conclude that 
\[\int_{\Sigma^{(i,j)}} \rho \, d\text{\rm vol}_{g^{(i)}} \leq (2\pi)^{n-2} \, \Big ( \prod_{k=0}^{n-3} b_k \Big ) \, \frac{r_j^{N-2}}{N-2} + C,\] 
where $C$ is independent of $i$ and $j$. This completes the proof of Proposition \ref{upper.bound.for.area.1}. \\

\begin{proposition}
\label{lower.bound.for.area}
We can find a large number $r_1 \geq \max \{2r_0,16r_{\text{\rm fol}}\}$ and a large constant $C$ with the following property. If $\bar{r} \in (r_1,r_j)$, then the $(g,\rho)$-area of $\Sigma^{(i,j)} \cap \{r > \bar{r}\}$ is bounded from below by 
\[(2\pi)^{n-2} \, \Big ( \prod_{k=0}^{n-3} b_k \Big ) \, \frac{r_j^{N-2}-\bar{r}^{N-2}}{N-2} - C \, \bar{r}^{-2}.\] 
\end{proposition}

\textbf{Proof.} 
For each $k \in \{0,1,\hdots,n-3\}$, we denote by $\Theta_k$ the pull-back of the volume form on $S^1$ under the map $\theta_k: (r_0,\infty) \times T^{n-1} \to S^1$. Here, we assume that the volume form on $S^1$ is normalized to have integral $2\pi$. For each $k \in \{0,1,\hdots,n-3\}$, $\Theta_k$ is a closed one-form on $(r_0,\infty) \times T^{n-1}$.

We can find a large constant $C_1$ and a large number $r_1 \geq \max \{2r_0,16r_{\text{\rm fol}}\}$ with the following properties: 
\begin{itemize} 
\item $(1 - \frac{C_1}{n} \, r^{-N}) \, r^{N-n} \leq \rho$ at each point in $[r_1,\infty) \times T^{n-1}$. 
\item $(1 - \frac{C_1}{n} \, r^{-N}) \, r^{-1} \, |dr|_g \leq 1$ at each point in $[r_1,\infty) \times T^{n-1}$.
\item For each $k \in \{0,1,\hdots,n-3\}$, we have $(1 - \frac{C_1}{n} \, r^{-N}) \, r \, |\Theta_k|_g \leq b_k^{-1}$ at each point in $[r_1,\infty) \times T^{n-1}$. 
\item $C_1 \, r_1^{-N} < 1$.
\end{itemize}
Then 
\[\Big ( 1 - \frac{C_1}{n} \, r^{-N} \Big )^n \, r^{N-3} \, |dr|_g \, \prod_{k=0}^{n-3} |\Theta_k|_g \leq \Big ( \prod_{k=0}^{n-3} b_k \Big )^{-1} \, \rho\] 
at each point in $[r_1,\infty) \times T^{n-1}$. Note that $1 - C_1 \, r^{-N} \leq (1 - \frac{C_1}{n} \, r^{-N})^n$ for $r \geq r_1$. Consequently, 
\begin{align*} 
&\int_{\Sigma^{(i,j)} \cap \{r > \bar{r}\}} (1 - C_1 \, r^{-N}) \, r^{N-3} \, dr \wedge \Theta_0 \wedge \Theta_1 \wedge \hdots \wedge \Theta_{n-3} \\ 
&\leq \Big ( \prod_{k=0}^{n-3} b_k \Big )^{-1} \, \int_{\Sigma^{(i,j)} \cap \{r > \bar{r}\}} \rho \, d\text{\rm vol}_g 
\end{align*} 
for each $\bar{r} \in (r_1,r_j)$. Using Stokes theorem, we compute 
\begin{align*} 
&\int_{\Sigma^{(i,j)} \cap \{r > \bar{r}\}} (1 - C_1 \, r^{-N}) \, r^{N-3} \, dr \wedge \Theta_0 \wedge \Theta_1 \wedge \hdots \wedge \Theta_{n-3} \\ 
&= \int_{\Sigma^{(i,j)} \cap \{r > \bar{r}\}} d \Big [ \Big ( \frac{r^{N-2}-\bar{r}^{N-2}}{N-2} + C_1 \, \frac{r^{-2}-\bar{r}^{-2}}{2} \Big ) \, \Theta_0 \wedge \Theta_1 \wedge \hdots \wedge \Theta_{n-3} \Big ] \\ 
&= \int_{\Gamma_{t_{i,j}}^{(j)}} \Big ( \frac{r^{N-2}-\bar{r}^{N-2}}{N-2} + C_1 \, \frac{r^{-2}-\bar{r}^{-2}}{2} \Big ) \, \Theta_0 \wedge \Theta_1 \wedge \hdots \wedge \Theta_{n-3} \\ 
&= (2\pi)^{n-2} \, \Big ( \frac{r_j^{N-2}-\bar{r}^{N-2}}{N-2} + C_1 \, \frac{r_j^{-2}-\bar{r}^{-2}}{2} \Big ) 
\end{align*} 
whenever $\bar{r} \in (r_1,r_j)$ is a regular value of the function $r|_{\Sigma^{(i,j)}}$. Putting these facts together, we conclude that 
\[(2\pi)^{n-2} \, \Big ( \frac{r_j^{N-2}-\bar{r}^{N-2}}{N-2} + C_1 \, \frac{r_j^{-2}-\bar{r}^{-2}}{2} \Big ) \leq \Big ( \prod_{k=0}^{n-3} b_k \Big )^{-1} \, \int_{\Sigma^{(i,j)} \cap \{r > \bar{r}\}} \rho \, d\text{\rm vol}_g\] 
whenever $\bar{r} \in (r_1,r_j)$ is a regular value of the function $r|_{\Sigma^{(i,j)}}$. Since the set of regular values is dense, the assertion is true for each $\bar{r} \in (r_1,r_j)$. This completes the proof of Proposition \ref{lower.bound.for.area}. \\

\begin{corollary}
\label{upper.bound.for.area.2}
If $\bar{r} \in (r_1,r_j)$, then the $(g^{(i)},\rho)$-area of $\Sigma^{(i,j)} \setminus \{r > \bar{r}\}$ is bounded from above by 
\[(2\pi)^{n-2} \, \Big ( \prod_{k=0}^{n-3} b_k \Big ) \, \frac{\bar{r}^{N-2}}{N-2} + C.\] 
The constant $C$ is independent of $i$ and $j$.
\end{corollary}

\textbf{Proof.} 
This follows by combining Proposition \ref{upper.bound.for.area.1} and Proposition \ref{lower.bound.for.area}. \\

In the next step, we establish a curvature bound for the hypersurface $\Sigma^{(i,j)}$.

\begin{proposition} 
\label{curvature.bound}
We have $|h_{\Sigma^{(i,j)}}| \leq C$ at each point in $\Sigma^{(i,j)} \setminus \{r > 2^{-\frac{1}{4N}} \, r_j\}$. The constant $C$ is independent of $i$ and $j$. Moreover, the higher order covariant derivatives of the second fundamental form of $\Sigma^{(i,j)}$ are uniformly bounded on the set $\Sigma^{(i,j)} \setminus \{r > 2^{-\frac{1}{2N}} \, r_j\}$.
\end{proposition}

\textbf{Proof.} 
Our assumptions imply that the injectivity radius of $(M,g)$ is bounded from below by a positive constant. Let us fix a real number $\alpha_0 \in (0,\frac{1}{8N})$ so that $\alpha_0$ is smaller than the injectivity radius of $(M,g)$. Let us consider an arbitrary point $q \in \Sigma^{(i,j)} \setminus \{r > 2^{-\frac{1}{4N}} \, r_j\}$. Then the ball $B_{(M,g)}(q,\alpha_0)$ is disjoint from $\partial \Sigma^{(i,j)}$. By Sard's lemma, we can find a real number $\alpha \in (\frac{\alpha_0}{4},\frac{\alpha_0}{2})$ such that $\Sigma^{(i,j)}$ intersects $\partial B_{(M,g)}(q,\alpha)$ transversally. 

By assumption, the higher order derivatives of the function $\log \rho$ satisfy $|\bar{D}^m \log \rho|_{\bar{g}} \leq C(m)$. In particular, $C^{-1} \, \rho(q) \leq \rho \leq C \, \rho(q)$ in $B_{(M,g)}(q,\alpha)$. We next observe that $B_{(M,g)}(q,\alpha)$ is diffeomorphic to a ball in $\mathbb{R}^n$. Since $\Sigma^{(i,j)}$ minimizes the $(g^{(i)},\rho)$-area, it follows that the $(g^{(i)},\rho)$-area of each connected component of $\Sigma^{(i,j)} \cap B_{(M,g)}(q,\alpha)$ is bounded from above by the $(g^{(i)},\rho)$-area of $\partial B_{(M,g)}(q,\alpha)$. Consequently, the $(g^{(i)},\rho)$-area of each connected component of $\Sigma^{(i,j)} \cap B_{(M,g)}(q,\alpha)$ is bounded from above by $C \, \rho(q)$, where $C$ is independent of $i$, $j$, and $q$. 

Using the Schoen-Simon curvature estimate (see \cite{Schoen-Simon}, Corollary 1, and \cite{Bellettini}, Theorem 2), we conclude that $|h_{\Sigma^{(i,j)}}| \leq C$ at the point $q$, where $C$ is independent of $i$, $j$, and $q$. To summarize, we have shown that the second fundamental form of $\Sigma^{(i,j)}$ is uniformly bounded on $\Sigma^{(i,j)} \setminus \{r > 2^{-\frac{1}{4N}} \, r_j\}$. Using standard interior estimates, we obtain bounds for the higher order covariant derivatives of the second fundamental form of $\Sigma^{(i,j)}$ on the set $\Sigma^{(i,j)} \setminus \{r > 2^{-\frac{1}{2N}} \, r_j\}$. This completes the proof of Proposition \ref{curvature.bound}. \\

\begin{proposition}
\label{mean.concavity}
Let $r_{\text{\rm barrier}}$ be chosen as in Corollary \ref{mean.curvature.of.Omega.g}. Suppose that $\sigma \geq \max \{r_{\text{\rm barrier}},8r_{\text{\rm fol}}\}$ and $\bar{t} \in S^1$. Then the domain $\Omega_{\sigma,\bar{t}}$ is strictly mean concave with respect to the conformal metric $\rho^{\frac{2}{n-1}} \, g^{(i)}$.
\end{proposition}

\textbf{Proof.} 
By Proposition \ref{perturbed.metrics}, the metric $g^{(i)}$ agrees with the metric $g$ in the region $\{r \geq 8r_{\text{\rm fol}}\}$. Therefore, the assertion follows from Corollary \ref{mean.curvature.of.Omega.g}. \\

\begin{proposition}
\label{sliding.barrier}
Let $r_{\text{\rm barrier}}$ be chosen as in Corollary \ref{mean.curvature.of.Omega.g}. Suppose that $\bar{\sigma} \geq \max \{r_{\text{\rm barrier}},8r_{\text{\rm fol}}\}$ and $\bar{t} \in S^1$. Moreover, suppose that $\partial \Sigma^{(i,j)}$ is disjoint from $\Omega_{\bar{\sigma},\bar{t}}$. Then $\Sigma^{(i,j)}$ is disjoint from $\Omega_{\bar{\sigma},\bar{t}}$. 
\end{proposition}

\textbf{Proof.} 
By assumption, $\partial \Sigma^{(i,j)}$ is disjoint from $\Omega_{\bar{\sigma},\bar{t}}$. By Proposition \ref{nested.family.of.sets}, $\partial \Sigma^{(i,j)}$ is disjoint from $\Omega_{\sigma,\bar{t}}$ for each $\sigma \in [\bar{\sigma},\infty)$. Moreover, if $\sigma$ is sufficiently large depending on $j$, then $\Sigma^{(i,j)}$ is disjoint from $\Omega_{\sigma,\bar{t}}$. Using Proposition \ref{mean.concavity} and a sliding barrier argument, we conclude that $\Sigma^{(i,j)}$ is disjoint from $\Omega_{\sigma,\bar{t}}$ for each $\sigma \in [\bar{\sigma},\infty)$. This completes the proof of Proposition \ref{sliding.barrier}. \\

\begin{corollary} 
\label{barrier.estimate}
We can find a large constant $r_2 \geq 16r_{\text{\rm fol}}$ and a large constant $L$ with the following properties: 
\begin{itemize}
\item $L \, r_2^{-N} \leq \frac{\pi}{2}$. 
\item If $j$ is sufficiently large, then $\Sigma^{(i,j)} \cap \{r \geq r_2\} \subset \{d_{S^1}(\theta_{n-2},t_{i,j}) \leq L \, r^{-N}\}$. 
\end{itemize}
Note that $r_2$ and $L$ are independent of $i$ and $j$.
\end{corollary}

\textbf{Proof.} 
This follows directly from Proposition \ref{sliding.barrier}. \\

In view of Corollary \ref{barrier.estimate}, we have $\Sigma^{(i,j)} \cap \{r \geq r_2\} \subset \{d_{S^1}(\theta_{n-2},t_{i,j}) \leq \frac{\pi}{2}\}$. For each pair $(i,j)$, we choose a smooth function $F^{(i,j)}: [r_2,\infty) \times T^{n-1} \to \mathbb{R}$ with the property that $F^{(i,j)}$ equals the signed distance on $S^1$ from $\theta_{n-2}$ to $t_{i,j}$ when $d_{S^1}(\theta_{n-2},t_{i,j}) \leq \frac{\pi}{2}$. In particular, $|F^{(i,j)}| = d_{S^1}(\theta_{n-2},t_{i,j})$ whenever $d_{S^1}(\theta_{n-2},t_{i,j}) \leq \frac{\pi}{2}$.

\begin{proposition}
\label{higher.derivatives.of.height.function}
For every positive integer $m$, we have $|D^{\Sigma^{(i,j)},m} F^{(i,j)}| \leq C(m) \, r^{-N}$ at each point in $\Sigma^{(i,j)} \cap \{2r_2 \leq r \leq 2^{-\frac{1}{N}} \, r_j\}$. Here, $D^{\Sigma^{(i,j)},m}$ denotes the covariant derivative of order $m$ with respect to the induced metric $g|_{\Sigma^{(i,j)}}$. Note that $C(m)$ may depend on $m$, but is independent of $i$ and $j$.
\end{proposition}

\textbf{Proof.} 
The Hessian of the function $F^{(i,j)}$ with respect to the hyperbolic metric $\bar{g}$ satisfies 
\begin{equation} 
\label{identity.1}
\bar{D}^2 F^{(i,j)} + r^{-1} \, (dr \otimes dF^{(i,j)} + dF^{(i,j)} \otimes dr) = 0 
\end{equation} 
for $d_{S^1}(\theta_{n-2},t_{i,j}) \leq \frac{\pi}{2}$. Moreover, 
\begin{equation} 
\label{identity.2}
\langle dr,dF^{(i,j)} \rangle_{\bar{g}} = 0 
\end{equation}
for $d_{S^1}(\theta_{n-2},t_{i,j}) \leq \frac{\pi}{2}$. 

We define a symmetric $(0,2)$-tensor $B^{(i,j)}$ on $[r_2,\infty) \times T^{n-1}$ by 
\[B^{(i,j)} = D^2 F^{(i,j)} + r^{-1} \, (dr \otimes dF^{(i,j)} + dF^{(i,j)} \otimes dr),\] 
where $D^2 F^{(i,j)}$ denotes the Hessian of the function $F^{(i,j)}$ with respect to the metric $g$. Moreover, we define a function $\beta^{(i,j)}: [r_2,\infty) \times T^{n-1} \to \mathbb{R}$ by 
\[\beta^{(i,j)} = \rho^{-1} \, \langle d\rho,dF^{(i,j)} \rangle,\] 
where the inner product is computed with respect to the metric $g$. Using (\ref{identity.1}) and (\ref{identity.2}), we obtain 
\[B^{(i,j)} = D^2 F^{(i,j)} - \bar{D}^2 F^{(i,j)}\] 
and 
\[\beta^{(i,j)} = \rho^{-1} \, \langle d\rho,dF^{(i,j)} \rangle - (N-n) \, r^{-1} \, \langle dr,dF^{(i,j)} \rangle_{\bar{g}}\]
for $d_{S^1}(\theta_{n-2},t_{i,j}) \leq \frac{\pi}{2}$. This implies 
\begin{equation} 
\label{estimate.1}
|D^m B^{(i,j)}| \leq C(m) \, r^{-N-1} 
\end{equation} 
and 
\begin{equation} 
\label{estimate.2} 
|D^m \beta^{(i,j)}| \leq C(m) \, r^{-N-1} 
\end{equation} 
for $d_{S^1}(\theta_{n-2},t_{i,j}) \leq \frac{\pi}{2}$. Here, $D^m$ denotes the covariant derivative of order $m$ with respect to the metric $g$ on $[r_2,\infty) \times T^{n-1}$. 

In the next step, we consider the restriction of the function $F^{(i,j)}$ to the hypersurface $\Sigma^{(i,j)} \cap \{r \geq r_2\}$. Corollary \ref{barrier.estimate} implies that $|F^{(i,j)}| = d_{S^1}(\theta_{n-2},t_{i,j}) \leq L \, r^{-N}$ at each point in $\Sigma^{(i,j)} \cap \{r \geq r_2\}$. Moreover, we compute 
\begin{align}
\label{Laplacian.of.F}
&\Delta_{\Sigma^{(i,j)}} F^{(i,j)} + 2 r^{-1} \, \langle \nabla^{\Sigma^{(i,j)}} r,\nabla^{\Sigma^{(i,j)}} F^{(i,j)} \rangle \notag \\ 
&= \text{\rm tr}_{\Sigma^{(i,j)}}(B^{(i,j)}) - H_{\Sigma^{(i,j)}} \, \langle \nabla F^{(i,j)},\nu_{\Sigma^{(i,j)}} \rangle 
\end{align} 
and 
\begin{align} 
\label{directional.derivative.of.F}
&\rho^{-1} \, \langle \nabla^{\Sigma^{(i,j)}} \rho,\nabla^{\Sigma^{(i,j)}} F^{(i,j)} \rangle \notag \\ 
&= \beta^{(i,j)} - \rho^{-1} \, \langle \nabla \rho,\nu_{\Sigma^{(i,j)}} \rangle \, \langle \nabla F^{(i,j)},\nu_{\Sigma^{(i,j)}} \rangle 
\end{align} 
at each point on $\Sigma^{(i,j)} \cap \{r \geq r_2\}$. In the next step, we add (\ref{Laplacian.of.F}) and (\ref{directional.derivative.of.F}). Since $\Sigma^{(i,j)}$ is $(g^{(i)},\rho)$-stationary and $g^{(i)} = g$ in the region $\{r \geq r_2\}$, it follows that 
\begin{align} 
\label{pde.for.F}
&\Delta_{\Sigma^{(i,j)}} F^{(i,j)} + 2 r^{-1} \, \langle \nabla^{\Sigma^{(i,j)}} r,\nabla^{\Sigma^{(i,j)}} F^{(i,j)} \rangle \notag \\ 
&+ \rho^{-1} \, \langle \nabla^{\Sigma^{(i,j)}} \rho,\nabla^{\Sigma^{(i,j)}} F^{(i,j)} \rangle \\ 
&= \text{\rm tr}_{\Sigma^{(i,j)}}(B^{(i,j)}) + \beta^{(i,j)} \notag
\end{align} 
at each point on $\Sigma^{(i,j)} \cap \{r \geq r_2\}$. Using (\ref{estimate.1}), (\ref{estimate.2}), and Proposition \ref{curvature.bound}, we obtain 
\begin{equation} 
\label{bound.for.source.term}
\big | D^{\Sigma^{(i,j)},m} \big ( \text{\rm tr}_{\Sigma^{(i,j)}}(B^{(i,j)}) + \beta^{(i,j)} \big ) \big | \leq C(m) \, r^{-N-1} 
\end{equation}
at each point on $\Sigma^{(i,j)} \cap \{r \geq r_2\}$. Here, $D^{\Sigma^{(i,j)},m}$ denotes the covariant derivative of order $m$ with respect to the metric $g|_{\Sigma^{(i,j)}}$. 

Suppose now that $q$ is a point in $\Sigma^{(i,j)} \cap \{2r_2 \leq r \leq 2^{-\frac{1}{N}} \, r_j\}$. By Proposition \ref{curvature.bound}, we control the geometry of $\Sigma^{(i,j)}$ in a ball around $q$ of some fixed radius. Moreover, we know that $|F^{(i,j)}| \leq L \, r^{-N}$ at each point in $\Sigma^{(i,j)} \cap \{r \geq r_2\}$. Using (\ref{pde.for.F}), (\ref{bound.for.source.term}), and standard interior estimates for elliptic PDE, we conclude that $|D^{\Sigma^{(i,j)},m} F^{(i,j)}| \leq C(m) \, r^{-N}$ at the point $q$. This completes the proof of Proposition \ref{higher.derivatives.of.height.function}. \\

\begin{corollary}
\label{bound.for.slope}
We have $|(\frac{\partial}{\partial \theta_{n-2}})^{\text{\rm tan}}| \leq C \, r^{2-N}$ at each point in $\Sigma^{(i,j)} \cap \{2r_2 \leq r \leq 2^{-\frac{1}{N}} \, r_j\}$. Here, $C$ is a large constant which is independent of $i$ and $j$.
\end{corollary}

\textbf{Proof.} 
It follows from Proposition \ref{higher.derivatives.of.height.function} that $|(\nabla F^{(i,j)})^{\text{\rm tan}}| \leq C \, r^{-N}$ at each point in $\Sigma^{(i,j)} \cap \{2r_2 \leq r \leq 2^{-\frac{1}{N}} \, r_j\}$. On the other hand, 
\[\Big | \nabla F^{(i,j)} - b_{n-2}^{-2} \, r^{-2} \, \frac{\partial}{\partial \theta_{n-2}} \Big | \leq C \, r^{-N-1}\] 
for $d_{S^1}(\theta_{n-2},t_{i,j}) \leq \frac{\pi}{2}$. Putting these facts together, the assertion follows. \\

\begin{corollary}
\label{transversality}
We can find a large number $r_3 \geq 2r_2$ such that $\frac{\partial}{\partial \theta_{n-2}} \notin T\Sigma^{(i,j)}$ at each point in $\Sigma^{(i,j)} \cap \{r_3 \leq r \leq 2^{-\frac{1}{N}} \, r_j\}$. Note that $r_3$ is independent of $i$ and $j$.
\end{corollary}

\begin{proposition} 
\label{graphical}
Let us consider the map $\pi^{(i,j)}: \Sigma^{(i,j)} \cap \{r_3 < r < 2^{-\frac{1}{N}} \, r_j\} \to (r_3,2^{-\frac{1}{N}} \, r_j) \times T^{n-2}$ which maps a point $(r,\theta_0,\hdots,\theta_{n-3},\theta_{n-2}) \in \Sigma^{(i,j)} \cap \{r_3 < r < 2^{-\frac{1}{N}} \, r_j\}$ to $(r,\theta_0,\hdots,\theta_{n-3}) \in (r_3,2^{-\frac{1}{N}} \, r_j) \times T^{n-2}$. If $j$ is sufficiently large, then the map $\pi^{(i,j)}$ is bijective.
\end{proposition} 

\textbf{Proof.} 
In view of Corollary \ref{transversality}, the differential of $\pi^{(i,j)}$ is invertible at each point in $\Sigma^{(i,j)} \cap \{r_3 < r < 2^{-\frac{1}{N}} \, r_j\}$. In particular, each point in $(r_3,2^{-\frac{1}{N}} \, r_j) \times T^{n-2}$ has the same number of pre-images under the map $\pi^{(i,j)}$. There are three possibilities: 

\textit{Case 1:} Each point in $(r_3,2^{-\frac{1}{N}} \, r_j) \times T^{n-2}$ has no pre-images under the map $\pi^{(i,j)}$. In this case, the hypersurface $\Sigma^{(i,j)}$ is disjoint from the region $\{r_3 < r < 2^{-\frac{1}{N}} \, r_j\}$. Consequently, $\Gamma_{t_{i,j}}^{(j)}$ bounds a hypersurface in $[r_3,\infty) \times T^{n-1}$, which is impossible.

\textit{Case 2:} Each point in $(r_3,2^{-\frac{1}{N}} \, r_j) \times T^{n-2}$ has exactly one pre-image under the map $\pi^{(i,j)}$. In this case, we are done.

\textit{Case 3:} Each point in $(r_3,2^{-\frac{1}{N}} \, r_j) \times T^{n-2}$ has at least two pre-images under the map $\pi^{(i,j)}$. We can find a large constant $C_2$ and a large number $r_4 \geq r_3$ such that $C_2 \, r_4^{-N} \leq \frac{1}{2}$ and 
\begin{align*} 
&\int_{\Sigma^{(i,j)} \cap \{r_3 < r < 2^{-\frac{1}{N}} \, r_j\}} \rho \, d\text{\rm vol}_g \\ 
&\geq \int_{\Sigma^{(i,j)} \cap \{r_4 < r < 2^{-\frac{1}{N}} \, r_j\}} (1-C_2 \, r^{-N}) \, r^{N-n} \, d\text{\rm vol}_{\bar{g}} 
\end{align*} 
for each $j$. This implies 
\begin{align*} 
&\int_{\Sigma^{(i,j)} \cap \{r_3 < r < 2^{-\frac{1}{N}} \, r_j\}} \rho \, d\text{\rm vol}_g \\ 
&\geq 2 \cdot (2\pi)^{n-2} \, \Big ( \prod_{k=0}^{n-3} b_k \Big ) \, \int_{r_4}^{2^{-\frac{1}{N}} \, r_j} (1-C_2 \, r^{-N}) \, r^{N-3} \, dr, 
\end{align*} 
for each $j$. This inequality contradicts Proposition \ref{upper.bound.for.area.1} if $j$ is sufficiently large. This completes the proof of Corollary \ref{graphical}. \\

\begin{definition} 
Given two integers $i,j$, let us consider the map $\pi^{(i,j)}: \Sigma^{(i,j)} \cap \{r_3 < r < 2^{-\frac{1}{N}} \, r_j\} \to (r_3,2^{-\frac{1}{N}} \, r_j) \times T^{n-2}$ defined in Proposition \ref{graphical}. We denote by $g_{\text{\rm hyp}}$ the pull-back of the hyperbolic metric $r^{-2} \, dr \otimes dr + \sum_{k=0}^{n-3} b_k^2 \, r^2 \, d\theta_k \otimes d\theta_k$ on $(r_3,2^{-\frac{1}{N}} \, r_j) \times T^{n-2}$ under the map $\pi^{(i,j)}$. Note that $g_{\text{\rm hyp}}$ is a hyperbolic metric on $\Sigma^{(i,j)} \cap \{r_3 < r < 2^{-\frac{1}{N}} \, r_j\}$. The metric $g_{\text{\rm hyp}}$ is obtained by restricting the $(0,2)$-tensor $\bar{g} - b_{n-2}^2 \, r^2 \, d\theta_{n-2} \otimes d\theta_{n-2}$ on $(r_3,2^{-\frac{1}{N}} \, r_j) \times T^{n-1}$ to $\Sigma^{(i,j)} \cap \{r_3 < r < 2^{-\frac{1}{N}} \, r_j\}$.
\end{definition}

\begin{proposition}
\label{higher.derivatives.of.height.function.with.respect.to.g_hyp}
For every positive integer $m$, we have $|D_{\text{\rm hyp}}^{\Sigma^{(i,j)},m} F^{(i,j)}|_{g_{\text{\rm hyp}}} \leq C(m) \, r^{-N}$ at each point in $\Sigma^{(i,j)} \cap \{r_3 < r < 2^{-\frac{1}{N}} \, r_j\}$. Here, $D_{\text{\rm hyp}}^{\Sigma^{(i,j)},m}$ denotes the covariant derivative of order $m$ with respect to the metric $g_{\text{\rm hyp}}$. Note that $C(m)$ may depend on $m$, but is independent of $i$ and $j$.
\end{proposition}

\textbf{Proof.} 
It follows from Proposition \ref{curvature.bound} that $|D^{\Sigma^{(i,j)},m} g_{\text{\rm hyp}}| \leq C(m)$ at each point in $\Sigma^{(i,j)} \cap \{r_3 < r < 2^{-\frac{1}{N}} \, r_j\}$. Moreover, Proposition \ref{higher.derivatives.of.height.function} implies that $|D^{\Sigma^{(i,j)},m} F^{(i,j)}| \leq C(m) \, r^{-N}$ at each point in $\Sigma^{(i,j)} \cap \{r_3 < r < 2^{-\frac{1}{N}} \, r_j\}$. Putting these facts together, the assertion follows. \\

\begin{proposition}
\label{bound.for.intrinsic.distance}
Let $d_{(\Sigma^{(i,j)},g^{(i)})}$ denote the intrinsic distance on $(\Sigma^{(i,j)},g^{(i)})$. Then 
\[\sup_{q \in \Sigma^{(i,j)} \setminus \{r \geq 2r_3\}} d_{(\Sigma^{(i,j)},g^{(i)})} \big ( q,\Sigma^{(i,j)} \cap \{r=2r_3\} \big ) \leq C.\] 
Here, $C$ is a large constant that may depend on $r_3$, but is independent of $i$ and $j$.
\end{proposition}

\textbf{Proof.} 
For each point $q \in \Sigma^{(i,j)} \setminus \{r \geq 2r_3\}$, we denote by $B_{(\Sigma^{(i,j)},g^{(i)})}(q,1)$ the intrinsic ball around $q$ of radius $1$ in $(\Sigma^{(i,j)},g^{(i)})$. Using the curvature bound in Proposition \ref{curvature.bound}, it is easy to see that the $(g^{(i)},\rho)$-area of $B_{(\Sigma^{(i,j)},g^{(i)})}(q,1)$ is bounded from below by a positive constant that may depend on $r_3$, but is independent of $i$ and $j$. Consequently, the $(g^{(i)},\rho)$-area of $\Sigma^{(i,j)} \setminus \{r \geq 2r_3\}$ is bounded from below by 
\[c \, \bigg ( \sup_{q \in \Sigma^{(i,j)} \setminus \{r \geq 2r_3\}} d_{(\Sigma^{(i,j)},g^{(i)})} \big ( q,\Sigma^{(i,j)} \cap \{r=2r_3\} \big ) - 4 \bigg ),\] 
where $c$ is a positive constant that may depend on $r_3$, but is independent of $i$ and $j$. On the other hand, Corollary \ref{upper.bound.for.area.2} gives an upper bound for the $(g^{(i)},\rho)$-area of $\Sigma^{(i,j)} \setminus \{r \geq 2r_3\}$. This completes the proof of Proposition \ref{bound.for.intrinsic.distance}. \\

The hypersurfaces $\Sigma^{(i,j)}$ satisfy the following stability inequality.

\begin{proposition}
\label{stability.inequality.1}
Let $a$ be a real number and let $V$ be a smooth vector field on $M$ with the property that $V = a \, \frac{\partial}{\partial \theta_{n-2}}$ in a neighborhood of the set $\{r=r_j\}$. Then 
\begin{align*} 
&\frac{1}{2} \int_{\Sigma^{(i,j)}} \rho \, \sum_{k=1}^{n-1} (\mathscr{L}_V \mathscr{L}_V g^{(i)})(e_k,e_k) \, d\text{\rm vol}_{g^{(i)}} + \int_{\Sigma^{(i,j)}} V(V(\rho)) \, d\text{\rm vol}_{g^{(i)}} \\ 
&- \frac{1}{2} \int_{\Sigma^{(i,j)}} \rho \, \sum_{k,l=1}^{n-1} (\mathscr{L}_V g^{(i)})(e_k,e_l) \, (\mathscr{L}_V g^{(i)})(e_k,e_l) \, d\text{\rm vol}_{g^{(i)}} \\ 
&+ \frac{1}{4} \int_{\Sigma^{(i,j)}} \rho \, \sum_{k,l=1}^{n-1} (\mathscr{L}_V g^{(i)})(e_k,e_k) \, (\mathscr{L}_V g^{(i)})(e_l,e_l) \, d\text{\rm vol}_{g^{(i)}} \\ 
&+ \int_{\Sigma^{(i,j)}} V(\rho) \, \sum_{k=1}^{n-1} (\mathscr{L}_V g^{(i)})(e_k,e_k) \, d\text{\rm vol}_{g^{(i)}} \\ 
&\geq -a^2 \int_{T^{n-2} \times \{t_{i,j}\}} \frac{\partial^2 u}{\partial \theta_{n-2}^2} \, d\text{\rm vol}_\gamma. 
\end{align*} 
Here, $\{e_1,\hdots,e_{n-1}\}$ denotes a local orthonormal frame on $(\Sigma^{(i,j)},g^{(i)})$.
\end{proposition}

\textbf{Proof.} 
Let $\varphi_s: M \to M$ denote the flow generated by the vector field $V$. We compute 
\[\frac{\partial}{\partial s} \varphi_s^*(g^{(i)}) \Big |_{s=0} = \mathscr{L}_V g^{(i)}, \quad \frac{\partial^2}{\partial s^2} \varphi_s^*(g^{(i)}) \Big |_{s=0} = \mathscr{L}_V \mathscr{L}_V g^{(i)}\] 
and 
\[\frac{\partial}{\partial s} (\rho \circ \varphi_s) \Big |_{s=0} = V(\rho), \quad \frac{\partial^2}{\partial s^2} (\rho \circ \varphi_s) \Big |_{s=0} = V(V(\rho)).\] 
This implies 
\begin{align*} 
&\frac{d^2}{ds^2} \bigg ( \int_{\varphi_s(\Sigma^{(i,j)})} \rho \, d\text{\rm vol}_{g^{(i)}} \bigg ) \bigg |_{s=0} \\ 
&= \frac{d^2}{ds^2} \bigg ( \int_{\Sigma^{(i,j)}} (\rho \circ \varphi_s) \, d\text{\rm vol}_{\varphi_s^*(g^{(i)})} \bigg ) \bigg |_{s=0} \\ 
&= \frac{1}{2} \int_{\Sigma^{(i,j)}} \rho \, \sum_{k=1}^{n-1} (\mathscr{L}_V \mathscr{L}_V g^{(i)})(e_k,e_k) \, d\text{\rm vol}_{g^{(i)}} + \int_{\Sigma^{(i,j)}} V(V(\rho)) \, d\text{\rm vol}_{g^{(i)}} \\ 
&- \frac{1}{2} \int_{\Sigma^{(i,j)}} \rho \, \sum_{k,l=1}^{n-1} (\mathscr{L}_V g^{(i)})(e_k,e_l) \, (\mathscr{L}_V g^{(i)})(e_k,e_l) \, d\text{\rm vol}_{g^{(i)}} \\ 
&+ \frac{1}{4} \int_{\Sigma^{(i,j)}} \rho \, \sum_{k,l=1}^{n-1} (\mathscr{L}_V g^{(i)})(e_k,e_k) \, (\mathscr{L}_V g^{(i)})(e_l,e_l) \, d\text{\rm vol}_{g^{(i)}} \\ 
&+ \int_{\Sigma^{(i,j)}} V(\rho) \, \sum_{k=1}^{n-1} (\mathscr{L}_V g^{(i)})(e_k,e_k) \, d\text{\rm vol}_{g^{(i)}}. 
\end{align*} 
Here, $\{e_1,\hdots,e_{n-1}\}$ denotes a local orthonormal frame on $(\Sigma^{(i,j)},g^{(i)})$. On the other hand, since the function (\ref{modified.area}) attains its minimum at $t_{i,j}$, we obtain 
\[\frac{d^2}{ds^2} \bigg ( \int_{\varphi_s(\Sigma^{(i,j)})} \rho \, d\text{\rm vol}_{g^{(i)}} \bigg ) \bigg |_{s=0} \geq -a^2 \int_{T^{n-2} \times \{t_{i,j}\}} \frac{\partial^2 u}{\partial \theta_{n-2}^2} \, d\text{\rm vol}_\gamma.\] 
Putting these facts together, the assertion follows. This completes the proof of Proposition \ref{stability.inequality.1}. \\

\begin{corollary}
\label{stability.inequality.2}
Let $\bar{r} \in (r_1,r_j)$. Moreover, let $a$ be a real number and let $V$ be a smooth vector field on $M$ with the property that $V = a \, \frac{\partial}{\partial \theta_{n-2}}$ in a neighborhood of the set $\{r=\bar{r}\}$. Then 
\begin{align*} 
&\frac{1}{2} \int_{\Sigma^{(i,j)} \setminus \{r>\bar{r}\}} \rho \, \sum_{k=1}^{n-1} (\mathscr{L}_V \mathscr{L}_V g^{(i)})(e_k,e_k) \, d\text{\rm vol}_{g^{(i)}} + \int_{\Sigma^{(i,j)} \setminus \{r>\bar{r}\}} V(V(\rho)) \, d\text{\rm vol}_{g^{(i)}} \\ 
&- \frac{1}{2} \int_{\Sigma^{(i,j)} \setminus \{r>\bar{r}\}} \rho \, \sum_{k,l=1}^{n-1} (\mathscr{L}_V g^{(i)})(e_k,e_l) \, (\mathscr{L}_V g^{(i)})(e_k,e_l) \, d\text{\rm vol}_{g^{(i)}} \\ 
&+ \frac{1}{4} \int_{\Sigma^{(i,j)} \setminus \{r>\bar{r}\}} \rho \, \sum_{k,l=1}^{n-1} (\mathscr{L}_V g^{(i)})(e_k,e_k) \, (\mathscr{L}_V g^{(i)})(e_l,e_l) \, d\text{\rm vol}_{g^{(i)}} \\ 
&+ \int_{\Sigma^{(i,j)} \setminus \{r>\bar{r}\}} V(\rho) \, \sum_{k=1}^{n-1} (\mathscr{L}_V g^{(i)})(e_k,e_k) \, d\text{\rm vol}_{g^{(i)}} \\ 
&\geq -a^2 \int_{T^{n-2} \times \{t_{i,j}\}} \frac{\partial^2 u}{\partial \theta_{n-2}^2} \, d\text{\rm vol}_\gamma - C \, a^2 \, \bar{r}^{-\delta}. 
\end{align*} 
Here, $\{e_1,\hdots,e_{n-1}\}$ denotes a local orthonormal frame on $(\Sigma^{(i,j)},g^{(i)})$, and $C$ is a positive constant which is independent of $i$, $j$, and $\bar{r}$.
\end{corollary}

\textbf{Proof.} 
We may assume that $V = a \, \frac{\partial}{\partial \theta_{n-2}}$ in the region $\{r>\bar{r}\}$. Using Lemma \ref{Lie.derivative.of.metric.along.angular.vector.field}, we obtain 
\[|\mathscr{L}_V g| \leq C \, |a| \, r^{1-N-\delta}, \quad |\mathscr{L}_V \mathscr{L}_V g| \leq C \, |a|^2 \, r^{2-N-\delta}\] 
in the region $\{r>\bar{r}\}$. Moreover, Lemma \ref{derivative.of.rho.along.angular.vector.field} gives 
\[|V(\rho)| \leq C \, |a| \, \rho \, r^{1-N-\delta}, \quad |V(V(\rho))| \leq C \, |a|^2 \, \rho \, r^{2-N-\delta}\] 
in the region $\{r>\bar{r}\}$. Consequently, 
\begin{align*} 
&\frac{1}{2} \int_{\Sigma^{(i,j)} \cap \{\bar{r} < r \leq r_j\}} \rho \, \sum_{k=1}^{n-1} (\mathscr{L}_V \mathscr{L}_V g)(e_k,e_k) \, d\text{\rm vol}_g + \int_{\Sigma^{(i,j)} \cap \{\bar{r} < r \leq r_j\}} V(V(\rho)) \, d\text{\rm vol}_g \\ 
&- \frac{1}{2} \int_{\Sigma^{(i,j)} \cap \{\bar{r} < r \leq r_j\}} \rho \, \sum_{k,l=1}^{n-1} (\mathscr{L}_V g)(e_k,e_l) \, (\mathscr{L}_V g)(e_k,e_l) \, d\text{\rm vol}_g \\ 
&+ \frac{1}{4} \int_{\Sigma^{(i,j)} \cap \{\bar{r} < r \leq r_j\}} \rho \, \sum_{k,l=1}^{n-1} (\mathscr{L}_V g)(e_k,e_k) \, (\mathscr{L}_V g)(e_l,e_l) \, d\text{\rm vol}_g \\ 
&+ \int_{\Sigma^{(i,j)} \cap \{\bar{r} < r \leq r_j\}} V(\rho) \, \sum_{k=1}^{n-1} (\mathscr{L}_V g)(e_k,e_k) \, d\text{\rm vol}_g \\ 
&\leq C \, |a|^2 \int_{\Sigma^{(i,j)} \cap \{\bar{r} < r \leq r_j\}} \rho \, r^{2-N-\delta}. 
\end{align*} 
Finally, Corollary \ref{upper.bound.for.area.2} implies 
\[\int_{\Sigma^{(i,j)} \cap \{\bar{r} < r \leq r_j\}} \rho \, r^{2-N-\delta} \leq C \, \bar{r}^{-\delta}.\] 
The assertion follows by combining these estimates with Proposition \ref{stability.inequality.1}. This completes the proof of Corollary \ref{stability.inequality.2}. \\

\begin{proposition}
\label{Sigma.passes.through.a.small.ball.around.p_*}
Suppose that the condition $(\star_{N,n-1})$ is satisfied. Then, for every integer $i$, there exists an integer $j \geq i$ with the property that $\Sigma^{(i,j)} \cap B_{(M,g)}(p_*,2\varepsilon_i) \neq \emptyset$.
\end{proposition}

\textbf{Proof.} 
We argue by contradiction. Suppose that there exists an integer $i$ with the property that $\Sigma^{(i,j)} \cap B_{(M,g)}(p_*,2\varepsilon_i) = \emptyset$ for all $j \geq i$. In the following, we fix an integer $i$ with this property. After passing to a subsequence, we may assume that the limit $\lim_{j \to \infty} t_{i,j} = t_*$ exists. The hypersurfaces $\Sigma^{(i,j)}$ satisfy local area bounds (see Corollary \ref{upper.bound.for.area.2}) and local curvature bounds (see Proposition \ref{curvature.bound}). Therefore, we may take a subsequential limit of the hypersurfaces $\Sigma^{(i,j)}$ as $j \to \infty$. In the limit, we obtain a properly embedded, connected, orientable hypersurface $\Sigma$ which is $t_*$-tame. The fact that $\Sigma$ is $t_*$-tame follows from Proposition \ref{higher.derivatives.of.height.function.with.respect.to.g_hyp}; the fact that $\Sigma$ is connected is a consequence of Proposition \ref{bound.for.intrinsic.distance}. Since $\Sigma^{(i,j)} \cap B_{(M,g)}(p_*,2\varepsilon_i) = \emptyset$ for all $j \geq i$, it follows that $\Sigma \cap B_{(M,g)}(p_*,\varepsilon_i) = \emptyset$. Since $\Gamma_{t_{i,j}}^{(j)}$ does not bound a hypersurface in $[4r_{\text{\rm fol}},\infty) \times T^{n-1}$, it follows that the set $\Sigma^{(i,j)} \cap \{r=4r_{\text{\rm fol}}\}$ is non-empty if $j$ is sufficiently large. Passing to the limit as $j \to \infty$, we conclude that the set $\Sigma \cap \{r=4r_{\text{\rm fol}}\}$ is non-empty. 

Clearly, $\Sigma$ is $(g^{(i)},\rho)$-stationary. We claim that $\Sigma$ is $(g^{(i)},\rho,u)$-stable in the sense of Definition \ref{definition.stability}. To see this, we fix a real number $a$ and a smooth vector field $V$ on $M$ such that $V = a \, \frac{\partial}{\partial \theta_{n-2}}$ outside a compact set. We apply Corollary \ref{stability.inequality.2} to the vector field $V$. In the first step, we pass to the limit as $j \to \infty$, keeping $i$ and $\bar{r}$ fixed. In the second step, we send $\bar{r} \to \infty$, keeping $i$ fixed. This shows that $\Sigma$ is $(g^{(i)},\rho,u)$-stable in the sense of Definition \ref{definition.stability}. 

As above, we define $U = M \setminus \{r \geq 8r_{\text{\rm fol}}\}$. If $n=N$, then Proposition \ref{perturbed.metrics} implies that $R_{g^{(i)}} + N(N-1) > 0$ at each point in $U \setminus B_{(M,g)}(p_*,\varepsilon_i)$. If $n<N$, then Proposition \ref{perturbed.metrics} implies that 
\[-2 \, \Delta_{g^{(i)}} \log \rho - \frac{N-n+1}{N-n} \, |d\log \rho|_{g^{(i)}}^2 + R_{g^{(i)}} + N(N-1) > 0\] 
at each point in $U \setminus B_{(M,g)}(p_*,\varepsilon_i)$. 

Let $\check{g}^{(i)}$ denote the restriction of the metric $g^{(i)}$ to $\Sigma$. The results in Section \ref{properties.of.stable.hypersurfaces} imply that we can find a positive smooth function $\check{v}: \Sigma \to \mathbb{R}$ with the property that $(\Sigma,\check{g}^{(i)},\check{\rho})$ is an $(N,n-1)$-dataset, where $\check{\rho}$ is defined by $\check{\rho} = b_{n-2}^{-1} \, \check{v} \, \rho$. Since $g^{(i)} = g$ on $M \setminus U$, Proposition \ref{scalar.curvature.of.hypersurface} implies that 
\[-2 \, \Delta_{\check{g}^{(i)}} \log \check{\rho} - \frac{N-n+2}{N-n+1} \, |d\log \check{\rho}|_{\check{g}^{(i)}}^2 + R_{\check{g}^{(i)}} + N(N-1) \geq 0\] 
at each point in $\Sigma \setminus U$. Moreover, since $\Sigma \cap U \subset U \setminus B_{(M,g)}(p_*,\varepsilon_i)$, it follows from the proof of Proposition \ref{scalar.curvature.of.hypersurface} that 
\[-2 \, \Delta_{\check{g}^{(i)}} \log \check{\rho} - \frac{N-n+2}{N-n+1} \, |d\log \check{\rho}|_{\check{g}^{(i)}}^2 + R_{\check{g}^{(i)}} + N(N-1) > 0\] 
at each point in $\Sigma \cap U$. Using condition $(\star_{N,n-1})$, we conclude that $(\Sigma,\check{g}^{(i)},\check{\rho})$ is a model $(N,n-1)$-dataset. Proposition \ref{scalar.curvature.of.HM.metrics} now implies that 
\[-2 \, \Delta_{\check{g}^{(i)}} \log \check{\rho} - \frac{N-n+2}{N-n+1} \, |d\log \check{\rho}|_{\check{g}^{(i)}}^2 + R_{\check{g}^{(i)}} + N(N-1) = 0\] 
at each point on $\Sigma$. Since the set $\Sigma \cap U$ is non-empty, we arrive at a contradiction. This completes the proof of Proposition \ref{Sigma.passes.through.a.small.ball.around.p_*}. \\

\begin{proposition}
\label{existence.of.totally.geodesic.hypersurface.passing.through.a.given.point}
Suppose that the condition $(\star_{N,n-1})$ is satisfied. Let us fix an arbitrary point $p_* \in M \setminus \{r \geq 2r_{\text{\rm fol}}\}$. Then we can find a properly embedded, connected, orientable hypersurface $\Sigma$ passing through $p_*$ and a positive smooth function $\check{v}: \Sigma \to \mathbb{R}$ with the following properties: 
\begin{itemize} 
\item The hypersurface $\Sigma$ is tame. 
\item The hypersurface $\Sigma$ is totally geodesic and the normal derivative of $\rho$ vanishes along $\Sigma$.
\item The function $\check{v}$ satisfies $\mathbb{L}_\Sigma \check{v} = 0$, where $\mathbb{L}_\Sigma$ denotes the weighted Jacobi operator of $\Sigma$ (see Definition \ref{definition.weighted.Jacobi.operator}). Moreover, $|\check{v} - b_{n-2} \, r| \leq O(r^{1-N})$.
\item If $n<N$, then the function $\check{v}^{-(N-n)} \, \rho$ is constant along $\Sigma$.
\item Let $\check{g}$ denote the restriction of $g$ to $\Sigma$ and let $\check{\rho} = b_{n-2}^{-1} \, \check{v} \, \rho$. Then $(\Sigma,\check{g},\check{\rho})$ is a model $(N,n-1)$-dataset.
\item For each $\bar{r} > 8r_{\text{\rm fol}}$, the $(g,\rho)$-area of $\Sigma \setminus \{r > \bar{r}\}$ is bounded by $C \, \bar{r}^{N-2}$ for some uniform constant $C$. 
\end{itemize}
\end{proposition}

\textbf{Proof.} 
Proposition \ref{Sigma.passes.through.a.small.ball.around.p_*} implies that, for each integer $i$, we can find an integer $j_0(i) \geq i$ with the property that $\Sigma^{(i,j_0(i))} \cap B_{(M,g)}(p_*,2\varepsilon_i) \neq \emptyset$. After passing to a subsequence, we may assume that the limit $\lim_{i \to \infty} t_{i,j_0(i)} = t_*$ exists. The hypersurfaces $\Sigma^{(i,j_0(i))}$ satisfy local area bounds (see Corollary \ref{upper.bound.for.area.2}) and local curvature bounds (see Proposition \ref{curvature.bound}). Therefore, we may take a subsequential limit of the hypersurfaces $\Sigma^{(i,j_0(i))}$ as $i \to \infty$. In the limit, we obtain a properly embedded, connected, orientable hypersurface $\Sigma$ which is $t_*$-tame. The fact that $\Sigma$ is $t_*$-tame follows from Proposition \ref{higher.derivatives.of.height.function.with.respect.to.g_hyp}; the fact that $\Sigma$ is connected is a consequence of Proposition \ref{bound.for.intrinsic.distance}. Since $\Sigma^{(i,j_0(i))} \cap B_{(M,g)}(p_*,2\varepsilon_i) \neq \emptyset$ for each $i$, we conclude that $\Sigma$ passes through the point $p_*$. 

Clearly, $\Sigma$ is $(g,\rho)$-stationary. We claim that $\Sigma$ is $(g,\rho,u)$-stable in the sense of Definition \ref{definition.stability}. To see this, we fix a real number $a$ and a smooth vector field $V$ on $M$ such that $V = a \, \frac{\partial}{\partial \theta_{n-2}}$ outside a compact set. We apply Corollary \ref{stability.inequality.2} to the vector field $V$ and the hypersurface $\Sigma^{(i,j_0(i))}$. In the first step, we pass to the limit as $i \to \infty$, keeping $\bar{r}$ fixed. In the second step, we send $\bar{r} \to \infty$. This shows that $\Sigma$ is $(g,\rho,u)$-stable in the sense of Definition \ref{definition.stability}. 

If $n=N$, then $R + N(N-1) \geq 0$ at each point on $\Sigma$. If $n<N$, then 
\[-2 \, \Delta \log \rho - \frac{N-n+1}{N-n} \, |d\log \rho|^2 + R + N(N-1) \geq 0\] 
at each point on $\Sigma$. 

Let $\check{g}$ denote the restriction of the metric $g$ to $\Sigma$. The results in Section \ref{properties.of.stable.hypersurfaces} imply that we can find a positive smooth function $\check{v}: \Sigma \to \mathbb{R}$ with the property that $(\Sigma,\check{g},\check{\rho})$ is an $(N,n-1)$-dataset, where $\check{\rho}$ is defined by $\check{\rho} = b_{n-2}^{-1} \, \check{v} \, \rho$. Using Proposition \ref{scalar.curvature.of.hypersurface}, we obtain 
\[-2 \, \Delta_{\check{g}} \log \check{\rho} - \frac{N-n+2}{N-n+1} \, |d\log \check{\rho}|_{\check{g}}^2 + R_{\check{g}} + N(N-1) \geq 0\] 
at each point on $\Sigma$. Using condition $(\star_{N,n-1})$, we conclude that $(\Sigma,\check{g},\check{\rho})$ is a model $(N,n-1)$-dataset. Proposition \ref{scalar.curvature.of.HM.metrics} now implies that 
\[-2 \, \Delta_{\check{g}} \log \check{\rho} - \frac{N-n+2}{N-n+1} \, |d\log \check{\rho}|_{\check{g}}^2 + R_{\check{g}} + N(N-1) = 0\] 
at each point on $\Sigma$. Thus, equality holds in Proposition \ref{scalar.curvature.of.hypersurface}. From this, we deduce that $\Sigma$ is totally geodesic and $\mathbb{L}_\Sigma \check{v} = 0$ at each point on $\Sigma$. Moreover, if $n<N$, then the function $\check{v}^{-(N-n)} \, \rho$ is constant along $\Sigma$. Since $\Sigma$ is $(g,\rho)$-stationary and totally geodesic, it follows that the normal derivative of $\rho$ vanishes at each point on $\Sigma$. Finally, Corollary \ref{upper.bound.for.area.2} gives an upper bound for the $(g,\rho)$-area of $\Sigma \setminus \{r > \bar{r}\}$ for each $\bar{r} > 8r_{\text{\rm fol}}$. This completes the proof of Proposition \ref{existence.of.totally.geodesic.hypersurface.passing.through.a.given.point}. \\

\section{Proof of Theorem \ref{main.theorem}}

\label{proof.of.main.theorem}

In this section, we establish the following theorem.

\begin{theorem}
\label{property.star.holds}
Let us fix an integer $N$ with $3 \leq N \leq 7$. Then property $(\star_{N,n})$ is satisfied for each $2 \leq n \leq N$. 
\end{theorem}

Theorem \ref{main.theorem} follows by putting $n=N$ in Theorem \ref{property.star.holds}.

To prove Theorem \ref{property.star.holds}, we fix an integer $N$ with $3 \leq N \leq 7$. We argue by induction on $n$. Theorem \ref{2D.case} implies that $(\star_{N,2})$ holds. Suppose next that $3 \leq n \leq N$ and $(\star_{N,n-1})$ holds. Our goal is to show that $(\star_{N,n})$ holds. To that end, suppose that $(M,g,\rho)$ is an $(N,n)$-dataset. If $n=N$, we assume that $\rho=1$ and $R + N(N-1) \geq 0$ at each point in $M$. If $n<N$, we assume that  
\[-2 \, \Delta \log \rho - \frac{N-n+1}{N-n} \, |\nabla \log \rho|^2 + R + N(N-1) \geq 0\] 
at each point in $M$. We will show that $(M,g,\rho)$ is a model $(N,n)$-dataset. 

\begin{proposition}
\label{existence.of.totally.geodesic.hypersurface.with.given.asymptotics}
For each $\bar{t} \in S^1$, we can find a properly embedded, connected, orientable hypersurface $\Sigma$ and a positive smooth function $\check{v}: \Sigma \to \mathbb{R}$ with the following properties: 
\begin{itemize} 
\item The hypersurface $\Sigma$ is $\bar{t}$-tame. 
\item The hypersurface $\Sigma$ is totally geodesic and the normal derivative of $\rho$ vanishes along $\Sigma$.
\item The function $\check{v}$ satisfies $\mathbb{L}_\Sigma \check{v} = 0$, where $\mathbb{L}_\Sigma$ denotes the weighted Jacobi operator of $\Sigma$ (see Definition \ref{definition.weighted.Jacobi.operator}). Moreover, $|\check{v} - b_{n-2} \, r| \leq O(r^{1-N})$.
\item If $n<N$, then the function $\check{v}^{-(N-n)} \, \rho$ is constant along $\Sigma$.
\item Let $\check{g}$ denote the restriction of $g$ to $\Sigma$ and let $\check{\rho} = b_{n-2}^{-1} \, \check{v} \, \rho$. Then $(\Sigma,\check{g},\check{\rho})$ is a model $(N,n-1)$-dataset.
\item For each $\bar{r} > 8r_{\text{\rm fol}}$, the $(g,\rho)$-area of $\Sigma \setminus \{r > \bar{r}\}$ is bounded by $C \, \bar{r}^{N-2}$ for some uniform constant $C$. 
\item We have $\Sigma \cap \{r > r_{\text{\rm fol}}\} = \mathcal{Z}_{\bar{t}}$.
\end{itemize}
\end{proposition}

\textbf{Proof.} 
Let us fix an element $\bar{t} \in S^1$, and let $p_*$ be an arbitrary point in $\mathcal{Z}_{\bar{t}} \cap \{r_{\text{\rm fol}} < r < 2r_{\text{\rm fol}}\}$. By the inductive hypothesis, property $(\star_{N,n-1})$ is satisfied. By Proposition \ref{existence.of.totally.geodesic.hypersurface.passing.through.a.given.point}, we can find a properly embedded, connected, orientable hypersurface $\Sigma$ passing through $p_*$ and a positive smooth function $\check{v}: \Sigma \to \mathbb{R}$ with the following properties: 
\begin{itemize} 
\item The hypersurface $\Sigma$ is $t_*$-tame for some element $t_* \in S^1$. 
\item The hypersurface $\Sigma$ is totally geodesic and the normal derivative of $\rho$ vanishes along $\Sigma$.
\item The function $\check{v}$ satisfies $\mathbb{L}_\Sigma \check{v} = 0$, where $\mathbb{L}_\Sigma$ denotes the weighted Jacobi operator of $\Sigma$. Moreover, $|\check{v} - b_{n-2} \, r| \leq O(r^{1-N})$.
\item If $n<N$, then the function $\check{v}^{-(N-n)} \, \rho$ is constant along $\Sigma$.
\item Let $\check{g}$ denote the restriction of $g$ to $\Sigma$ and let $\check{\rho} = b_{n-2}^{-1} \, \check{v} \, \rho$. Then $(\Sigma,\check{g},\check{\rho})$ is a model $(N,n-1)$-dataset.
\item For each $\bar{r} > 8r_{\text{\rm fol}}$, the $(g,\rho)$-area of $\Sigma \setminus \{r > \bar{r}\}$ is bounded by $C \, \bar{r}^{N-2}$ for some uniform constant $C$. 
\end{itemize}
Proposition \ref{tame.and.totally.geodesic} implies that $\Sigma \cap \{r > r_{\text{\rm fol}}\} = \mathcal{Z}_{t_*}$. Since $p_* \in \Sigma \cap \{r > r_{\text{\rm fol}}\}$, it follows that $p_* \in \mathcal{Z}_{t_*}$. On the other hand, $p_* \in \mathcal{Z}_{\bar{t}}$ by assumption. Consequently, $t_* = \bar{t}$. This completes the proof of Proposition \ref{existence.of.totally.geodesic.hypersurface.with.given.asymptotics}. \\

\begin{definition}
For each $\bar{t} \in S^1$, we denote by $\Sigma_{\bar{t}}$ the unique hypersurface satisfying the conclusion of Proposition \ref{existence.of.totally.geodesic.hypersurface.with.given.asymptotics}. We denote by $\check{g}_{\bar{t}}$ the induced metric on $\Sigma_{\bar{t}}$. Moreover, we denote by $\check{v}_{\bar{t}}: \Sigma_{\bar{t}} \to \mathbb{R}$ and $\check{\rho}_{\bar{t}}: \Sigma_{\bar{t}} \to \mathbb{R}$ the functions constructed in Proposition \ref{existence.of.totally.geodesic.hypersurface.with.given.asymptotics}. Note that $\check{\rho}_{\bar{t}} = b_{n-2}^{-1} \, \check{v}_{\bar{t}} \, \rho$ at each point on $\Sigma_{\bar{t}}$.
\end{definition}

\begin{proposition}
\label{union.of.Sigma_t}
We have $\bigcup_{t \in S^1} \Sigma_t = M$. 
\end{proposition}

\textbf{Proof.} 
Let us fix an arbitrary point $p_* \in M \setminus \{r \geq 2r_{\text{\rm fol}}\}$. By Proposition \ref{existence.of.totally.geodesic.hypersurface.passing.through.a.given.point}, we can find a properly embedded, connected, orientable hypersurface $\Sigma$ passing through $p_*$ with the property that $\Sigma$ is tame and totally geodesic. By Proposition \ref{tame.and.totally.geodesic}, we can find an element $t_* \in S^1$ such that $\Sigma \cap \{r > r_{\text{\rm fol}}\} = \mathcal{Z}_{t_*}$. Proposition \ref{uniqueness} now implies that $\Sigma = \Sigma_{t_*}$. Since $p_* \in \Sigma$, we conclude that $p_* \in \Sigma_{t_*}$. Thus, 
\[M \setminus \{r \geq 2r_{\text{\rm fol}}\} \subset \bigcup_{t \in S^1} \Sigma_t.\] 
On the other hand, using Proposition \ref{existence.of.totally.geodesic.hypersurface.with.given.asymptotics}, we obtain 
\[[2r_{\text{\rm fol}},\infty) \times T^{n-1} \subset \bigcup_{t \in S^1} \mathcal{Z}_t \subset \bigcup_{t \in S^1} \Sigma_t.\] 
This completes the proof of Proposition \ref{union.of.Sigma_t}. \\

\begin{proposition}
\label{Jacobi.field}
Let us fix an element $\bar{t} \in S^1$. Let us consider a sequence $t_j \in S^1$ such that $t_j \neq \bar{t}$ for each $j$ and $t_j \to \bar{t}$ as $j \to \infty$. After passing to a subsequence, we can find a sequence of positive real numbers $\delta_j \to 0$ and a sequence of smooth functions $w^{(j)}: \Sigma_{\bar{t}} \setminus \{r>\delta_j^{-1}\} \to \mathbb{R}$ with the following properties: 
\begin{itemize}
\item For every nonnegative integer $m$, we have 
\[\sup_{\Sigma_{\bar{t}} \setminus \{r>\delta_j^{-1}\}} |D^{\Sigma_{\bar{t}},m} w^{(j)}| \to 0\] 
as $j \to \infty$.
\item If $j$ is sufficiently large, then  
\[\exp_x(w^{(j)}(x) \, \nu_{\Sigma_{\bar{t}}}(x)) \in \Sigma_{t_j}\] 
for all points $x \in \Sigma_{\bar{t}} \setminus \{r>\delta_j^{-1}\}$. 
\item The rescaled functions $d_{S^1}(t_j,\bar{t})^{-1} \, w^{(j)}$ converge in $C_{\text{\rm loc}}^\infty$ to a smooth function $w: \Sigma_{\bar{t}} \to \mathbb{R}$.
\item The function $\frac{w}{\check{v}_{\bar{t}}}$ is equal to a non-zero constant.
\end{itemize} 
\end{proposition} 

\textbf{Proof.} 
It follows from Proposition \ref{existence.of.totally.geodesic.hypersurface.with.given.asymptotics} that the hypersurfaces $\Sigma_{t_j}$ satisfy local area bounds. After passing to a subsequence, the hypersurfaces $\Sigma_{t_j}$ converge, in the sense of measures, to an integer multiplicity rectifiable varifold. The support of this limiting varifold is a closed subset of $M$ which we denote by $\hat{\Sigma}$. Since the hypersurfaces $\Sigma_{t_j}$ are totally geodesic, it follows that $\hat{\Sigma}$ is a smooth (possibly disconnected) submanifold of $M$. After passing to a subsequence, the hypersurfaces $\Sigma_{t_j}$ converge to $\hat{\Sigma}$ smoothly on compact subsets of $M$. Note that the convergence may be multi-sheeted. 

Since $\Sigma_{t_j} \cap \{r > r_{\text{\rm fol}}\} = \mathcal{Z}_{t_j} \cap \{r > r_{\text{\rm fol}}\}$ for each $j$, it follows that $\hat{\Sigma} \cap \{r > 2r_{\text{\rm fol}}\} = \mathcal{Z}_{\bar{t}} \cap \{r > 2r_{\text{\rm fol}}\}$. Let $\tilde{\Sigma}$ denote the connected component of $\hat{\Sigma}$ that contains the set $\mathcal{Z}_{\bar{t}} \cap \{r > 2r_{\text{\rm fol}}\}$. Using Proposition \ref{uniqueness}, we conclude that $\tilde{\Sigma} = \Sigma_{\bar{t}}$. 

Note that the multiplicity of the limiting varifold is locally constant, and is equal to $1$ on $\tilde{\Sigma}$. Consequently, we can find a sequence of positive real numbers $\delta_j \to 0$ and a sequence of smooth functions $w^{(j)}: \Sigma_{\bar{t}} \setminus \{r>\delta_j^{-1}\} \to \mathbb{R}$ such that 
\[\sup_{\Sigma_{\bar{t}} \setminus \{r>\delta_j^{-1}\}} |D^{\Sigma_{\bar{t}},m} w^{(j)}| \to 0\] 
for every nonnegative integer $m$ and 
\[\exp_x(w^{(j)}(x) \, \nu_{\Sigma_{\bar{t}}}(x)) \in \Sigma_{t_j}\] 
for all points $x \in \Sigma_{\bar{t}} \setminus \{r>\delta_j^{-1}\}$. 

For each $t \in S^1$, the hypersurface $\Sigma_t$ is $(g,\rho)$-stationary. Let $\mathbb{L}_{\Sigma_{\bar{t}}}$ denote the weighted Jacobi operator of $\Sigma_{\bar{t}}$ (see Definition \ref{definition.weighted.Jacobi.operator}). The function $w^{(j)}$ satisfies an equation of the form $\tilde{\mathbb{L}}^{(j)} w^{(j)} = 0$. Here, $\tilde{\mathbb{L}}^{(j)}$ is a linear differential operator of second order on $\Sigma_{\bar{t}} \setminus \{r>\delta_j^{-1}\}$ with coefficients that may depend on $w^{(j)}$ and its first derivatives. Since $w^{(j)} \to 0$ in $C_{\text{\rm loc}}^\infty$, the coefficients of $\tilde{\mathbb{L}}^{(j)}$ converge to the corresponding coefficients of $\mathbb{L}_{\Sigma_{\bar{t}}}$ in $C_{\text{\rm loc}}^\infty$. 

We next observe that $\Sigma_t \cap \{r > r_{\text{\rm fol}}\} = \mathcal{Z}_t$ for each $t \in S^1$. From this, we deduce that 
\begin{equation} 
\label{bound.for.height.function}
\limsup_{j \to \infty} \sup_{\Sigma_{\bar{t}} \cap \{2r_{\text{\rm fol}} \leq r \leq \bar{r}\}} d_{S^1}(t_j,\bar{t})^{-1} \, |w^{(j)}| < \infty 
\end{equation} 
for each $\bar{r} \in (2r_{\text{\rm fol}},\infty)$. For each $j$, we define 
\begin{equation} 
\label{definition.of.alpha_j}
\alpha_j = \sup_{\Sigma_{\bar{t}} \setminus \{r>4r_{\text{\rm fol}}\}} |w^{(j)}|. 
\end{equation}
Clearly, $\alpha_j \to 0$ as $j \to \infty$. 

We claim that $\limsup_{j \to \infty} d_{S^1}(t_j,\bar{t})^{-1} \, \alpha_j < \infty$. To prove this, we argue by contradiction. Suppose that $\limsup_{j \to \infty} d_{S^1}(t_j,\bar{t})^{-1} \, \alpha_j = \infty$. After passing to a subsequence, we may assume that $\liminf_{j \to \infty} d_{S^1}(t_j,\bar{t})^{-1} \, \alpha_j = \infty$. Using (\ref{bound.for.height.function}) and (\ref{definition.of.alpha_j}), we obtain 
\[\limsup_{j \to \infty} \sup_{\Sigma_{\bar{t}} \setminus \{r>\bar{r}\}} \alpha_j^{-1} \, |w^{(j)}| < \infty\] 
for each $\bar{r} \in (2r_{\text{\rm fol}},\infty)$. After passing to a subsequence, the rescaled functions $\alpha_j^{-1} \, w^{(j)}$ converge in $C_{\text{\rm loc}}^\infty$ to a smooth function $\tilde{w}: \Sigma_{\bar{t}} \to \mathbb{R}$ satisfying $\mathbb{L}_{\Sigma_{\bar{t}}} \tilde{w} = 0$. Moreover, it follows from (\ref{bound.for.height.function}) that the function $\tilde{w}$ vanishes identically outside some compact set. Since $\Sigma_{\bar{t}}$ is connected, standard unique continuation theorems for elliptic PDE (see e.g. \cite{Aronszajn}) imply that $\tilde{w}$ vanishes identically. This leads a contradiction. Thus, $\limsup_{j \to \infty} d_{S^1}(t_j,\bar{t})^{-1} \, \alpha_j < \infty$. 

Using (\ref{bound.for.height.function}) and (\ref{definition.of.alpha_j}), we obtain 
\[\limsup_{j \to \infty} \sup_{\Sigma_{\bar{t}} \setminus \{r>\bar{r}\}} d_{S^1}(t_j,\bar{t})^{-1} \, |w^{(j)}| < \infty\] 
for each $\bar{r} \in (2r_{\text{\rm fol}},\infty)$. After passing to a subsequence, the rescaled functions $d_{S^1}(t_j,\bar{t})^{-1} \, w^{(j)}$ converge in $C_{\text{\rm loc}}^\infty$ to a smooth function $w: \Sigma_{\bar{t}} \to \mathbb{R}$ satisfying $\mathbb{L}_{\Sigma_{\bar{t}}} w = 0$. On the other hand, $\check{v}_{\bar{t}}$ is a positive smooth function on $\Sigma_{\bar{t}}$ satisfying $\mathbb{L}_{\Sigma_{\bar{t}}} \check{v}_{\bar{t}} = 0$. Putting these facts together, we conclude that 
\begin{equation} 
\label{quotient}
-\check{v}_{\bar{t}} \, \text{\rm div}_{\Sigma_{\bar{t}}} \Big ( \rho \, \nabla^{\Sigma_{\bar{t}}} \Big ( \frac{w}{\check{v}_{\bar{t}}} \Big ) \Big ) - 2\rho \, \Big \langle \nabla^{\Sigma_{\bar{t}}} \check{v}_{\bar{t}},\nabla^{\Sigma_{\bar{t}}} \Big ( \frac{w}{\check{v}_{\bar{t}}} \Big ) \Big \rangle = 0 
\end{equation} 
at each point on $\Sigma_{\bar{t}}$. 

Finally, since $\Sigma_t \cap \{r > r_{\text{\rm fol}}\} = \mathcal{Z}_t$ for each $t \in S^1$, it follows that $w = |d\Xi|_g^{-1}$ outside a compact set, where the function $\Xi$ is defined as in Section \ref{foliation}. Note that $r \, |d\Xi|_g$ converges to a non-zero constant at infinity. Consequently, the function $\frac{w}{r}$ converges to a non-zero constant at infinity. On the other hand, the function $\frac{\check{v}_{\bar{t}}}{r}$ also converges to a non-zero constant at infinity. Putting these facts together, it follows that the function $\frac{w}{\check{v}_{\bar{t}}}$ converges to a non-zero constant at infinity. Using (\ref{quotient}) and the maximum principle, we conclude that the function $\frac{w}{\check{v}_{\bar{t}}}$ is equal to a non-zero constant. This completes the proof of Proposition \ref{Jacobi.field}. \\

\begin{corollary}
\label{Hessian.of.check.v}
Let us fix an element $\bar{t} \in S^1$. If $\{e_1,\hdots,e_{n-1}\}$ is a local orthonormal frame on $\Sigma_{\bar{t}}$, then 
\[(D^{\Sigma_{\bar{t}},2} \check{v}_{\bar{t}})(e_i,e_k) + R(e_i,\nu_{\Sigma_{\bar{t}}},e_k,\nu_{\Sigma_{\bar{t}}}) \, \check{v}_{\bar{t}} = 0\] 
for all $i,k \in \{1,\hdots,n-1\}$.
\end{corollary} 

\textbf{Proof.} 
Let us consider a sequence $t_j \in S^1$ such that $t_j \neq \bar{t}$ for each $j$ and $t_j \to \bar{t}$ as $j \to \infty$. Let $w: \Sigma_{\bar{t}} \to \mathbb{R}$ be defined as in Proposition \ref{Jacobi.field}. The hypersurfaces $\Sigma_{t_j}$ are totally geodesic, and the hypersurface $\Sigma_{\bar{t}}$ is totally geodesic as well. Consequently, the function $w$ satisfies 
\[(D^{\Sigma_{\bar{t}},2} w)(e_i,e_k) + R(e_i,\nu_{\Sigma_{\bar{t}}},e_k,\nu_{\Sigma_{\bar{t}}}) \, w = 0\] 
for all $i,k \in \{1,\hdots,n-1\}$. Finally, $w$ is a non-zero multiple of $\check{v}_{\bar{t}}$ by Proposition \ref{Jacobi.field}. This completes the proof of Corollary \ref{Hessian.of.check.v}. \\

\begin{proposition}
\label{curvature.tensor.of.M}
Let $p$ be an arbitrary point in $M$. There exists a symmetric bilinear form $T: T_p M \times T_p M \to \mathbb{R}$ and a real number $\Upsilon \in [1,\infty)$ with the following properties: 
\begin{itemize}
\item The eigenvalues of $T$ are $1$ and $0$, and the corresponding multiplicities are $2$ and $n-2$, respectively. 
\item The Riemann curvature tensor of $(M,g)$ at the point $p$ is given by 
\[-\frac{1}{2} \, (1-\Upsilon^{-N}) \, g \owedge g - \frac{N}{2} \, \Upsilon^{-N} \, T \owedge g + \frac{N(N-1)}{4} \, \Upsilon^{-N} \, T \owedge T.\] 
\item The Ricci tensor of $(M,g)$ at the point $p$ is given by 
\[-(n-1) \, g - (N-n+1) \, \Upsilon^{-N} \, g + \frac{1}{2} \, N \, (N-n+1) \, \Upsilon^{-N} \, T.\] 
\end{itemize}
\end{proposition}

\textbf{Proof.} 
By Proposition \ref{union.of.Sigma_t}, we can find an element $\bar{t} \in S^1$ such that $p \in \Sigma_{\bar{t}}$. By Proposition \ref{existence.of.totally.geodesic.hypersurface.with.given.asymptotics}, $(\Sigma_{\bar{t}},\check{g}_{\bar{t}},\check{\rho}_{\bar{t}})$ is a model $(N,n-1)$-dataset. In view of Proposition \ref{properties.of.HM.metrics}, there exists a symmetric bilinear form $\check{T}: T_p \Sigma_{\bar{t}} \times T_p \Sigma_{\bar{t}} \to \mathbb{R}$ and a real number $\Upsilon$ with the following properties: 
\begin{itemize}
\item The eigenvalues of $\check{T}$ are $1$ and $0$, and the corresponding multiplicities are $2$ and $n-3$, respectively. 
\item The Hessian of the function $\check{\rho}_{\bar{t}}^{\frac{1}{N-n+1}}: \Sigma_{\bar{t}} \to \mathbb{R}$ at the point $p$ is given by 
\[D^{\Sigma_{\bar{t}},2} \check{\rho}_{\bar{t}}^{\frac{1}{N-n+1}} = \check{\rho}_{\bar{t}}^{\frac{1}{N-n+1}} \, (1-\Upsilon^{-N}) \, \check{g}_{\bar{t}} + \frac{N}{2} \, \check{\rho}_{\bar{t}}^{\frac{1}{N-n+1}} \, \Upsilon^{-N} \, \check{T}.\] 
\item The Riemann curvature tensor of $\Sigma_{\bar{t}}$ at the point $p$ is given by 
\[-\frac{1}{2} \, (1-\Upsilon^{-N}) \, \check{g}_{\bar{t}} \owedge \check{g}_{\bar{t}} - \frac{N}{2} \, \Upsilon^{-N} \, \check{T} \owedge \check{g}_{\bar{t}} + \frac{N(N-1)}{4} \, \Upsilon^{-N} \, \check{T} \owedge \check{T}.\] 
\end{itemize}
Let $\{e_1,\hdots,e_{n-1}\}$ be an orthonormal basis of $T_p \Sigma_{\bar{t}}$. Since $\Sigma_{\bar{t}}$ is totally geodesic, the Gauss equations imply that $R(e_i,e_j,e_k,e_l) = R_{\Sigma_{\bar{t}}}(e_i,e_j,e_k,e_l)$ for all $i,j,k,l \in \{1,\hdots,n-1\}$. Therefore, 
\begin{align} 
\label{component.a}
&R(e_i,e_j,e_k,e_l) \notag \\ 
&= -(1-\Upsilon^{-N}) \, (\delta_{ik} \, \delta_{jl} - \delta_{il} \, \delta_{jk}) \notag \\ 
&- \frac{N}{2} \, \Upsilon^{-N} \, (\check{T}(e_i,e_k) \, \delta_{jl} - \check{T}(e_i,e_l) \, \delta_{jk} - \check{T}(e_j,e_k) \, \delta_{il} + \check{T}(e_j,e_l) \, \delta_{ik}) \\ 
&+ \frac{N(N-1)}{2} \, \Upsilon^{-N} \, (\check{T}(e_i,e_k) \, \check{T}(e_j,e_l) - \check{T}(e_i,e_l) \, \check{T}(e_j,e_k)) \notag
\end{align}
for all $i,j,k,l \in \{1,\hdots,n-1\}$. Moreover, since $\Sigma_{\bar{t}}$ is totally geodesic, the Codazzi equations imply that 
\begin{equation} 
\label{component.b}
R(e_i,e_j,e_k,\nu_{\Sigma_{\bar{t}}}) = 0 
\end{equation} 
for all $i,j,k \in \{1,\hdots,n-1\}$. 

By Proposition \ref{existence.of.totally.geodesic.hypersurface.with.given.asymptotics}, the function $\check{v}_{\bar{t}}^{-(N-n)} \, \rho$ is constant along $\Sigma_{\bar{t}}$. Moreover, it follows from the definition of $\check{\rho}_{\bar{t}}$ that the function $\check{v}_{\bar{t}}^{-1} \, \rho^{-1} \, \check{\rho}_{\bar{t}}$ is constant along $\Sigma_{\bar{t}}$. Consequently, the function $\check{v}_{\bar{t}}^{-(N-n+1)} \, \check{\rho}_{\bar{t}}$ is constant along $\Sigma_{\bar{t}}$. Using Corollary \ref{Hessian.of.check.v}, we obtain 
\[(D^{\Sigma_{\bar{t}},2} \check{\rho}_{\bar{t}}^{\frac{1}{N-n+1}})(e_i,e_k) + R(e_i,\nu_{\Sigma_{\bar{t}}},e_k,\nu_{\Sigma_{\bar{t}}}) \, \check{\rho}_{\bar{t}}^{\frac{1}{N-n+1}} = 0\] 
for all $i,k \in \{1,\hdots,n-1\}$. This implies 
\begin{equation} 
\label{component.c}
R(e_i,\nu_{\Sigma_{\bar{t}}},e_k,\nu_{\Sigma_{\bar{t}}}) = -(1-\Upsilon^{-N}) \, \delta_{ik} - \frac{N}{2} \, \Upsilon^{-N} \, \check{T}(e_i,e_k) 
\end{equation}
for all $i,k \in \{1,\hdots,n-1\}$. Combining (\ref{component.a}), (\ref{component.b}), and (\ref{component.c}), it follows that the Riemann curvature tensor of $(M,g)$ at the point $p$ is given by 
\[-\frac{1}{2} \, (1-\Upsilon^{-N}) \, g \owedge g - \frac{N}{2} \, \Upsilon^{-N} \, T \owedge g + \frac{N(N-1)}{4} \, \Upsilon^{-N} \, T \owedge T,\] 
where $T: T_p M \times T_p M \to \mathbb{R}$ denotes the trivial extension of $\check{T}: T_p \Sigma_{\bar{t}} \times T_p \Sigma_{\bar{t}} \to \mathbb{R}$. The formula for the Ricci tensor of $(M,g)$ at the point $p$ follows by taking the trace. This completes the proof of Proposition \ref{curvature.tensor.of.M}. \\

By Proposition \ref{curvature.tensor.of.M}, the norm of the traceless Ricci tensor of $(M,g)$ is given by 
\[\sqrt{\frac{n-2}{2n}} \, N \, (N-n+1) \, \Upsilon^{-N}.\] 
Since $3 \leq n \leq N$, we conclude that $\Upsilon$ defines a smooth function on $M$, which takes values in the interval $[1,\infty)$. Moreover, it follows from Proposition \ref{curvature.tensor.of.M} that the traceless Ricci tensor of $(M,g)$ is given by 
\[\frac{1}{2n} \, N \, (N-n+1) \, \Upsilon^{-N} \, (n \, T - 2 \, g).\] 
Since $3 \leq n \leq N$, we conclude that $T$ defines a smooth $(0,2)$-tensor field on $M$. For each point $p \in M$, the tangent space $T_p M$ can be decomposed as a direct sum $T_p M = \mathcal{E}_p \oplus \mathcal{F}_p$, where $\mathcal{E}_p$ denotes the eigenspace of $T$ with eigenvalue $1$ and $\mathcal{F}_p$ denotes the eigenspace of $T$ with eigenvalue $0$. Clearly, $\mathcal{E}$ is a smooth subbundle of $TM$ of rank $2$ and $\mathcal{F}$ is a smooth subbundle of $TM$ of rank $n-2$. Note that $\mathcal{E}_p \subset T_p \Sigma_{\bar{t}}$ whenever $\bar{t} \in S^1$ and $p \in \Sigma_{\bar{t}}$.

\begin{lemma} 
\label{2D.submanifold.in.M}
There exists a smooth map $\Phi_{\text{\rm 2D}}: \mathbb{R}^2 \to M$ such that $\Phi_{\text{\rm 2D}}^* g = g_{\text{\rm HM},N,2}$ and $\Upsilon \circ \Phi_{\text{\rm 2D}} = \Upsilon_{\text{\rm HM},N,2}$. Moreover, for each point in $\mathbb{R}^2$, the differential of $\Phi_{\text{\rm 2D}}$ takes values in the subbundle $\mathcal{E}$.
\end{lemma}

\textbf{Proof.} 
We consider an arbitrary element $\bar{t} \in S^1$. Since $(\Sigma_{\bar{t}},\check{g}_{\bar{t}},\check{\rho}_{\bar{t}})$ is a model $(N,n-1)$-dataset, we can find a smooth map $\check{\Phi}: \mathbb{R}^2 \times \mathbb{R}^{n-3} \to \Sigma_{\bar{t}}$ such that $\check{\Phi}^* \check{g}_{\bar{t}} = g_{\text{\rm HM},N,n-1}$. This implies $R_{\Sigma_{\bar{t}}} \circ \check{\Phi} = -(n-1)(n-2) + (N-n+2)(N-n+1) \, \Upsilon_{\text{\rm HM},N,n-1}^{-N}$, where $R_{\Sigma_{\bar{t}}}$ denotes the scalar curvature of $\Sigma_{\bar{t}}$. On the other hand, using Proposition \ref{curvature.tensor.of.M} and the fact that $\Sigma_{\bar{t}}$ is totally geodesic, we obtain $R_{\Sigma_{\bar{t}}} = -(n-1)(n-2) + (N-n+2)(N-n+1) \, \Upsilon^{-N}$ at each point on $\Sigma_{\bar{t}}$. Since $3 \leq n \leq N$, it follows that $\Upsilon \circ \check{\Phi} = \Upsilon_{\text{\rm HM},N,n-1}$. Finally, we define $\Phi_{\text{\rm 2D}}$ to be the restriction of the map $\check{\Phi}: \mathbb{R}^2 \times \mathbb{R}^{n-3} \to \Sigma_{\bar{t}}$ to $\mathbb{R}^2 \times \{0\} \subset \mathbb{R}^2 \times \mathbb{R}^{n-3}$. This completes the proof of Lemma \ref{2D.submanifold.in.M}. \\

\begin{lemma}
\label{normal.derivative.of.Upsilon.and.T}
Let us fix an element $\bar{t} \in S^1$. Then $\langle \nabla \Upsilon,\nu_{\Sigma_{\bar{t}}} \rangle = 0$ at each point on $\Sigma_{\bar{t}}$. Moreover, if $\{e_1,\hdots,e_{n-1}\}$ is a local orthonormal frame on $\Sigma_{\bar{t}}$, then $(D_{\nu_{\Sigma_{\bar{t}}}} T)(e_i,e_k) = 0$ for all $i,k \in \{1,\hdots,n-1\}$.
\end{lemma}

\textbf{Proof.}
Let us consider a sequence $t_j \in S^1$ such that $t_j \neq \bar{t}$ for each $j$ and $t_j \to \bar{t}$ as $j \to \infty$. Let $\delta_j$, $w^{(j)}$, and $w$ be defined as in Proposition \ref{Jacobi.field}. By Proposition \ref{Jacobi.field}, $w$ is a non-zero multiple of $\check{v}_{\bar{t}}$. In particular, $w$ is non-zero at each point on $\Sigma_{\bar{t}}$.

For each $j$, we define a smooth map $\Psi^{(j)}: \Sigma_{\bar{t}} \setminus \{r>\delta_j^{-1}\} \to \Sigma_{t_j}$ by $\Psi^{(j)}(x) = \exp_x(w^{(j)}(x) \, \nu_{\Sigma_{\bar{t}}}(x))$ for $x \in \Sigma_{\bar{t}} \setminus \{r>\delta_j^{-1}\}$. Let $(\Psi^{(j)})^* \check{g}_{t_j}$ denote the pull-back of the metric $\check{g}_{t_j}$ under the map $\Psi^{(j)}$. Recall that $d_{S^1}(t_j,\bar{t})^{-1} \, w^{(j)} \to w$ in $C_{\text{\rm loc}}^\infty$. Since $\Sigma_{\bar{t}}$ is totally geodesic, it follows that 
\begin{equation} 
\label{metric.unchanged.to.leading.order}
d_{S^1}(t_j,\bar{t})^{-1} \, \big ( (\Psi^{(j)})^* \check{g}_{t_j} - \check{g}_{\bar{t}} \big ) \to 0 
\end{equation}
in $C_{\text{\rm loc}}^\infty$. Using (\ref{metric.unchanged.to.leading.order}), we obtain 
\begin{equation} 
\label{scalar.curvature.unchanged.to.leading.order}
d_{S^1}(t_j,\bar{t})^{-1} \, (R_{\Sigma_{t_j}} \circ \Psi^{(j)} - R_{\Sigma_{\bar{t}}}) \to 0 
\end{equation}
at each point on $\Sigma_{\bar{t}}$. On the other hand, for each $t \in S^1$, the scalar curvature of $\Sigma_t$ is given by $R_{\Sigma_t} = -(n-1)(n-2) + (N-n+2)(N-n+1) \, \Upsilon^{-N}$. Since $3 \leq n \leq N$, the relation (\ref{scalar.curvature.unchanged.to.leading.order}) implies that 
\begin{equation} 
\label{Upsilon.unchanged.to.leading.order}
d_{S^1}(t_j,\bar{t})^{-1} \, (\Upsilon^{-N} \circ \Psi^{(j)} - \Upsilon^{-N}) \to 0 
\end{equation}
at each point on $\Sigma_{\bar{t}}$. Thus, we conclude that $\langle \nabla \Upsilon,w \, \nu_{\Sigma_{\bar{t}}} \rangle = 0$ at each point on $\Sigma_{\bar{t}}$. 

Using (\ref{metric.unchanged.to.leading.order}) again, we obtain 
\begin{equation} 
\label{Ricci.tensor.unchanged.to.leading.order}
d_{S^1}(t_j,\bar{t})^{-1} \, (\text{\rm Ric}_{\Sigma_{t_j}}(\Psi^{(j)}_* e_i,\Psi^{(j)}_* e_k) - \text{\rm Ric}_{\Sigma_{\bar{t}}}(e_i,e_k)) \to 0 
\end{equation} 
at each point on $\Sigma_{\bar{t}}$. On the other hand, for each $t \in S^1$, the Ricci tensor of $\Sigma_t$ is given by the restriction of the tensor $-(n-2) \, g - (N-n+2) \, \Upsilon^{-N} \, g + \frac{1}{2} \, N \, (N-n+2) \, \Upsilon^{-N} \, T$ to the tangent bundle of $\Sigma_t$. Since $3 \leq n \leq N$, the relations (\ref{metric.unchanged.to.leading.order}), (\ref{Upsilon.unchanged.to.leading.order}), and (\ref{Ricci.tensor.unchanged.to.leading.order}) imply that 
\begin{equation} 
\label{T.unchanged.to.leading.order}
d_{S^1}(t_j,\bar{t})^{-1} \, (T(\Psi^{(j)}_* e_i,\Psi^{(j)}_* e_k) - T(e_i,e_k)) \to 0 
\end{equation}
at each point on $\Sigma_{\bar{t}}$. Since $T$ is a smooth $(0,2)$-tensor field on the ambient manifold $M$, the relation (\ref{T.unchanged.to.leading.order}) gives 
\[(D_{w \, \nu_{\Sigma_{\bar{t}}}} T)(e_i,e_k) + T(D_{e_i}(w \, \nu_{\Sigma_{\bar{t}}}),e_k) + T(e_i,D_{e_k}(w \, \nu_{\Sigma_{\bar{t}}})) = 0\] 
at each point on $\Sigma_{\bar{t}}$. The last two terms on the left hand side vanish since $\Sigma_{\bar{t}}$ is totally geodesic and $T(\cdot,\nu_{\Sigma_{\bar{t}}}) = 0$. Thus, $(D_{w \, \nu_{\Sigma_{\bar{t}}}} T)(e_i,e_k) = 0$ at each point on $\Sigma_{\bar{t}}$. This completes the proof of Lemma \ref{normal.derivative.of.Upsilon.and.T}. \\

\begin{lemma}
\label{norm.of.gradient.of.Upsilon}
We have $|\nabla \Upsilon|^2 = \Upsilon^2 \, (1-\Upsilon^{-N})$ at each point on $M$.
\end{lemma}

\textbf{Proof.} 
Let us consider an arbitrary point $p \in M$. By Proposition \ref{union.of.Sigma_t}, we can find an element $\bar{t} \in S^1$ such that $p \in \Sigma_{\bar{t}}$. Since $(\Sigma_{\bar{t}},\check{g}_{\bar{t}},\check{\rho}_{\bar{t}})$ is a model $(N,n-1)$-dataset, we know that $|\nabla^{\Sigma_{\bar{t}}} \Upsilon|^2 = \Upsilon^2 \, (1-\Upsilon^{-N})$ at each point on $\Sigma_{\bar{t}}$. On the other hand, Lemma \ref{normal.derivative.of.Upsilon.and.T} implies $\langle \nabla \Upsilon,\nu_{\Sigma_{\bar{t}}} \rangle = 0$ at each point on $\Sigma_{\bar{t}}$. Putting these facts together, we conclude that $|\nabla \Upsilon|^2 = \Upsilon^2 \, (1-\Upsilon^{-N})$ at each point on $\Sigma_{\bar{t}}$. This completes the proof of Lemma \ref{norm.of.gradient.of.Upsilon}. \\

\begin{lemma}
\label{gradient.of.Upsilon}
Let $p \in M$, and let $\{e_1,\hdots,e_n\}$ be an orthonormal basis of $T_p M$ such that $\mathcal{E}_p = \text{\rm span}\{e_1,e_2\}$. Then $\langle \nabla \Upsilon,e_k \rangle = 0$ for all $k \in \{3,\hdots,n\}$. 
\end{lemma}

\textbf{Proof.} 
By Proposition \ref{union.of.Sigma_t}, we can find an element $\bar{t} \in S^1$ such that $p \in \Sigma_{\bar{t}}$. Note that $\mathcal{E}_p \subset T_p \Sigma_{\bar{t}}$. Let $\{e_1,\hdots,e_n\}$ be an orthonormal basis of $T_p M$ such that $\mathcal{E}_p = \text{\rm span}\{e_1,e_2\}$ and $\nu_{\Sigma_{\bar{t}}} = e_n$. To prove the assertion, we distinguish two cases: 

\textit{Case 1:} We first consider the case $k \in \{3,\hdots,n-1\}$. Since $(\Sigma_{\bar{t}},\check{g}_{\bar{t}},\check{\rho}_{\bar{t}})$ is a model $(N,n-1)$-dataset, it follows that $\langle \nabla \Upsilon,e_k \rangle = 0$ for all $k \in \{3,\hdots,n-1\}$. 

\textit{Case 2:} We now consider the case $k=n$. Using Lemma \ref{normal.derivative.of.Upsilon.and.T}, we obtain $\langle \nabla \Upsilon,e_n \rangle = 0$. This completes the proof of Lemma \ref{gradient.of.Upsilon}. \\

\begin{lemma}
\label{covariant.derivative.of.T}
Let $p \in M$, and let $\{e_1,\hdots,e_n\}$ be an orthonormal basis of $T_p M$ such that $\mathcal{E}_p = \text{\rm span}\{e_1,e_2\}$. Then 
\[(D_{e_l} T)(e_i,e_k) = \Upsilon^{-1} \, \langle \nabla \Upsilon,e_i \rangle \, \delta_{kl}\] 
for all $i \in \{1,2\}$, $k \in \{3,\hdots,n\}$, and $l \in \{1,\hdots,n\}$.
\end{lemma}

\textbf{Proof.} 
By Proposition \ref{union.of.Sigma_t}, we can find an element $\bar{t} \in S^1$ such that $p \in \Sigma_{\bar{t}}$. Note that $\mathcal{E}_p \subset T_p \Sigma_{\bar{t}}$. Let $\{e_1,\hdots,e_n\}$ be an orthonormal basis of $T_p M$ such that $\mathcal{E}_p = \text{\rm span}\{e_1,e_2\}$ and $\nu_{\Sigma_{\bar{t}}} = e_n$. To prove the assertion, we distinguish four cases: 

\textit{Case 1:} We first consider the case when $k \in \{3,\hdots,n-1\}$ and $l \in \{1,\hdots,n-1\}$. Let $\check{T}$ denote the restriction of $T$ to the tangent bundle of $\Sigma_{\bar{t}}$. Since $(\Sigma_{\bar{t}},\check{g}_{\bar{t}},\check{\rho}_{\bar{t}})$ is a model $(N,n-1)$-dataset, we know that 
\[(D_{e_l}^{\Sigma_{\bar{t}}} \check{T})(e_i,e_k) = \Upsilon^{-1} \, \langle \nabla^{\Sigma_{\bar{t}}} \Upsilon,e_i \rangle \, \delta_{kl}\] 
for all $i \in \{1,2\}$, $k \in \{3,\hdots,n-1\}$, and $l \in \{1,\hdots,n-1\}$. Since $\Sigma_{\bar{t}}$ is totally geodesic, it follows that 
\[(D_{e_l} T)(e_i,e_k) = \Upsilon^{-1} \, \langle \nabla \Upsilon,e_i \rangle \, \delta_{kl}\] 
for all $i \in \{1,2\}$, $k \in \{3,\hdots,n-1\}$, and $l \in \{1,\hdots,n-1\}$.

\textit{Case 2:} We next consider the case when $k=n$ and $l \in \{1,\hdots,n-1\}$. At each point on $\Sigma_{\bar{t}}$, we have $T(\cdot,\nu_{\Sigma_{\bar{t}}}) = 0$. We differentiate this identity in tangential direction. Since $\Sigma_{\bar{t}}$ is totally geodesic, we conclude that $(D_{e_l} T)(e_i,\nu_{\Sigma_{\bar{t}}}) = 0$ for all $i \in \{1,2\}$ and $l \in \{1,\hdots,n-1\}$.

\textit{Case 3:} We now consider the case when $k \in \{3,\hdots,n-1\}$ and $l=n$. Using Lemma \ref{normal.derivative.of.Upsilon.and.T}, we obtain $(D_{e_n} T)(e_i,e_k) = 0$ for all $i \in \{1,2\}$ and $k \in \{3,\hdots,n-1\}$. 

\textit{Case 4:} Finally, we consider the case when $k=n$ and $l=n$. By Proposition \ref{curvature.tensor.of.M}, the Ricci tensor of $(M,g)$ is given by $-(n-1) \, g - (N-n+1) \, \Upsilon^{-N} \, g + \frac{1}{2} \, N \, (N-n+1) \, \Upsilon^{-N} \, T$. Using the contracted second Bianchi identity on $M$, we obtain 
\begin{equation} 
\label{Bianchi.1} 
\sum_{m=1}^n (D_{e_m} T)(e_i,e_m) = (n-2) \, \Upsilon^{-1} \, \langle \nabla \Upsilon,e_i \rangle 
\end{equation}
for all $i \in \{1,2\}$. We next observe that the Ricci tensor of $\Sigma_{\bar{t}}$ is given by $-(n-2) \, \check{g}_{\bar{t}} - (N-n+2) \, \Upsilon^{-N} \, \check{g}_{\bar{t}} + \frac{1}{2} \, N \, (N-n+2) \, \Upsilon^{-N} \, \check{T}$, where $\check{T}$ denotes the restriction of $T$ to the tangent bundle of $\Sigma_{\bar{t}}$. Using the contracted second Bianchi identity on $\Sigma_{\bar{t}}$, we obtain 
\begin{equation} 
\label{Bianchi.2}
\sum_{m=1}^{n-1} (D_{e_m}^{\Sigma_{\bar{t}}} \check{T})(e_i,e_m) = (n-3) \, \Upsilon^{-1} \, \langle \nabla^{\Sigma_{\bar{t}}} \Upsilon,e_i \rangle 
\end{equation}
for all $i \in \{1,2\}$. Since $\Sigma_{\bar{t}}$ is totally geodesic, the identity (\ref{Bianchi.2}) can be rewritten as 
\begin{equation} 
\label{Bianchi.3}
\sum_{m=1}^{n-1} (D_{e_m} T)(e_i,e_m) = (n-3) \, \Upsilon^{-1} \, \langle \nabla \Upsilon,e_i \rangle 
\end{equation}
for all $i \in \{1,2\}$. Subtracting (\ref{Bianchi.3}) from (\ref{Bianchi.1}), we conclude that 
\[(D_{e_n} T)(e_i,e_n) = \Upsilon^{-1} \, \langle \nabla \Upsilon,e_i \rangle\] 
for all $i \in \{1,2\}$. This completes the proof of Lemma \ref{covariant.derivative.of.T}. \\

\begin{proposition}
\label{connection.coefficients}
Suppose that $\{e_1,\hdots,e_n\}$ is a local orthonormal frame on $M$ such that $\mathcal{E} = \text{\rm span}\{e_1,e_2\}$. Then 
\[\langle D_{e_l} e_k,e_i \rangle = -\Upsilon^{-1} \, \langle \nabla \Upsilon,e_i \rangle \, \delta_{kl}\] 
for all $i \in \{1,2\}$, $k \in \{3,\hdots,n\}$, and $l \in \{1,\hdots,n\}$. 
\end{proposition} 

\textbf{Proof.} 
Note that $T(e_i,e_k) = 0$ for all $i \in \{1,2\}$ and $k \in \{3,\hdots,n\}$. Differentiating this identity gives 
\[(D_{e_l} T)(e_i,e_k) + T(e_i,D_{e_l} e_k) + T(D_{e_l} e_i,e_k) = 0\] 
for all $i \in \{1,2\}$, $k \in \{3,\hdots,n\}$, and $l \in \{1,\hdots,n\}$. Since $T(\cdot,e_k) = 0$ for $k \in \{3,\hdots,n\}$, it follows that 
\begin{equation} 
\label{derivative.of.T}
(D_{e_l} T)(e_i,e_k) + T(e_i,D_{e_l} e_k) = 0 
\end{equation}
for all $i \in \{1,2\}$, $k \in \{3,\hdots,n\}$, and $l \in \{1,\hdots,n\}$. Using the identity (\ref{derivative.of.T}) together with Lemma \ref{covariant.derivative.of.T}, we obtain 
\[\langle e_i,D_{e_l} e_k \rangle = T(e_i,D_{e_l} e_k) = -(D_{e_l} T)(e_i,e_k) = -\Upsilon^{-1} \, \langle \nabla \Upsilon,e_i \rangle \, \delta_{kl}\] 
for all $i \in \{1,2\}$, $k \in \{3,\hdots,n\}$, and $l \in \{1,\hdots,n\}$. This completes the proof of Proposition \ref{connection.coefficients}. \\

\begin{proposition}
\label{Hessian.of.Upsilon}
Suppose that $\{e_1,\hdots,e_n\}$ is a local orthonormal frame on $M$ such that $\mathcal{E} = \text{\rm span}\{e_1,e_2\}$. Then 
\[(D^2 \Upsilon)(e_k,e_l) = \Upsilon \, (1-\Upsilon^{-N}) \, \delta_{kl}\] 
for all $k,l \in \{3,\hdots,n\}$. 
\end{proposition} 

\textbf{Proof.} 
Lemma \ref{gradient.of.Upsilon} implies that $\langle \nabla \Upsilon,e_k \rangle = 0$ for all $k \in \{3,\hdots,n\}$. Differentiating this identity gives 
\begin{equation} 
\label{Hessian.of.Upsilon.2}
(D^2 \Upsilon)(e_k,e_l) + \langle \nabla \Upsilon,D_{e_l} e_k \rangle = 0. 
\end{equation}
for all $k,l \in \{3,\hdots,n\}$. Using the identity (\ref{Hessian.of.Upsilon.2}) together with Lemma \ref{gradient.of.Upsilon} and Proposition \ref{connection.coefficients}, we obtain 
\[(D^2 \Upsilon)(e_k,e_l) = -\sum_{i=1}^2 \langle \nabla \Upsilon,e_i \rangle \, \langle D_{e_l} e_k,e_i \rangle = \Upsilon^{-1} \, |\nabla \Upsilon|^2 \, \delta_{kl}\]
for all $k,l \in \{3,\hdots,n\}$. The assertion follows now from Lemma \ref{norm.of.gradient.of.Upsilon}. This completes the proof of Proposition \ref{Hessian.of.Upsilon}. \\

\begin{corollary} 
\label{parallel.subbundles}
The bundles $\mathcal{E}$ and $\mathcal{F}$ are invariant under parallel transport with respect to the metric $\Upsilon^{-2} \, g$. 
\end{corollary}

\textbf{Proof.} 
Suppose that $\{e_1,\hdots,e_n\}$ is a local orthonormal frame on $M$ such that $\mathcal{E} = \text{\rm span}\{e_1,e_2\}$. Proposition \ref{connection.coefficients} implies that 
\[D_{e_l} e_k + \delta_{kl} \, \Upsilon^{-1} \, \nabla \Upsilon \in \mathcal{F}\]
for all $k \in \{3,\hdots,n\}$ and $l \in \{1,\hdots,n\}$. Using Lemma \ref{gradient.of.Upsilon}, we obtain 
\begin{equation} 
\label{modified.covariant.derivative}
D_{e_l} e_k + \delta_{kl} \, \Upsilon^{-1} \, \nabla \Upsilon - \Upsilon^{-1} \, \langle \nabla \Upsilon,e_k \rangle \, e_l - \Upsilon^{-1} \, \langle \nabla \Upsilon,e_l \rangle \, e_k \in \mathcal{F} 
\end{equation}
for all $k \in \{3,\hdots,n\}$ and $l \in \{1,\hdots,n\}$. The expression in (\ref{modified.covariant.derivative}) is equal to the covariant derivative of $e_k$ along $e_l$ with respect to the metric $\Upsilon^{-2} \, g$ (see \cite{Besse}, Theorem 1.159). This completes the proof of Corollary \ref{parallel.subbundles}. \\

\begin{corollary}
\label{curvature.F}
The restriction of the Riemann curvature tensor of the metric $\Upsilon^{-2} \, g$ to the subbundle $\mathcal{F}$ vanishes.
\end{corollary}

\textbf{Proof.} 
Suppose that $\{e_1,\hdots,e_n\}$ is a local orthonormal frame on $M$ such that $\mathcal{E} = \text{\rm span}\{e_1,e_2\}$. Proposition \ref{curvature.tensor.of.M} gives 
\[R(e_i,e_j,e_k,e_l) = -(1-\Upsilon^{-N}) \, (\delta_{ik} \, \delta_{jl} - \delta_{il} \, \delta_{jk})\] 
for all $i,j,k,l \in \{3,\hdots,n\}$. Using Lemma \ref{norm.of.gradient.of.Upsilon} and Proposition \ref{Hessian.of.Upsilon}, we obtain 
\begin{align*} 
&R(e_i,e_j,e_k,e_l) - \Upsilon^{-2} \, |\nabla \Upsilon|^2 \, (\delta_{ik} \, \delta_{jl} - \delta_{il} \, \delta_{jk}) \\ 
&+ \Upsilon^{-1} \, (D^2 \Upsilon)(e_i,e_k) \, \delta_{jl} - \Upsilon^{-1} \, (D^2 \Upsilon)(e_i,e_l) \, \delta_{jk} \\ 
&- \Upsilon^{-1} \, (D^2 \Upsilon)(e_j,e_k) \, \delta_{il} + \Upsilon^{-1} \, (D^2 \Upsilon)(e_j,e_l) \, \delta_{ik} = 0 
\end{align*}
for all $i,j,k,l \in \{3,\hdots,n\}$. The assertion follows now from the standard formula for the change of the Riemann curvature tensor under a conformal change of the metric (see \cite{Besse}, Theorem 1.159). This completes the proof of Corollary \ref{curvature.F}. \\

\begin{proposition}
\label{weight}
The function $\Upsilon^{-(N-n)} \, \rho$ is constant on $M$.
\end{proposition}

\textbf{Proof.} 
It suffices to show that the gradient of $\Upsilon^{-(N-n)} \, \rho$ vanishes identically. To prove this, we fix an arbitrary point $p \in M$. By Proposition \ref{union.of.Sigma_t}, we can find an element $\bar{t} \in S^1$ such that $p \in \Sigma_{\bar{t}}$. By Proposition \ref{existence.of.totally.geodesic.hypersurface.with.given.asymptotics}, the function $\check{v}_{\bar{t}}^{-(N-n)} \, \rho$ is constant along $\Sigma_{\bar{t}}$. Moreover, it follows from the definition of $\check{\rho}_{\bar{t}}$ that the function $\check{v}_{\bar{t}} \, \rho \, \check{\rho}_{\bar{t}}^{-1}$ is constant along $\Sigma_{\bar{t}}$. Consequently, the function 
\[\rho^{N-n+1} \, \check{\rho}_{\bar{t}}^{-(N-n)} = \check{v}_{\bar{t}}^{-(N-n)} \, \rho \, (\check{v}_{\bar{t}} \, \rho \, \check{\rho}_{\bar{t}}^{-1})^{N-n}\] 
is constant along $\Sigma_{\bar{t}}$. On the other hand, since $(\Sigma_{\bar{t}},\check{g}_{\bar{t}},\check{\rho}_{\bar{t}})$ is a model $(N,n-1)$-dataset, the function $\Upsilon^{-(N-n+1)} \, \check{\rho}_{\bar{t}}$ is constant along $\Sigma_{\bar{t}}$. Putting these facts together, we conclude that the function 
\[\Upsilon^{-(N-n)} \, \rho = (\Upsilon^{-(N-n+1)} \, \check{\rho}_{\bar{t}})^{\frac{N-n}{N-n+1}} \, (\rho^{N-n+1} \, \check{\rho}_{\bar{t}}^{-(N-n)})^{\frac{1}{N-n+1}}\] 
is constant along $\Sigma_{\bar{t}}$. Thus, 
\begin{equation} 
\label{tangential.component.of.gradient}
\nabla^{\Sigma_{\bar{t}}}(\Upsilon^{-(N-n)} \, \rho) = 0 
\end{equation}
at each point on $\Sigma_{\bar{t}}$. On the other hand, Lemma \ref{normal.derivative.of.Upsilon.and.T} implies that $\langle \nabla \Upsilon,\nu_{\Sigma_{\bar{t}}} \rangle = 0$ at each point on $\Sigma_{\bar{t}}$. Moreover, Proposition \ref{existence.of.totally.geodesic.hypersurface.with.given.asymptotics} gives $\langle \nabla \rho,\nu_{\Sigma_{\bar{t}}} \rangle = 0$ at each point on $\Sigma_{\bar{t}}$. Therefore, 
\begin{equation} 
\label{normal.component.of.gradient}
\langle \nabla(\Upsilon^{-(N-n)} \, \rho),\nu_{\Sigma_{\bar{t}}} \rangle = 0 
\end{equation}
at each point on $\Sigma_{\bar{t}}$. Combining (\ref{tangential.component.of.gradient}) and (\ref{normal.component.of.gradient}), we obtain $\nabla(\Upsilon^{-(N-n)} \, \rho) = 0$ at each point on $\Sigma_{\bar{t}}$. This completes the proof of Proposition \ref{weight}. \\

Using Corollary \ref{parallel.subbundles} and de Rham's decomposition theorem (see \cite{Besse}, Theorem 10.43), we conclude that the universal cover of $(M,\Upsilon^{-2} \, g)$ is isometric to a product of a two-dimensional manifold (corresponding to the subbundle $\mathcal{E}$) with an $(n-2)$-dimensional manifold (corresponding to the subbundle $\mathcal{F}$). Corollary \ref{curvature.F} implies that the second factor is flat. By Lemma \ref{gradient.of.Upsilon}, the directional derivative of $\Upsilon$ along every vector in $\mathcal{F}$ vanishes.

We next consider the map $\Phi_{\text{\rm 2D}}: \mathbb{R}^2 \to M$ constructed in Lemma \ref{2D.submanifold.in.M}. Then $\Phi_{\text{\rm 2D}}^*(\Upsilon^{-2} \, g) = \Upsilon_{\text{\rm HM},N,2}^{-2} \, g_{\text{\rm HM},N,2}$ and $\Upsilon \circ \Phi_{\text{\rm 2D}} = \Upsilon_{\text{\rm HM},N,2}$. Moreover, for each point in $\mathbb{R}^2$, the differential of $\Phi_{\text{\rm 2D}}$ takes values in the subbundle $\mathcal{E}$. Using the normal exponential map of $\Phi_{\text{\rm 2D}}$ with respect to the metric $\Upsilon^{-2} \, g$, we can construct a smooth map $\Phi: \mathbb{R}^2 \times \mathbb{R}^{n-2} \to M$ such that $\Phi^*(\Upsilon^{-2} \, g) = \Upsilon_{\text{\rm HM},N,n}^{-2} \, g_{\text{\rm HM},N,n}$ and $\Upsilon \circ \Phi = \Upsilon_{\text{\rm HM},N,n}$. Putting these facts together, we conclude that $\Phi^* g = g_{\text{\rm HM},N,n}$. By Proposition \ref{weight}, the function $\rho \circ \Phi$ is a constant multiple of $\Upsilon^{N-n} \circ \Phi = \Upsilon_{\text{\rm HM},N,n}^{N-n} = \rho_{\text{\rm HM},N,n}$. Thus, $(M,g,\rho)$ is a model $(N,n)$-dataset. This completes the proof of Theorem \ref{property.star.holds}.

\appendix

\section{Some auxiliary identities}

\label{key.identities}

In this appendix, we derive several identities involving the weighted Jacobi operator and the second variation of the $(g,\rho)$-area. The calculations are lengthy, but standard; see \cite{Ambrozio-Carlotto-Sharp} for related work.

\begin{proposition}
\label{formula.weighted.Jacobi.operator} 
Let $(M,g)$ be an orientable Riemannian manifold and let $\rho$ be a positive function on $M$. Let $\Sigma$ be an orientable hypersurface in $M$ satisfying $H_\Sigma + \langle \nabla \log \rho,\nu_\Sigma \rangle = 0$. Let $V$ be a smooth vector field on $M$. We define a function $v$ on $\Sigma$ by $v = \langle V,\nu_\Sigma \rangle$. Then 
\begin{align*} 
&-\text{\rm div}_\Sigma(\rho \, \nabla^\Sigma v) - \rho \, (\text{\rm Ric}(\nu_\Sigma,\nu_\Sigma) + |h_\Sigma|^2) \, v \\ 
&+ (D^2 \rho)(\nu_\Sigma,\nu_\Sigma) \, v - \rho^{-1} \, \langle \nabla \rho,\nu_\Sigma \rangle^2 \, v \\ 
&= -\rho \sum_{k=1}^{n-1} (D_{e_k} (\mathscr{L}_V g))(e_k,\nu_\Sigma) + \frac{1}{2} \, \rho \sum_{k=1}^{n-1} (D_{\nu_\Sigma} (\mathscr{L}_V g))(e_k,e_k) \\ 
&- \rho \sum_{k,l=1}^{n-1} h_\Sigma(e_k,e_l) \, (\mathscr{L}_V g)(e_k,e_l) - (\mathscr{L}_V g)(\nabla \rho,\nu_\Sigma) + \rho \, \big \langle \nabla(V(\log \rho)),\nu_\Sigma \big \rangle 
\end{align*} 
at each point on $\Sigma$.
\end{proposition}

\textbf{Proof.} 
Using the Codazzi equations, we compute 
\begin{align*} 
&-\Delta_\Sigma v - (\text{\rm Ric}(\nu_\Sigma,\nu_\Sigma) + |h_\Sigma|^2) \, v \\ 
&= -\sum_{k=1}^{n-1} \langle D_{e_k,e_k}^2 V,\nu_\Sigma \rangle - \sum_{k=1}^{n-1} R(e_k,V,e_k,\nu_\Sigma) - 2 \sum_{k,l=1}^{n-1} h_\Sigma(e_k,e_l) \, \langle D_{e_k} V,e_l \rangle \\ 
&- \langle V^{\text{\rm tan}},\nabla^\Sigma H_\Sigma \rangle + H_\Sigma \, \langle D_{\nu_\Sigma} V,\nu_\Sigma \rangle 
\end{align*} 
and 
\begin{align*} 
&-\langle \nabla^\Sigma \rho,\nabla^\Sigma v \rangle + (D^2 \rho)(\nu_\Sigma,\nu_\Sigma) \, v - \rho^{-1} \, \langle \nabla \rho,\nu_\Sigma \rangle^2 \, v \\ 
&= -\langle D_{\nabla \rho} V,\nu_\Sigma \rangle - \langle D_{\nu_\Sigma} V,\nabla \rho \rangle + \rho \, \big \langle \nabla(V(\log \rho)),\nu_\Sigma \big \rangle \\ 
&- \rho \, \big \langle V^{\text{\rm tan}},\nabla^\Sigma (\langle \nabla \log \rho,\nu_\Sigma \rangle) \big \rangle + \rho \, \langle \nabla \log \rho,\nu_\Sigma \rangle \, \langle D_{\nu_\Sigma} V,\nu_\Sigma \rangle 
\end{align*}
at each point on $\Sigma$. Using the identity $H_\Sigma + \langle \nabla \log \rho,\nu_\Sigma \rangle = 0$, we obtain 
\begin{align*} 
&-\text{\rm div}_\Sigma(\rho \, \nabla^\Sigma v) - \rho \, (\text{\rm Ric}(\nu_\Sigma,\nu_\Sigma) + |h_\Sigma|^2) \, v \\ 
&+ (D^2 \rho)(\nu_\Sigma,\nu_\Sigma) \, v - \rho^{-1} \, \langle \nabla \rho,\nu_\Sigma \rangle^2 \, v \\ 
&= -\rho \sum_{k=1}^{n-1} \langle D_{e_k,e_k}^2 V,\nu_\Sigma \rangle - \rho \sum_{k=1}^{n-1} R(e_k,V,e_k,\nu_\Sigma) - 2\rho \sum_{k,l=1}^{n-1} h_\Sigma(e_k,e_l) \, \langle D_{e_k} V,e_l \rangle \\ 
&- \langle D_{\nabla \rho} V,\nu_\Sigma \rangle - \langle D_{\nu_\Sigma} V,\nabla \rho \rangle + \rho \, \big \langle \nabla(V(\log \rho)),\nu_\Sigma \big \rangle 
\end{align*} 
at each point on $\Sigma$. From this, the assertion follows easily. \\

\begin{proposition} 
\label{formula.second.variation} 
Let $(M,g)$ be an orientable Riemannian manifold and let $\rho$ be a positive function on $M$. Let $\Sigma$ be an orientable hypersurface in $M$ satisfying $H_\Sigma + \langle \nabla \log \rho,\nu_\Sigma \rangle = 0$. Let $V$ be a smooth vector field on $M$, and let $W = D_V V$. We define a function $v$ on $\Sigma$ by $v = \langle V,\nu_\Sigma \rangle$. Moreover, we define a tangential vector field $Z$ along $\Sigma$ by 
\[Z = D_{V^{\text{\rm tan}}}^\Sigma (V^{\text{\rm tan}}) - \text{\rm div}_\Sigma(V^{\text{\rm tan}}) \, V^{\text{\rm tan}} + 2 \sum_{k=1}^{n-1} h_\Sigma(V^{\text{\rm tan}},e_k) \, \langle V,\nu_\Sigma \rangle \, e_k.\] 
Then 
\begin{align*} 
&\rho \, |\nabla^\Sigma v|^2 - \rho \, (\text{\rm Ric}(\nu_\Sigma,\nu_\Sigma) + |h_\Sigma|^2) \, v^2 + (D^2 \rho)(\nu_\Sigma,\nu_\Sigma) \, v^2 - \rho^{-1} \, \langle \nabla \rho,\nu_\Sigma \rangle^2 \, v^2 \\ 
&+ \text{\rm div}_\Sigma(\rho \, W^{\text{\rm tan}}) - \text{\rm div}_\Sigma(\rho \, Z) + \text{\rm div}_\Sigma(\langle V^{\text{\rm tan}},\nabla^\Sigma \rho \rangle \, V^{\text{\rm tan}}) \\ 
&= \frac{1}{2} \, \rho \sum_{k=1}^{n-1} (\mathscr{L}_V \mathscr{L}_V g)(e_k,e_k) + V(V(\rho)) \\ 
&- \frac{1}{2} \, \rho \sum_{k,l=1}^{n-1} (\mathscr{L}_V g)(e_k,e_l) \, (\mathscr{L}_V g)(e_k,e_l) \\ 
&+ \frac{1}{4} \, \rho \sum_{k,l=1}^{n-1} (\mathscr{L}_V g)(e_k,e_k) \, (\mathscr{L}_V g)(e_l,e_l) \\ 
&+ V(\rho) \sum_{k=1}^{n-1} (\mathscr{L}_V g)(e_k,e_k) 
\end{align*} 
at each point on $\Sigma$. Here, $\{e_1,\hdots,e_{n-1}\}$ denotes a local orthonormal frame on $\Sigma$.
\end{proposition}

\textbf{Proof.} 
We write $Z = Z^{(1)} + Z^{(2)}$, where 
\[Z^{(1)} = D_{V^{\text{\rm tan}}}^\Sigma (V^{\text{\rm tan}}) - \text{\rm div}_\Sigma(V^{\text{\rm tan}}) \, V^{\text{\rm tan}}\] 
and 
\[Z^{(2)} = 2 \sum_{k=1}^{n-1} h_\Sigma(V^{\text{\rm tan}},e_k) \, \langle V,\nu_\Sigma \rangle \, e_k.\] 
Using the Gauss equations and the identity $\langle D_{e_k} V,e_l \rangle = \langle D_{e_k} V^{\text{\rm tan}},e_l \rangle + h_\Sigma(e_k,e_l) \, \langle V,\nu_\Sigma \rangle$ for $k,l \in \{1,\hdots,n-1\}$, we compute 
\begin{align*} 
&\text{\rm div}_\Sigma(Z^{(1)}) \\ 
&= \sum_{k,l=1}^{n-1} \langle D_{e_k} V^{\text{\rm tan}},e_l \rangle \, \langle D_{e_l} V^{\text{\rm tan}},e_k \rangle - \sum_{k,l=1}^{n-1} \langle D_{e_k} V^{\text{\rm tan}},e_k \rangle \, \langle D_{e_l} V^{\text{\rm tan}},e_l \rangle \\ 
&+ \text{\rm Ric}_\Sigma(V^{\text{\rm tan}},V^{\text{\rm tan}}) \\ 
&= \sum_{k,l=1}^{n-1} \langle D_{e_k} V,e_l \rangle \, \langle D_{e_l} V,e_k \rangle - \sum_{k,l=1}^{n-1} \langle D_{e_k} V,e_k \rangle \, \langle D_{e_l} V,e_l \rangle \\ 
&- 2 \sum_{k,l=1}^{n-1} h_\Sigma(e_k,e_l) \, \langle D_{e_k} V^{\text{\rm tan}},e_l \rangle \, \langle V,\nu_\Sigma \rangle + 2 \, H_\Sigma \sum_{k=1}^{n-1} \langle D_{e_k} V,e_k \rangle \, \langle V,\nu_\Sigma \rangle \\ 
&- H_\Sigma^2 \, \langle V,\nu_\Sigma \rangle^2 - |h_\Sigma|^2 \, \langle V,\nu_\Sigma \rangle^2 + H_\Sigma \, h_\Sigma(V^{\text{\rm tan}},V^{\text{\rm tan}}) - h_\Sigma^2(V^{\text{\rm tan}},V^{\text{\rm tan}}) \\ 
&+ \sum_{k=1}^{n-1} R(V^{\text{\rm tan}},e_k,V^{\text{\rm tan}},e_k). 
\end{align*}
Using the Codazzi equations, we obtain 
\begin{align*} 
\text{\rm div}_\Sigma(Z^{(2)}) 
&= 2 \sum_{k,l=1}^{n-1} h_\Sigma(e_k,e_l) \, \langle D_{e_k} V^{\text{\rm tan}},e_l \rangle \, \langle V,\nu_\Sigma \rangle \\ 
&+ 2 \sum_{k=1}^{n-1} h_\Sigma(V^{\text{\rm tan}},e_k) \, \langle D_{e_k} V,\nu_\Sigma \rangle + 2 \, h_\Sigma^2(V^{\text{\rm tan}},V^{\text{\rm tan}}) \\ 
&+ 2 \sum_{k=1}^{n-1} R(V^{\text{\rm tan}},e_k,\nu_\Sigma,e_k) \, \langle V,\nu_\Sigma \rangle + 2 \, \langle \nabla^\Sigma H_\Sigma,V^{\text{\rm tan}} \rangle \, \langle V,\nu_\Sigma \rangle. 
\end{align*}
Moreover, 
\[|\nabla^\Sigma v|^2 = \sum_{k=1}^{n-1} \langle D_{e_k} V,\nu_\Sigma \rangle^2 + 2 \sum_{k=1}^{n-1} h_\Sigma(V^{\text{\rm tan}},e_k) \, \langle D_{e_k} V,\nu_\Sigma \rangle + h_\Sigma^2(V^{\text{\rm tan}},V^{\text{\rm tan}}).\] 
Putting these facts together, we obtain 
\begin{align*} 
&\text{\rm div}_\Sigma Z - |\nabla^\Sigma v|^2 + (\text{\rm Ric}(\nu_\Sigma,\nu_\Sigma) + |h_\Sigma|^2) \, v^2 \\ 
&= \sum_{k,l=1}^{n-1} \langle D_{e_k} V,e_l \rangle \, \langle D_{e_l} V,e_k \rangle - \sum_{k,l=1}^{n-1} \langle D_{e_k} V,e_k \rangle \, \langle D_{e_l} V,e_l \rangle \\ 
&- \sum_{k=1}^{n-1} \langle D_{e_k} V,\nu_\Sigma \rangle^2 + \sum_{k=1}^{n-1} R(V,e_k,V,e_k) \\ 
&+ H_\Sigma \, h_\Sigma(V^{\text{\rm tan}},V^{\text{\rm tan}}) - H_\Sigma^2 \, \langle V,\nu_\Sigma 
\rangle^2 \\ 
&+ 2 \, H_\Sigma \sum_{k=1}^{n-1} \langle D_{e_k} V,e_k \rangle \, \langle V,\nu_\Sigma \rangle + 2 \, \langle \nabla^\Sigma H_\Sigma,V^{\text{\rm tan}} \rangle \, \langle V,\nu_\Sigma \rangle. 
\end{align*} 
Using the identity $H_\Sigma + \langle \nabla \log \rho,\nu_\Sigma \rangle = 0$, it follows that 
\begin{align*} 
&\text{\rm div}_\Sigma(\rho \, Z) - \text{\rm div}_\Sigma(\langle V^{\text{\rm tan}},\nabla^\Sigma \rho \rangle \, V^{\text{\rm tan}}) \\ 
&- \rho \, |\nabla^\Sigma v|^2 + \rho \, (\text{\rm Ric}(\nu_\Sigma,\nu_\Sigma) + |h_\Sigma|^2) \, v^2 \\ 
&- (D^2 \rho)(\nu_\Sigma,\nu_\Sigma) \, v^2 + \rho^{-1} \, \langle \nabla \rho,\nu_\Sigma \rangle^2 \, v^2 \\ 
&= \rho \sum_{k,l=1}^{n-1} \langle D_{e_k} V,e_l \rangle \, \langle D_{e_l} V,e_k \rangle - \rho \sum_{k,l=1}^{n-1} \langle D_{e_k} V,e_k \rangle \, \langle D_{e_l} V,e_l \rangle \\ 
&- \rho \sum_{k=1}^{n-1} \langle D_{e_k} V,\nu_\Sigma \rangle^2 + \rho \sum_{k=1}^{n-1} R(V,e_k,V,e_k) \\ 
&- 2 \, V(\rho) \sum_{k=1}^{n-1} \langle D_{e_k} V,e_k \rangle - (D^2 \rho)(V,V). 
\end{align*} 
Finally, a straightforward calculation gives 
\[(\mathscr{L}_V \mathscr{L}_V g)(X,Y) - (\mathscr{L}_W g)(X,Y) = 2 \, \langle D_X V,D_Y V \rangle - 2 \, R(V,X,V,Y)\] 
for all vector fields $X,Y$ on $M$. Moreover, 
\[V(V(\rho)) - W(\rho) = (D^2 \rho)(V,V).\] 
Using these identities together with the identity $H_\Sigma + \langle \nabla \log \rho,\nu_\Sigma \rangle = 0$, we obtain 
\begin{align*} 
&\frac{1}{2} \, \rho \sum_{k=1}^{n-1} (\mathscr{L}_V \mathscr{L}_V g)(e_k,e_k) + V(V(\rho)) - \text{\rm div}_\Sigma(\rho \, W^{\text{\rm tan}}) \\ 
&= \rho \sum_{k=1}^{n-1} |D_{e_k} V|^2 - \rho \sum_{k=1}^{n-1} R(V,e_k,V,e_k) + (D^2 \rho)(V,V) \\  
&= \rho \sum_{k,l=1}^{n-1} \langle D_{e_k} V,e_l \rangle^2 + \rho \sum_{k=1}^{n-1} \langle D_{e_k} V,\nu_\Sigma \rangle^2 - \rho \sum_{k=1}^{n-1} R(V,e_k,V,e_k) + (D^2 \rho)(V,V). 
\end{align*} 
Putting these facts together, we conclude that 
\begin{align*} 
&\frac{1}{2} \, \rho \sum_{k=1}^{n-1} (\mathscr{L}_V \mathscr{L}_V g)(e_k,e_k) + V(V(\rho)) - \text{\rm div}_\Sigma(\rho \, W^{\text{\rm tan}}) \\ 
&+ \text{\rm div}_\Sigma(\rho \, Z) - \text{\rm div}_\Sigma(\langle V^{\text{\rm tan}},\nabla^\Sigma \rho \rangle \, V^{\text{\rm tan}}) \\ 
&- \rho \, |\nabla^\Sigma v|^2 + \rho \, (\text{\rm Ric}(\nu_\Sigma,\nu_\Sigma) + |h_\Sigma|^2) \, v^2 \\ 
&- (D^2 \rho)(\nu_\Sigma,\nu_\Sigma) \, v^2 + \rho^{-1} \, \langle \nabla \rho,\nu_\Sigma \rangle^2 \, v^2 \\ 
&= \rho \sum_{k,l=1}^{n-1} \langle D_{e_k} V,e_l \rangle^2 + \rho \sum_{k,l=1}^{n-1} \langle D_{e_k} V,e_l \rangle \, \langle D_{e_l} V,e_k \rangle \\ 
&- \rho \sum_{k,l=1}^{n-1} \langle D_{e_k} V,e_k \rangle \, \langle D_{e_l} V,e_l \rangle - 2 \, V(\rho) \sum_{k=1}^{n-1} \langle D_{e_k} V,e_k \rangle. 
\end{align*}
From this, the assertion follows easily. This completes the proof of Proposition \ref{formula.second.variation}. \\

\section{Asymptotic behavior of solutions of linear PDEs}

\label{linear.PDE}

\begin{theorem}
\label{asymptotics.for.solutions.of.linear.PDE}
Let $N$ and $n$ be two integers with $3 \leq n \leq N$, and let $\check{\gamma}$ be a flat metric on the torus $T^{n-2}$. We define a hyperbolic metric $g_{\text{\rm hyp}}$ on $[1,\infty) \times T^{n-2}$ by $g_{\text{\rm hyp}} = r^{-2} \, dr \otimes dr + r^2 \, \check{\gamma}$. Consider a sequence $r_j \to \infty$. For each $j$, we assume that $w^{(j)}$ and $\zeta^{(j)}$ are smooth functions defined on the domain $[1,r_j] \times T^{n-2}$ such that 
\[-\text{\rm div}_{g_{\text{\rm hyp}}}(r^{N-n} \, dw^{(j)}) + (N-1) \, r^{N-n} \, w^{(j)} = \zeta^{(j)}\] 
on the domain $[2,r_j] \times T^{n-2}$ and $w^{(j)} = 0$ on the set $\{r_j\} \times T^{n-2}$. We further assume that there exists a real number $\delta \in (0,\frac{1}{2}]$ such that $|w^{(j)}| \leq r^{1-N}$, $|\zeta^{(j)}| \leq r^{1-n-\delta}$, and $|d\zeta^{(j)}|_{g_{\text{\rm hyp}}} \leq r^{1-n}$. Finally, we assume that $w$ is a function defined on $[2,\infty) \times T^{n-2}$ such that $w^{(j)} \to w$ in $C_{\text{\rm loc}}^2([2,\infty) \times T^{n-2})$. Then there exists a function $A \in C^{\frac{\delta}{10}}(T^{n-2},\check{\gamma})$ such that 
\[|w - r^{1-N} \, A| \leq C \, r^{1-N-\frac{\delta}{10}}\] 
and 
\[|\langle dr,dw \rangle_{g_{\text{\rm hyp}}} + (N-1) \, r^{2-N} \, A| \leq C \, r^{2-N-\frac{\delta}{10}}\] 
in the region $[2,\infty) \times T^{n-2}$.
\end{theorem}

The proof of Theorem \ref{asymptotics.for.solutions.of.linear.PDE} relies on several lemmata.

\begin{lemma}
\label{C1.bound.for.w}
We have $|dw^{(j)}|_{g_{\text{\rm hyp}}} \leq C \, r^{1-N}$ for $2 \leq r \leq \frac{r_j}{2}$. The constant $C$ is independent of $j$ and $r$.
\end{lemma} 

\textbf{Proof.} 
By assumption, $|w^{(j)}| \leq r^{1-N}$ and $|\zeta^{(j)}| \leq r^{1-n-\delta}$ for $1 \leq r \leq r_j$. Therefore, the assertion follows from standard interior estimates for elliptic PDE. \\

\begin{lemma}
\label{Holder.bound.for.zeta}
Let $\varphi_s: T^{n-2} \to T^{n-2}$ denote the flow generated by a parallel unit vector field on $(T^{n-2},\check{\gamma})$. For each $s \in \mathbb{R}$, we have 
\[|\zeta^{(j)} \circ \varphi_s - \zeta^{(j)}| \leq C \, r^{1-n-\frac{\delta}{4}} \, |s|^{\frac{\delta}{2}}\] 
for $2 \leq r \leq r_j$. The constant $C$ is independent of $j$, $r$, and $s$.
\end{lemma}

\textbf{Proof.} 
By assumption, $|\zeta^{(j)}| \leq r^{1-n-\delta}$ for $2 \leq r \leq r_j$. This implies 
\[|\zeta^{(j)} \circ \varphi_s - \zeta^{(j)}| \leq C \, r^{1-n-\delta}\] 
for $2 \leq r \leq r_j$ and $s \in \mathbb{R}$. On the other hand, using the estimate $|d\zeta^{(j)}|_{g_{\text{\rm hyp}}} \leq r^{1-n}$, we obtain 
\[|\zeta^{(j)} \circ \varphi_s - \zeta^{(j)}| \leq C \, r^{2-n} \, |s|\] 
for $2 \leq r \leq r_j$ and $s \in \mathbb{R}$. Putting these facts together, we obtain 
\[|\zeta^{(j)} \circ \varphi_s - \zeta^{(j)}| \leq C \, r^{1-n-\delta} \, \min \{1,r^{1+\delta} \, |s|\}\] 
for $2 \leq r \leq r_j$ and $s \in \mathbb{R}$. Thus, we conclude that 
\[|\zeta^{(j)} \circ \varphi_s - \zeta^{(j)}| \leq C \, r^{1-n-\delta} \, (r^{1+\delta} \, |s|)^{\frac{\delta}{2}}\] 
for $2 \leq r \leq r_j$ and $s \in \mathbb{R}$. Since $\delta - \frac{(1+\delta)\delta}{2} \geq \frac{\delta}{4}$, the assertion follows. This completes the proof of Lemma \ref{Holder.bound.for.zeta}. \\

\begin{lemma}
\label{Holder.bound.for.w.1}
Let $\varphi_s: T^{n-2} \to T^{n-2}$ denote the flow generated by a parallel unit vector field on $(T^{n-2},\check{\gamma})$. For each $s \in \mathbb{R}$, we have 
\[|w^{(j)} \circ \varphi_s - w^{(j)}| \leq C \, r^{1-N} \, |s|^{\frac{\delta}{2}}\] 
for $2 \leq r \leq r_j$. The constant $C$ is independent of $j$, $r$, and $s$.
\end{lemma}

\textbf{Proof.} 
Let $s$ be an arbitrary real number. Clearly, $w^{(j)} \circ \varphi_s - w^{(j)} = 0$ for $r=r_j$. Moreover,  Lemma \ref{C1.bound.for.w} implies that $|w^{(j)} \circ \varphi_s - w^{(j)}| \leq C_0 \, |s|^{\frac{\delta}{2}}$ for $r=2$. Here, $C_0$ is independent of $j$ and $s$. Using Lemma \ref{Holder.bound.for.zeta}, we obtain 
\begin{align*} 
&|-\text{\rm div}_{g_{\text{\rm hyp}}}(r^{N-n} \, d(w^{(j)} \circ \varphi_s - w^{(j)})) + (N-1) \, r^{N-n} \, (w^{(j)} \circ \varphi_s - w^{(j)})| \\ 
&= |\zeta^{(j)} \circ \varphi_s - \zeta^{(j)}| \leq C_1 \, r^{1-n-\frac{\delta}{4}} \, |s|^{\frac{\delta}{2}} 
\end{align*}
for $2 \leq r \leq r_j$. Here, $C_1$ is independent of $j$, $r$, and $s$. On the other hand, a straightforward calculation shows that 
\begin{align*} 
&-\text{\rm div}_{g_{\text{\rm hyp}}}(r^{N-n} \, d(r^{1-N} - r^{1-N-\frac{\delta}{4}})) + (N-1) \, r^{N-n} \, (r^{1-N} - r^{1-N-\frac{\delta}{4}}) \\ 
&= \frac{\delta}{4} \, \Big ( N+\frac{\delta}{4} \Big ) \, r^{1-n-\frac{\delta}{4}} 
\end{align*}
for $r \geq 2$. In the next step, we choose a large constant $C_2$ such that $C_0 < (2^{1-N}-2^{1-N-\frac{\delta}{4}}) \, C_2$ and $C_1 < \frac{\delta}{4} \, (N+\frac{\delta}{4}) \, C_2$. Using a standard comparison principle (cf. Theorem 3.3 in \cite{Gilbarg-Trudinger}), we conclude that 
\[|w^{(j)} \circ \varphi_s - w^{(j)}| \leq C_2 \, (r^{1-N} - r^{1-N-\frac{\delta}{4}}) \, |s|^{\frac{\delta}{2}}\] 
for $2 \leq r \leq r_j$. This completes the proof of Lemma \ref{Holder.bound.for.w.1}. \\

\begin{lemma}
\label{Holder.bound.for.w.2}
We have $\|w^{(j)}(r,\cdot)\|_{C^{\frac{\delta}{2}}(T^{n-2},\check{\gamma})} \leq C \, r^{1-N}$ for $2 \leq r \leq r_j$. The constant $C$ is independent of $j$ and $r$.
\end{lemma}

\textbf{Proof.} 
This follows immediately from Lemma \ref{Holder.bound.for.w.1}. \\

\begin{lemma}
\label{C2.bound}
We have $\|w^{(j)}(r,\cdot)\|_{C^2(T^{n-2},\check{\gamma})} \leq C \, r^{3-N-\frac{\delta}{10}}$ for $2 \leq r \leq \frac{r_j}{2}$. The constant $C$ is independent of $j$ and $r$.
\end{lemma} 

\textbf{Proof.} 
By assumption, $|w^{(j)}| \leq r^{1-N}$, $|\zeta^{(j)}| \leq r^{1-n}$, and $|d\zeta^{(j)}|_{g_{\text{\rm hyp}}} \leq r^{1-n}$. Using standard interior estimates for elliptic PDE, we obtain 
\[\|w^{(j)}(r,\cdot)\|_{C^{2,\frac{1}{2}}(T^{n-2},\check{\gamma})} \leq C \, r^{\frac{7}{2}-N}\] 
for $2 \leq r \leq \frac{r_j}{2}$. Moreover, Lemma \ref{Holder.bound.for.w.2} implies 
\[\|w^{(j)}(r,\cdot)\|_{C^{\frac{\delta}{2}}(T^{n-2},\check{\gamma})} \leq C \, r^{1-N}\] 
for $2 \leq r \leq \frac{r_j}{2}$. A standard interpolation inequality (cf. \cite{Lunardi}, Corollary 1.2.7 and Corollary 1.2.19) gives 
\[\|w^{(j)}(r,\cdot)\|_{C^{2,\frac{\delta^2}{50}}(T^{n-2},\check{\gamma})} \leq C \, \|w^{(j)}(r,\cdot)\|_{C^{\frac{\delta}{2}}(T^{n-2},\check{\gamma})}^{\frac{1}{5}+\frac{\delta}{25}} \, \|w^{(j)}(r,\cdot)\|_{C^{2,\frac{1}{2}}(T^{n-2},\check{\gamma})}^{\frac{4}{5}-\frac{\delta}{25}}\] 
for $2 \leq r \leq \frac{r_j}{2}$. Putting these facts together, we conclude that 
\[\|w^{(j)}(r,\cdot)\|_{C^{2,\frac{\delta^2}{50}}(T^{n-2},\check{\gamma})} \leq C \, r^{3-N-\frac{\delta}{10}}\] 
for $2 \leq r \leq \frac{r_j}{2}$. This completes the proof of Lemma \ref{C2.bound}. \\

By assumption, $w^{(j)} \to w$ in $C_{\text{\rm loc}}^2([2,\infty) \times T^{n-2})$. In view of Lemma \ref{C1.bound.for.w}, the limiting function $w$ satisfies $|w| \leq r^{1-N}$ and $|dw|_{g_{\text{\rm hyp}}} \leq C \, r^{1-N}$ for $r \geq 2$. This implies 
\begin{equation} 
\label{radial.derivative.of.w.1} 
\Big | \frac{\partial}{\partial r} w + (N-1) \, r^{-1} \, w \Big | \leq C \, r^{-N} 
\end{equation} 
for $r \geq 2$. Moreover, the function $w$ satisfies 
\begin{equation} 
\label{pde.for.w.1} 
|\text{\rm div}_{g_{\text{\rm hyp}}}(r^{N-n} \, dw) - (N-1) \, r^{N-n} \, w| \leq r^{1-n-\delta} 
\end{equation} 
for $r \geq 2$. The inequality (\ref{pde.for.w.1}) can be rewritten as 
\begin{equation} 
\label{pde.for.w.2}
\Big | r^2 \, \frac{\partial^2}{\partial r^2} w + (N-1) \, r \, \frac{\partial}{\partial r} w + r^{-2} \, \Delta_{\check{\gamma}} w - (N-1) \, w \Big | \leq r^{1-N-\delta} 
\end{equation} 
for $r \geq 2$. Using Lemma \ref{C2.bound}, we obtain $|\Delta_{\check{\gamma}} w| \leq C \, r^{3-N-\frac{\delta}{10}}$ for all $r \geq 2$. Using this estimate together with (\ref{pde.for.w.2}), we conclude that 
\begin{equation} 
\label{pde.for.w.3}
\Big | r^2 \, \frac{\partial^2}{\partial r^2} w + (N-1) \, r \, \frac{\partial}{\partial r} w - (N-1) \, w \Big | \leq C \, r^{1-N-\frac{\delta}{10}} 
\end{equation}
for $r \geq 2$. The inequality (\ref{pde.for.w.3}) can be rewritten as 
\begin{equation} 
\label{pde.for.w.4}
\Big | \frac{\partial}{\partial r} \Big ( \frac{\partial}{\partial r} w + (N-1) \, r^{-1} \, w \Big ) \Big | \leq C \, r^{-N-1-\frac{\delta}{10}} 
\end{equation} 
for $r \geq 2$. In the next step, we integrate the inequality (\ref{pde.for.w.4}) along radial curves. Using (\ref{radial.derivative.of.w.1}), we conclude that 
\begin{equation} 
\label{radial.derivative.of.w.2} 
\Big | \frac{\partial}{\partial r} w + (N-1) \, r^{-1} \, w \Big | \leq C \, r^{-N-\frac{\delta}{10}} 
\end{equation}  
for $r \geq 2$. The inequality (\ref{radial.derivative.of.w.2}) can be rewritten as 
\begin{equation} 
\label{radial.derivative.of.w.3}
\Big | \frac{\partial}{\partial r}(r^{N-1} \, w) \Big | \leq C \, r^{-1-\frac{\delta}{10}} 
\end{equation} 
for $r \geq 2$. It follows from (\ref{radial.derivative.of.w.3}) that the functions $r^{N-1} \, w(r,\cdot) \in C^0(T^{n-2},\check{\gamma})$ converge uniformly to a function $A \in C^0(T^{n-2},\check{\gamma})$ as $r \to \infty$. Moreover, 
\begin{equation} 
\label{asymptotic.expansion.for.w}
|r^{N-1} \, w - A| \leq C \, r^{-\frac{\delta}{10}} 
\end{equation}
for $r \geq 2$. Combining (\ref{radial.derivative.of.w.2}) and (\ref{asymptotic.expansion.for.w}), we obtain 
\begin{equation} 
\label{asymptotic.expansion.for.radial.derivative.of.w}
\Big | r^N \, \frac{\partial}{\partial r} w + (N-1) \, A \Big | \leq C \, r^{-\frac{\delta}{10}} 
\end{equation} 
for $r \geq 2$. Finally, Lemma \ref{Holder.bound.for.w.2} implies that $\|r^{N-1} \, w(r,\cdot)\|_{C^{\frac{\delta}{2}}(T^{n-1},\check{\gamma})} \leq C$ for all $r \geq 2$. Consequently, the function $A$ belongs to the H\"older space $C^{\frac{\delta}{2}}(T^{n-2},\check{\gamma})$. This completes the proof of Theorem \ref{asymptotics.for.solutions.of.linear.PDE}.

\end{document}